\documentclass{amsart}

\usepackage{geometry}
\usepackage{amsmath}
\usepackage{amssymb}
\usepackage{amsthm}
\usepackage{color}
\usepackage[shortlabels]{enumitem}
\usepackage{hyperref}

\usepackage{hyperref}

\newtheorem{prop}{Proposition} 
\newtheorem{cor}{Corollary}
\newtheorem{lemma}{Lemma}
\newtheorem{thm}{Theorem}
\newtheorem{defn}{Definition}

\numberwithin{equation}{section}

\usepackage[dvipsnames]{xcolor}

\begin{document}

\begin{abstract}
We prove that the primes below $x$ are, on average, equidistributed in arithmetic progressions to smooth moduli of size up to  $x^{1/2+1/40-\epsilon}$. The exponent of distribution $\frac{1}{2} + \frac{1}{40}$ improves on a result of Polymath~\cite{Polymath:2014:EDZ}, who had previously obtained the exponent  $\frac{1}{2} + \frac{7}{300}$. 
As a consequence, we improve results on intervals of bounded length which contain many primes, showing that 
\begin{align*}
\liminf_{n \rightarrow \infty} (p_{n+m}-p_n)  = O( \exp (3.8075 m)).
\end{align*}
The main new ingredient of our proof is a modification of the $q$-van der Corput process. It allows us to exploit additional averaging for the exponential sums which appear in the  Type I estimates of~\cite{Polymath:2014:EDZ}.
\end{abstract}

\title{On primes in arithmetic progressions and bounded gaps between many primes}
\date{}
\author{Julia Stadlmann}
\maketitle

\vspace{-5mm}

\section{Introduction}\label{sec:intro}

In this paper we study equidistribution estimates for primes in arithmetic progressions. Given some $\theta >0$ and  sets $\mathcal{Q}(x) \subseteq \mathbb{N}$, we consider the following statement: For $x>1$,  $a \in \mathbb{N}$,  $A>0$ and $\varepsilon>0$, 
\begin{align} \label{equidistribution}
\sum_{\substack{ q \leq x^{ \theta -\varepsilon } \\ q  \in \mathcal{Q}(x) \\ (q,a)=1 }} \Bigg| \sum_{\substack{ p \leq x \\ p \equiv a (q) }} 1 - \dfrac{1}{\phi(q)}\sum_{\substack{ p \leq x \\ (p,q)=1 }} 1  \,\,\,\, \Bigg| \ll_{A, \varepsilon} \dfrac{x}{\log(x)^A}.
\end{align} 
Bombieri~\cite{Bombieri:1965:OLS} and Vinogradov~\cite{Vinogradov:1965:DHD} showed in 1965 that (\ref{equidistribution}) holds with the exponent $\theta = \frac{1}{2}$ when $\mathcal{Q}(x) =\mathbb{N}$ for all $x$.  The Elliott–Halberstam conjecture suggests that (\ref{equidistribution}) also holds with $\theta = 1$,  but to date the statement (\ref{equidistribution})  remains unproven for any $\theta > \frac{1}{2}$ when $\mathcal{Q}(x) =\mathbb{N}$. 

\vspace{3mm}
In the 1980s, Bombieri, Fouvry, Friedlander and Iwaniec obtained better exponents for   variants of (\ref{equidistribution})  in which the absolute value signs had been replaced by suitable coefficients. For instance, Fouvry and Iwaniec~\cite{Fouvry:1983:PAP} obtained the exponent of distribution $\frac{9}{17}$ for well-factorable coefficients and Bombieri, Friedlander and Iwaniec~\cite{Bombieri:2986:PAP} improved this exponent to $\frac{4}{7}$. Additionally, Bombieri, Friedlander and Iwaniec~\cite{Bombieri:2986:PAP} obtained the exponent  $\frac{29}{56}$  for a less restrictive set of coefficients when $Q(x) = \{q_1q_2:q_1 \leq x^{\theta_1}, q_2 \leq x^{\theta_2} \}$ with $\theta_1<\frac{1}{3}$, $\theta_2<\frac{1}{5}$ and $5\theta_1+2\theta_2<2$. More recently, Maynard~\cite{Maynard:2022:PA2} obtained the exponent of distribution  $\frac{3}{5}$ for triply well-factorable coefficients. This result  was improved to $\frac{66}{107}$ by Lichtman~\cite{Lichtman:2023:TWF} and to $\frac{5}{8}$ by
Pascadi~\cite{Pascadi:2024:TWF}.
 Maynard~\cite{Maynard:2022:PA1} also proved that (\ref{equidistribution}) holds with $\theta = \frac{11}{21}$ and $\mathcal{Q}(x) = \{q_1q_2: q_1 \leq x^{1/21}, q_2 \leq x^{10/21-\varepsilon}\}$, provided that the implied constant is allowed to depend on $a$. Further, Maynard~\cite{Maynard:2022:PA3}  proved that (\ref{equidistribution}) holds with $\theta = \frac{1}{2} + \delta$ for $\delta \in (0, \frac{1}{1000})$  when $\mathcal{Q}(x) = \{q_1q_2q_3: q_1 \leq Q_1, q_2 \leq Q_2, q_3 \leq Q_3\}$ with  $x^{40\delta} < Q_2 < x^{1/20-7\delta}$ and $\frac{x^{1/10+12\delta}}{Q_2} <Q_3< \frac{x^{1/10-4\delta}}{Q_2^{3/5}}$.

\vspace{3mm} In this paper, we are interested in the case $\mathcal{Q}(x) = \{ q \in \mathbb{N}: q \mid P(x^\delta)\}$ where $P(x^\delta) = \prod_{p < x^\delta} p$. This is the set of squarefree, $x^\delta$-smooth integers. Equidistribution estimates for primes in APs to smooth moduli played a key role in Zhang's proof of bounded gaps between primes~\cite{Zhang:2015:BGP}. Zhang obtained the exponent of distribution $\frac{1}{2}+\frac{1}{584}$ for a variant of (\ref{equidistribution}). Polymath~\cite{Polymath:2014:EDZ} improved on these equidistribution estimates, showing that (\ref{equidistribution}) holds with $\theta = \frac{1}{2} + \frac{7}{300}$ when $\mathcal{Q}(x) = \{ q \in \mathbb{N}: q \mid P(x^\delta)\}$. We prove the following:

\begin{thm} \label{thm:equidistribution}
Let $P(z) = \prod_{p<z} p$. Let $\varepsilon>0$. There exists $\delta >0$, dependent on $\varepsilon$, such that the following is true: For $x>1$,  $a \in \mathbb{Z}$  and $A>0$, we have
\begin{align*}
\sum_{\substack{ q \leq x^{1/2 +1/40 -\varepsilon } \\ q \mid P(x^\delta) \\ (q,a)=1}} \Bigg| \sum_{\substack{ p \leq x \\ p \equiv a (q) }} 1 - \dfrac{1}{\phi(q)}\sum_{\substack{ p \leq x \\ (p,q)=1 }} 1  \,\,\,\, \Bigg| \ll_{A, \varepsilon} \dfrac{x}{\log(x)^A}.
\end{align*} 
\end{thm}
The key ingredient of our proof  is an improved exponential sum bound for what was the limiting case in previous works. By taking advantage of the precise shape of the exponential phases appearing in this problem, we are able to exploit the summation over a variable which was previously treated trivially. A sketch of this argument can be found in Section~\ref{ssec:sketch1}.

\vspace{3mm} Since equidistribution estimates of the form described in Theorem~\ref{thm:equidistribution} are closely linked to bounds on the infimum limit of gaps between primes, we also  study the quantity
\begin{align*}
H_m=\liminf_{n \rightarrow \infty} (p_{n+m}-p_n).
\end{align*}
For large $m$, a first upper bound on this quantity was given by Maynard~\cite{Maynard:2015:SGP}, who showed that $H_m \ll m^3 \exp(4m)$; and a similar bound was also obtained by Tao around the same time. The methods of Maynard and Tao were further optimized in the Polymath paper~\cite{Polymath:2014:VSS}, leading to the bound $H_m \ll \exp((4-\frac{28}{157})m)$. Note here that $4-\frac{28}{157} \approx 3.822$.  Baker and Irving~\cite{Baker:2017:BIP} combined the techniques of the previous papers with Harman's sieve, which gave  $H_m \ll \exp(3.815 m)$. Using the method of Polymath together with Theorem~\ref{thm:equidistribution}, we are immediately able to deduce that $H_m \ll \exp(3.81 m)$, improving on the result of Baker and Irving. However, if we additionally combine our new equidistribution estimates with Harman's sieve, like in~\cite{Baker:2017:BIP}, we arrive at the following bound:

\begin{thm} \label{thm:infinumboundedgaps}
The quantity $H_m=\liminf_{n \rightarrow \infty} (p_{n+m}-p_n)$ satisfies $ H_m \ll \exp (3.8075 m)$. 
\end{thm}

Our new equidistribution estimates can also be used  to improve bounds on $H_m$ for fixed small $m \geq 2$.  Denoting by $H(k)$ the diameter of the narrowest admissible $k$-tuple,  Polymath~\cite{Polymath:2014:VSS} showed that 
\begin{align*}
&H_2 \leq H(35 410) \leq 398 130,\\
&H_3 \leq H(1649821) \leq  24797814,\\
&H_4 \leq H(75845707) \leq  1431556072,\\
&H_5 \leq H(3473955908) \leq 80550202480.
\end{align*}
For given $m \in \{2,3,4,5\}$, we are now able to decrease the size of $k$ for which it is known that $H_m \leq H(k)$. Sutherland computed corresponding bounds on $H(k)$, giving us the following improvement:

\begin{cor}\label{cor:H2}
Write  $H_m=\liminf_{n \rightarrow \infty} (p_{n+m}-p_n)$. Denote by $H(k)$ the diameter of the narrowest admissible $k$-tuple. Then the following is true:
\begin{align*}
&H_2 \leq H(35265) \leq 396504, \\
&H_3 \leq H(1624545) \leq  24407016,\\
&H_4 \leq H(73807570) \leq  1391051532,\\
&H_5 \leq H(3340375663) \leq 77510685234.
\end{align*} 
\end{cor}

\subsection{A sketch of the proof of Theorem 1} \label{ssec:sketch1} We now outline the main arguments of the proof of Theorem~\ref{thm:equidistribution}.  Many details are suppressed for clarity.

\vspace{3mm}
 By the Heath-Brown identity, in order to prove that the primes have exponent of distribution $1/2+2\omega$ to smooth moduli, it suffices to show the following: If $\alpha_1, \dots, \alpha_j : \mathbb{N} \rightarrow \mathbb{C}$  with $\alpha_i$ supported on $[N_i,2N_i]$,  $\prod_{i=1}^j N_i \asymp x$ and ($\alpha_i(n)=1$ when $N_i \gg x^\delta$), then 
\begin{align} \label{equ:hb}
\sum_{\substack{ q \leq x^{1/2+2\omega -\varepsilon} \\ q \mid P(x^\delta) \\ (q,a)=1
 }} \Bigg| \sum_{ \substack{ n \equiv a (q) \\ n \in [x,2x]}} (\alpha_1 \star \dots \star \alpha_j)(n) -\dfrac{1}{\phi(q) } \sum_{\substack{ (n,q)=1\\ n \in [x,2x]}}  (\alpha_1 \star \dots \star \alpha_j)(n)  \Bigg|  \ll_{A,\delta} \dfrac{x}{\log(x)^A}.
\end{align}
The current techniques for establishing estimates of the form~(\ref{equ:hb}) are highly dependent on the values of $N_1, \dots, N_j$. For instance, if we use the Type I/II/III estimates of Polymath~\cite{Polymath:2014:EDZ} and take $\omega = 7/600$, then (\ref{equ:hb}) holds if there are distinct $i_1, i_2, i_3 \in \{1,\dots,j\}$ with $N_{i_1}N_{i_2}N_{i_3} \geq x^{0.9}$ and $N_{i_j} \gg x^\delta$ and (\ref{equ:hb}) also holds if there exists $I \subseteq \{1, \dots, j\}$ with $\prod_{i \in I} N_i \in [x^{0.4},x^{0.6}]$. Every choice of $\alpha_1 \star \dots \star \alpha_j$ under consideration satisfies one of these two properties, and so  Polymath~\cite{Polymath:2014:EDZ} obtained the exponent of distribution $\frac{1}{2}+\frac{7}{300}$. However,  if we increase $\omega$ just a little further, then the results of Polymath do not provide sufficiently good equidistribution  estimates whenever $N_1, \dots, N_5$ are all close to $x^{0.2}$.

\vspace{3mm} To establish a larger exponent of distribution, our goal is  to prove that (\ref{equ:hb}) still holds for $N_1, \dots, N_j$ with  $\prod_{i \in I} N_i \approx x^{0.4}$ if $\omega$ is a bit larger than $7/600$. (We thus improve the Type I equidistribution estimates of Polymath~\cite{Polymath:2014:EDZ}.)  
In~\cite{Polymath:2014:EDZ} it was shown that in order to obtain such a result, it is  enough to prove that
\begin{align} \label{equ:linnik}
\sum_{\substack{d \asymp \Delta}} \Bigg| \sum_{n \asymp N} \sum_{\substack{y \asymp Y \\ (y,m)=1}} e_{md}\left(\dfrac{By}{n+Cd} \right) \Bigg| \ll \dfrac{\Delta N}{x^\epsilon}
\end{align} 
for some choice of $\Delta$ with $1 \ll \Delta \ll N$ and some $\epsilon>0$   and for $ N \approx x^{0.4}$, $Y \approx x^{8\omega}$, $m \approx x^{1+8\omega}/(N\Delta)$ and $B, C \in \mathbb{N}$ with $(B,md)=1$. Here $e_{q}(\frac{a}{b}):= \exp(\frac{2\pi i n}{q})$ for $n \in \mathbb{Z}$  with $n \equiv a/b$ mod $q$, see also  Section~\ref{sssec:exponentials}. (Note: Inequality (\ref{equ:linnik}) is primarily a consequence of Linnik's dispersion method and completion of sums.) \vspace{3mm}

We now demonstrate two different approaches to  bounding (\ref{equ:linnik}). In both cases, the sum over $d$ in (\ref{equ:linnik}) will be incorporated into the $q$-van der Corput process to reduce the modulus to $m$.  The optimal choice of $\Delta$ balances the resulting diagonal and off-diagonal terms: larger choices of $\Delta$ correspond to better bounds for the off-diagonal terms, but worse bounds for the diagonal terms. 

\subsubsection{The q-van der Corput method - Polymath} We first recap the arguments of Polymath~\cite{Polymath:2014:EDZ}, paying particular attention to the diagonal terms. Using the Cauchy-Schwarz inequality, we have 
\begin{align*}
\sum_{\substack{d \asymp \Delta}} \Bigg| \sum_{n \asymp N} \sum_{\substack{y \asymp Y \\ (y,m)=1}} e_{md}\left(\dfrac{By}{n+Cd} \right) \Bigg| 
 &\ll \sum_{\substack{d \asymp \Delta}} \sum_{c(d)}\sum_{\substack{y \asymp Y \\ (y,m)=1}} \Bigg| \sum_{\substack{n \asymp N \\ n \equiv c(d)}}  e_{m}\left(\dfrac{By}{d(n+Cd)} \right) \Bigg| \\
 &\ll (\Delta^2 Y)^{1/2} \Bigg( \sum_{\substack{d \asymp \Delta}}\sum_{c(d)}\sum_{\substack{y \asymp Y \\ (y,m)=1}}  \sum_{\substack{n, \widetilde{n} \asymp N \\ n, \widetilde{n} \equiv c(d)}}  e_{m}\left(\dfrac{By (\widetilde{n}-n)}{d(n+Cd)(\widetilde{n}+Cd)}  \right)  \Bigg)^{1/2}
 \\
 &\ll (\Delta^2 Y)^{1/2} \Bigg(  \sum_{\substack{y \asymp Y \\ (y,m)=1}} \sum_{k \ll \frac{N}{\Delta}}\sum_{\substack{d \asymp \Delta}} \sum_{\substack{n \asymp N }} e_{m}\left(\dfrac{By k}{n(n+kd)}  \right)  \Bigg)^{1/2}
  \\
 &\ll  \Delta^{3/2} N^{1/2}Y + \Delta Y \!\!\sup_{\substack{y \asymp Y \\ (y,m)=1}}\!\! \Bigg(  \sum_{1 \leq k \ll \frac{N}{\Delta}} \sum_{\substack{d \asymp \Delta}} \sum_{\substack{n \asymp N }} e_{m}\left(\dfrac{By k}{n(n+kd)}  \right)  \Bigg)^{1/2}.
\end{align*} 
Here $ \Delta^{3/2} N^{1/2}Y$ is the contribution of the diagonal terms (which correspond to $n = \widetilde{n}$). We seek an upper bound of $\Delta N /x^\epsilon$, so this forces $ \Delta <  N/Y^2$.  For the off-diagonal terms, Polymath~\cite{Polymath:2014:EDZ} used the $q$-van der Corput process and  Deligne's work on the Riemann Hypothesis over finite fields to prove that 
\begin{align} \label{equ:introcorput}
\sum_{1 \leq k \ll \frac{N}{\Delta}} \sum_{\substack{d \asymp \Delta}} \sum_{\substack{n \asymp N }} e_{m}\left(\dfrac{By k}{n(n+kd)}  \right) \ll \dfrac{Nm^{2/3}}{\Delta^{1/2}}.
\end{align}
Recall that  $m \approx x^{1+8\omega}/(N \Delta)$, $Y \approx x^{8\omega}$ and $\Delta \approx N/Y^2$. To get an upper bound of $\Delta N/ x^\epsilon$ on the LHS of  (\ref{equ:linnik}), we  require $(Nm^{2/3}/\Delta^{1/2})^{1/2} \ll N/(x^\epsilon Y)$ and hence $N \gg x^{4/17+240\omega/17+12\epsilon/17}$. But recall that we are interested in the case $N \approx x^{0.4}$. This now forces $\omega \leq 7/600$.

\vspace{3mm} Larger $\Delta$ gave better bounds on the LHS of (\ref{equ:introcorput}), but the maximal possible size of $\Delta$ was determined by the diagonal terms, which contributed $ \Delta^{3/2}N^{1/2}Y$.  So  to get an improvement, the key idea is to change the shape of the diagonal terms, reducing their contribution. Observe that the argument above made no use of the summation over $y$. This summation can be exploited to get smaller diagonal terms.

\subsubsection{The q-van der Corput method - Variant} \label{sssec:variant} Now we introduce our new argument.  If we set $y/n \equiv c$ mod $d$ and apply the Chinese Remainder Theorem to deduce $e_{md}(\cdot) = e_d(\frac{\cdot}{m}) e_m(\frac{\cdot}{d})$, then 
\begin{align*}
e_{md}\left( \frac{By}{n+Cd} \right) = e_d\left( \frac{By}{m(n+Cd)} \right)e_m\left( \frac{By}{d(n+Cd)} \right)= e_d\left( \frac{Bc}{m} \right)e_m\left( \frac{By}{d(n+Cd)} \right)
\end{align*}
and the $d$-factor no longer depends on $n$ and $y$. This allows us to apply the $q$-van der Corput process with both $n$ and $y$ on the inside. More precisely, we have 
\begingroup
\allowdisplaybreaks
\begin{align*}
\sum_{\substack{d \asymp \Delta}} \Bigg| \sum_{n \asymp N} \sum_{\substack{y \asymp Y }} e_{md}\left(\dfrac{By}{n+Cd} \right) \Bigg| &\ll \sum_{\substack{d \asymp \Delta}} \sum_{c(d)} \Bigg| \sum_{\substack{n \asymp N }} \sum_{\substack{y \asymp Y \\ y/n \equiv c(d)}} e_{m}\left(\dfrac{By}{d(n+Cd)} \right) \Bigg| \\
 &\ll \Delta \Bigg( \sum_{\substack{d \asymp \Delta}}\sum_{c(d)} \sum_{\substack{n, \widetilde{n}\asymp N  \\ y, \widetilde{y} \asymp Y \\ y/n \equiv \widetilde{y}/\widetilde{n} \equiv c(d)}}  e_{m}\left(\dfrac{By }{d(n+Cd)} -  \dfrac{B \widetilde{y} }{d(\widetilde{n}+Cd)} \right)  \Bigg)^{1/2}
 \\
 &\ll \Delta \Bigg( \sum_{\substack{ y, \widetilde{y} \asymp Y }} \sum_{\substack{d \asymp \Delta}} \sum_{\substack{n, \widetilde{n}\asymp N  \\  y\widetilde{n} \equiv \widetilde{y} n (d)}}  e_{m}\left(\dfrac{B(y\widetilde{n} - \widetilde{y}n+C(y-\widetilde{y})d) }{d(n+Cd)(\widetilde{n}+Cd)} \right)  \Bigg)^{1/2}. 
\end{align*} 
\endgroup
To prove that (\ref{equ:linnik}) holds, it thus suffices to show that for all $y, \widetilde{y} \asymp Y$, 
\begin{align} \label{equ:intronewcorput}
   \sum_{\substack{d \asymp \Delta}} \sum_{\substack{n, \widetilde{n}\asymp N  \\  y\widetilde{n} \equiv \widetilde{y} n (d)}}  e_{m}\left(\dfrac{B(y\widetilde{n} - \widetilde{y}n+C(y-\widetilde{y})d) }{d(n+Cd)(\widetilde{n}+Cd)}   \right) \ll \dfrac{ N^2}{x^{2\epsilon} Y^2}.
\end{align}
Observe that $y\widetilde{n} \equiv \widetilde{y} n (d)$ if and only if $\widetilde{n} = (n \widetilde{y} + kd)/y$ for some $k \ll Y N/\Delta$ with $y \mid n \widetilde{y} +kd$.  Then 
\begin{align*}
e_{m}\left(\dfrac{B(y\widetilde{n} - \widetilde{y}n+C(y-\widetilde{y})d) }{d(n+Cd)(\widetilde{n}+Cd)}   \right) = e_{m}\left(\dfrac{B(k+C(y-\widetilde{y})) }{(n+Cd)(\frac{\widetilde{y}}{y}n + \frac{kd}{y}+Cd)}   \right).
\end{align*}
The diagonal terms of  (\ref{equ:intronewcorput}) now correspond to choices of $n$, $k$ and $d$ for which $k+C(y - \widetilde{y}) = 0$ and $y \mid n \widetilde{y} + kd$. Fix some $d \asymp \Delta$. The diagonal terms have $k = -C(y - \widetilde{y})$ with $y \mid n \widetilde{y} - d  C(y - \widetilde{y})$.  Assuming $(y, \widetilde{y})=1$, only $O(N/Y)$ choices of $n$ satisfy this divisibility condition. Overall, the diagonal terms thus contribute $O(\Delta N/Y)$ to (\ref{equ:intronewcorput}). This forces $\Delta < N/Y$. (The value of $\Delta$ is thus allowed to be larger than in the work of Polymath~\cite{Polymath:2014:EDZ}, where $\Delta < N/Y^2$.)

\vspace{3mm} For the off-diagonal terms, we first notice that the condition $ y \mid n \widetilde{y}+kd$ gives $n = y n_1 + F_k d$ for some $n_1 \in \mathbb{N}$ and some constant $F_k$ (dependent on $k$). We further write $\frac{\widetilde{y}}{y}n + \frac{kd}{y} \equiv \widetilde{y} n_1 + G_k d$ mod $m$ and have the relationship $\frac{F_k}{y} - \frac{G_k}{\widetilde{y}} \equiv - \frac{kd}{y \widetilde{y}}$ mod $m$.  Substituting and using  Deligne's bounds for the summation over $d$ and $n$, we obtain the following estimate: 
\begin{align*} 
  &\sum_{\substack{ k \ll Y N/\Delta \\ k+C(y- \widetilde{y}) \neq 0}} \sum_{\substack{d \asymp \Delta}} \sum_{\substack{n\asymp N  \\  y \mid n \widetilde{y}+kd}}  e_{m}\left(\dfrac{B(k+C(y-\widetilde{y})) }{(n+Cd)(\frac{\widetilde{y}}{y}n + \frac{kd}{y}+Cd)}   \right) \\
  =  &\sum_{\substack{ k \ll Y N/\Delta \\ k+C(y- \widetilde{y}) \neq 0}} \sum_{\substack{d \asymp \Delta}} \sum_{\substack{n_1 \asymp N }}  e_{m}\left(\dfrac{B(\frac{k}{y \widetilde{y}}+C(\frac{1}{\widetilde{y}}-\frac{1}{y})) }{( n_1 + (\frac{F_k}{y} +\frac{C}{y})d)(n_1 + (\frac{G_k}{\widetilde{y}}+\frac{C}{\widetilde{y}})d)}   \right) \ll \dfrac{ Y N m   }{\Delta }. \nonumber
\end{align*}
Then (\ref{equ:intronewcorput}) and hence (\ref{equ:linnik}) hold provided that $m \ll  N \Delta/(x^{2\epsilon}Y^3)$. Recalling that $\Delta \approx N/Y$ and $m \approx x^{1+8\omega}/(N\Delta)$, we require $N \gg x^{1/4+12\omega}$. At $\omega = 1/80$, this condition is satisfied by $N \approx x^{0.4}$, concluding the sketch of our proof. 

\vspace{3mm} In this argument, we exploited the simple shape of the quotient $(\frac{B y }{n +Cd})$ inside the exponential sums of (\ref{equ:linnik}). This is what gave us simple expressions for the diagonal terms and what allowed us to easily deal with arising modularity conditions via substitution of linear expressions.

\subsection{Notation} Below is a list of notational conventions and important definitions.

\subsubsection{Basics} 

In this paper, lower case Roman letters denote integers and upper case Roman letters denote real numbers (unless specified otherwise). In particular, all sums and products are over the integers, with the exception of sums and products over $p$ or $p_i$, which are over the primes.

\vspace{3mm}

We let $h \sim H$  denote the condition $h \in (H, 2H] \cup [2H,H)$. Here $H$ may be positive or negative. We also use Vinogradov  notation ($\ll$ and $\gg$) and write $X \asymp Y$ if $Y \ll X \ll Y$. If we write $X \ll_A Y$, then the implied constant may depend on $A$.

\vspace{3mm} For $\delta  > 0$, we set $P(x^\delta) = \prod_{p < x^\delta} p$. Divisors of $P(x^\delta)$ are squarefree and $x^\delta$-smooth.

\vspace{3mm}
As usually, Dirichlet convolutions are defined as follows: If $\alpha: \mathbb{N} \rightarrow \mathbb{C}$ and $\beta: \mathbb{N} \rightarrow \mathbb{C}$, then
\begin{align*}
(\alpha \star \beta)(n) = \sum_{d \mid n} \alpha(d) \beta\left(\dfrac{n}{d}\right).
\end{align*}

The function $\phi(n)$ is the Euler totient function, $\Lambda(n)$ is the von Mangoldt function and $\mu(n)$ is the M\"{o}bius function. If $Q(n)$ is a statement about integers, then  $1_{Q(n)}$ is the indicator function of that property. 

\vspace{3mm}  In all our proofs, we assume without loss of generality that   $x$ is very large.

\subsubsection{Exponentials}\label{sssec:exponentials} For $q \in \mathbb{N}$ and $n \in \mathbb{Z}$, we define $e_q(n) = \exp(2\pi i n/q)$. For $a, b \in \mathbb{Z}$, we set 
\begin{align*}
e_q\left(\dfrac{a}{b}\right) = \begin{cases} \,\,
e_q\big( \overline{a/b} \big) \quad \mbox{ if } a/b  \mbox{ is well-defined mod } q, \\ \,\, 0 \qquad\qquad\! \mbox{ otherwise,}
\end{cases}
\end{align*}
where $\overline{a/b}$ represents an integer $n$ which satisfies $n \equiv a/b$ mod $q$. \vspace{3mm}

One note of caution: generally, $e_q(\frac{a}{b}) \neq \exp(\frac{2\pi i a}{bq})$. For instance, $e_3(\frac{1}{2}) = \exp(\frac{4\pi i}{3})$, while $e_3(\frac{1}{3})=0$.

\vspace{3mm}
By the Chinese Remainder Theorem, we have the following  identity for any $n \in \mathbb{Z}$: If $(q_1,q_2)=1$, then
\begin{align*}
e_{q_1q_2}\left(n\right)=e_{q_1}\left(\dfrac{n}{q_2}\right)e_{q_2}\left(\dfrac{n}{q_1}\right).
\end{align*}
More generally, we will often use that for $a,b \in \mathbb{Z}$ and pairwise coprime natural numbers $q_1, \dots, q_k$,
\begin{align} \label{CRT}
e_{q_1\dots q_k}\left(\dfrac{a}{b}\right)=\prod_{i=1}^k e_{q_i}\left(\dfrac{a}{b\prod_{j \neq i} q_j}\right).
\end{align}

\subsubsection{Sequences}

For a given collection of functions from $\mathbb{N}$ to $\mathbb{C}$, we now formally define the  exponent of distribution to smooth moduli:

\begin{defn}[Exponent of distribution]
We say that a function $f: \mathbb{N} \times (1,\infty) \rightarrow \mathbb{C}$ has exponent of distribution $\theta$ to smooth moduli if the following holds: For any $\varepsilon >0$ there exists $\delta>0$ such that for any $x>1$,  $a \in \mathbb{Z}$ and  $A>0$, we have
\begin{align*}
\sum_{\substack{q \leq x^{\theta -\varepsilon } \\ q \mid P(x^\delta) \\ (q,a)=1}} \Bigg| \sum_{\substack{n \in [x,2x] \\ n \equiv a (q)}} f(n;x) - \dfrac{1}{\phi(q)} \sum_{\substack{n \in [x,2x] \\(n,q)=1 }}f(n;x) \Bigg| \ll_{A,\varepsilon} \dfrac{x}{\log(x)^A}.
\end{align*}
\end{defn}
(Note: This is essentially the same as  Definition~1 of~\cite{Baker:2017:BIP}, we have just rephrased the statement slightly to make it clearer that the functions under consideration may depend on $x$, whereas the implied constant must be independent of $x$.)

\vspace{3mm}

We will frequently work with coefficient sequences, which are defined as follows: 

\begin{defn}[Coefficient sequences]\label{def:sequence}  A  coefficient sequence  is  a function $\alpha: \mathbb{N} \times (1, \infty) \rightarrow \mathbb{R}$ which satisfies $|\alpha(n;x)| \ll \tau(n)^{O(1)}\log(x)^{O(1)}$. 
\begin{enumerate}[{\rm (i)}]
\item $\alpha$ is said to be located at scale $N(x)$ if there are some constants $1 \ll c \ll C \ll 1$ such that   $\alpha(\,\cdot\,;x)$ is supported on $[cN(x),CN(x)]$ for every $x>1$.
\item If $\alpha$ is located at scale $N(x)$, it is said to have the Siegel-Walfisz property if 
\begin{align*}
\Bigg| \sum_{\substack{n=a (q) \\ (n,r)=1 }} \alpha(n;x)  -\dfrac{1}{\phi(q)}\sum_{\substack{ (n,qr)=1 }} \alpha(n;x) \Bigg| \ll_A \tau(qr)^{O(1)}\dfrac{ N(x) }{\log(x)^{A}}
\end{align*}
for any $x>1$, any $q, r \geq 1$, any $A >1$, and any  residue class $a$ mod $q$ with $(a,q)=1$. 
\item $\alpha$ is said to be shifted smooth at scale $N(x)$ if there are some constants $1\ll c \ll C \ll 1$ so that for every $x>1$, there exists a constant $x_0 \in \mathbb{R}$ and a  smooth function $\psi:\mathbb{R} \rightarrow \mathbb{C}$ supported on $[c,C]$, with $|\psi^{(j)}(t)| \ll_j \log(x)^{O_j(1)}$ for all $t \in \mathbb{R}$, such that 
\begin{align*}
\alpha(n;x) = \psi\left(\dfrac{n-x_0}{N}\right).
\end{align*}
If we can take $x_0=0$ for all $x$, $\alpha$ is said to be smooth at scale $N(x)$. 
\end{enumerate}
\end{defn}

(Note 1: This definition is essentially the same as  Definition~2.5 of   Polymath~\cite{Polymath:2014:EDZ}, again we have just rephrased the statement slightly to highlight that coefficient sequences may depend on $x$, while implied constants must not depend on $x$.) \vspace{2mm}

(Note 2: Later on, we will adopt Polymath's convention of suppressing the dependence of $\alpha$ on $x$ in our notation. This means that  we will typically denote sequence entries by $\alpha(n)$ (rather than $\alpha(n;x)$) and will say that $\alpha$ is located at scale $N$ (rather than $N(x)$), since proofs would otherwise look quite cluttered.) 
\vspace{2mm}

(Note 3: In this paper, we will often consider a set of coefficient sequences at a given scale. In each such case, we  assume that the coefficient sequences under consideration all satisfy Definition~\ref{def:sequence} with the same choice of implied constants, and upper bounds we prove  will depend on this choice of implied constants.)

\subsubsection{Notation index}

Finally, for the reader's convenience, we provide a notation index which gives an overview over the properties of important variables used in the proofs of Section~\ref{sec:equidistribution}. The index is roughly sorted by type of  expression and chronological appearance in our proofs, with related variables placed closely together. (Of course all these expressions will also be defined more formally at a later stage.) 

\begin{center}
  \begin{tabular}{ | l | l | l |}
    \hline
    Expression & Role  & Properties \\ \hline
    $\omega$ & $1/2+2\omega$ is an exponent of distribution & $\omega>0$ \\ \hline
    $\varepsilon$, $\delta$ & Small constants & $\delta>0$, $\varepsilon \in (0,10^{-100}\delta)$ \\ \hline
    $\alpha$, $\beta$ & Arbitrary coefficient sequences & $\alpha$ is located at scale $M$ \\ & & $\beta$ is located at scale $N$ \\ & & $\beta$ has the Siegel-Walfisz property
    \\ \hline
        $\psi_M$, $\psi_N$, $\psi_{\Delta_1}$, & Smooth coefficient sequences & $\psi_C$ is smooth at scale $C$ \\  $\psi_N^*$, $\psi_{\Delta_1}^*$  & &  $\psi_C^*$ is also smooth at scale $C$
       \\ \hline
       $\psi_{\Delta^*/z_1 },  \psi_{(\lambda^*/\lambda )N}$, & Shifted smooth coefficient sequences & $\psi_C$ is shifted smooth at scale $C$ \\ $\psi_{\Delta_2/z_1}, \psi_{N_2}$ & &  with corresponding shift $x_0 \ll x$
        \\ \hline
    $N$, $M$ & Scales of coefficient sequences $\alpha$ and $\beta$ & $x \ll MN \ll x$, $N= x^{\gamma(x)}$ \\ & &  with $\gamma(x)  \in (12\omega+6\delta,\frac{1}{2}-2\omega-8\varepsilon)$ \\ \hline
     $H$, $H^*$ & Length of summation over $h$& $H=x^\varepsilon RQ^2 (q_0 M)^{-1}$ 
  \\ & &   
     $1 \leq H \ll (x^{4\omega + \delta+7\varepsilon})/q_0$
 \\ & & $1 \ll |H^*| \ll H$
  \\ \hline
    $Q$, $R$, $U$, $V$ & Sizes of factors of a smooth modulus $q$  & $x^{-4\varepsilon-\delta}N \ll R \ll x^{-2\varepsilon}N$\\ & &
$x^{1/2-\varepsilon} \ll QR \ll x^{1/2+2\omega+\varepsilon}$\\ & &
$q_0^{-1} x^{-\delta-5\varepsilon}\frac{Q}{H} \ll U \ll q_0^{-1}x^{-5\varepsilon} \frac{Q}{H}$\\ & &
 $x^{5\varepsilon}H \ll V \ll x^{\delta+5\varepsilon}H$\\ & &  
 $UV \asymp Q/q_0$ 
 \\ \hline 
 $D$ & Approximate size of factor $d$ of  $q$   & $D= N/(x^{50\varepsilon}H^2) $
 \\ \hline 
 $\Delta$, $\Delta_1$, $\Delta^*$ & Shortened summation ranges for $d$   & $x^{-\delta} D \ll \Delta \ll D$
 \\ &  & $\Delta_1 = \Delta/(x^{5\varepsilon})$ \\ & &  $ \frac{N}{x^{\delta + 55\varepsilon}H^2} \ll \Delta_1 \ll \frac{N}{x^{ 55\varepsilon}H^2}  $
  \\ & & $\Delta^* = \min\{N/(\Lambda x^{5\varepsilon}), \Delta_1\}$
  \\ \hline
  $W_1$ & Size of $w_1$ & $1 \ll W_1 \ll \Delta$
 \\ \hline
 $Y$ & Size of $y$ and $\widetilde{y}$ & $1\ll |Y| \ll \frac{H^* V}{(v_1,v_2)}$
  \\ \hline
  $\Lambda$ &  Size of  $\lambda$ and $\widetilde{\lambda}$ & $1 \ll |\Lambda| \ll \frac{1}{w_1 (v_1,v_2)}x^{\delta +5\varepsilon} H^2$
  \\ \hline
$N_2$, $\Delta_2$ & (Used to simplify $\Sigma_8$ in Lemma~\ref{lem:adjustments2}) & $0<\Delta_2 \leq \Delta^*$   and $0<N_2 \leq N$
   \\ \hline
    $q_0, q_2, r_1$,    & Factors of a smooth modulus $q$ &   $q_0, r_1, u_1, v_1, v_2, q_2 \in \mathbb{N}$  
 \\ $u_1, v_1, v_2$ & &  $q_0 \ll Q$, $q_2 \asymp Q/q_0$,  $r_1 \asymp R/\Delta$ \\ & & $u_1 \asymp U$, $v_1, v_2 \asymp V$,  \\ \hline
 $v^*$, $v_1^*$, $v_2^*$  & Factors of $v_1$ and $v_2$ & $v^* = (v_1,v_2)$, $v_1^*=v_1/v^*$,  $v_2^*=v_2/v^*$
 \\ \hline  
 $d_0$ & $d=d_0+d_1$ is a factor of  $q$ &  $d_0 \in \mathbb{N}$, $d_0 \asymp \Delta$
 \\ \hline
 $m$ & Modulus of exponential $e_m(\,\cdot \,)$ & $m =  r_1q_0u_1[v_1,v_2]q_2$ \\ & & $m \mid P(x^\delta)$, $(m,ab_1b_2)=1$ \\ & & $\frac{RQ^2H}{q_0 (v_1,v_2)\Delta_1} \ll m \ll \frac{x^\delta RQ^2H}{q_0 (v_1,v_2)\Delta_1}$
  \\ \hline 
        $c_1$, $c_2$, $q_3$ & Factors of $m$ & $c_1 \mid m$, $c_2 \mid m$, $q_3 \mid \frac{m}{q_0}$
        \\ \hline 
  $m^*$ & Power of $m$ which exceeds $x$ & $m^*=m^{\lfloor 2\log(x) \rfloor}$
  \\ \hline 
          $w_0$, $w_1$, $w_2$, $z_1$  & Factors of $(h_1v_2^*-h_2v_1^*)$ or $(h_1v_2^*-h_2v_1^*)^2$ & $w_0, w_1, w_2, z_1\in \mathbb{N}$,  $w_1$ squarefree\\ & &  $w_1 \asymp W_1$,  $z_1 \ll x^{5\varepsilon} \Delta_1$\\ &  &  $w_0 \mid w_1$, $w_1 \mid z_1$, 
       $w_2 \mid m^*$   \\ &  &   $(w_1,m)=1$, $(z_1,m)=1$
        \\ \hline
      $y$, $\widetilde{y}$ & $y, \widetilde{y} \in \{h_1v_2^*-h_2v_1^*: h_1, h_2 \sim H^*\}$ & $y, \widetilde{y} \in \mathbb{Z} \setminus \{0\}$ and $y, \widetilde{y}  \sim Y$ \\ & &
        $(y,m^*) = (\widetilde{y},m^*)=w_2$, $w_1 \mid (y,\widetilde{y})$ \\ \hline 
  $\lambda$, $\widetilde{\lambda}$ & $\lambda = y/w_1$ and $\widetilde{\lambda} = \widetilde{y}/w_1$ &  $\lambda, \widetilde{\lambda} \in \mathbb{Z}\setminus \{0\}$ and $\lambda, \widetilde{\lambda} \sim \Lambda$ \\ & &  $(\lambda,m^*) = (\widetilde{\lambda},m^*)=w_2$
   \\ \hline
 $\ell$, $l$ & Integer shifts & $\ell, l \in \mathbb{Z} \setminus \{0\}$ \\ & & $ |\ell| \ll \frac{N}{R}$ and $|l| \ll N$
 \\
    \hline
        $a$, $b_1$, $b_2$ & Residue classes modulo $q$ & $a,b_1,b_2 \in \mathbb{Z}$ and $(q_0,ab_1b_2)=1$
          \\ \hline
          $A$, $A_1$, $A_2$, $B$, $B_1$ & Residue classes modulo $m$ & $(A,m)=(A_1,m)=(A_2,m)=1$
           \\ \hline
          $F_k$, $G_k$ & Residue classes modulo $m$ & (See Lemma~\ref{lem:adjustments2a} for details)
    \\ \hline
  \end{tabular}
\end{center}

\section{Key propositions} \label{sec:keysteps}

We now state our two main propositions, each of which corresponds to a section of this paper. We also show how Theorem~\ref{thm:equidistribution}, Theorem~\ref{thm:infinumboundedgaps} and Corollary~\ref{cor:H2} follow from these results.  

\subsection{Equidistribution estimates}

The following result improves on part (iii) of Theorem~2.8 of~\cite{Polymath:2014:EDZ} when $\sigma>8\omega+2\delta$ and $\sigma<24 \omega$.

\begin{prop} \label{prop:typeI/IIestimate}
Let $\omega, \delta, \sigma >0$. Suppose the following three inequalities are satisfied:
\begin{align*}
\begin{cases}
\,\,\, 72\omega + 24 \delta < 1, \\
\,\,\, 48 \omega + 16\delta + 4\sigma <1, \\ 
\,\,\, 64 \omega +20\delta+2\sigma<1.
\end{cases}
\end{align*}
   Let $\alpha, \beta: \mathbb{N}\times (1,\infty) \rightarrow \mathbb{R}$ be   coefficient sequences. Suppose that $\alpha$ is located at  scale $M$ and $\beta$ is located at  scale $N$ with $N(x) \in [x^{1/2-\sigma},x^{1/2}]$ and $M(x) \asymp \frac{x}{N(x)}$. Suppose  that $\beta$ has  the Siegel--Walfisz property.
   
   \vspace{3mm}
Then  for $x>1$, $a \in \mathbb{Z}$ and $A>0$, we have
\begin{align} \label{equ:statementI/II}
\sum_{\substack{q \leq x^{1/2+2\omega } \\ q \mid P(x^\delta) \\ (q,a)=1}} \Bigg| \sum_{n \equiv a (q)} (\alpha \star \beta)(n;x) - \dfrac{1}{\phi(q)} \sum_{(n,q)=1}(\alpha \star \beta)(n;x) \Bigg| \ll_{A, \omega, \delta}  \dfrac{x}{\log(x)^A}.
\end{align} 
\end{prop}
(Note: Like in~\cite{Polymath:2014:EDZ}, the implied constant in (\ref{equ:statementI/II}) does not depend on the specific choice of $\alpha$ and $\beta$, but does depend on the values of the implied constants in Definition~\ref{def:sequence}.)

\subsubsection{Comment}
The bound (\ref{equ:statementI/II}) is of a  form similar to  Polymath's $\mbox{Type}_I^{(i)}[\omega,\delta,\sigma]$ and $\mbox{Type}_{II}^{(i)}[\omega,\delta,\sigma]$ estimates (see~\cite{Polymath:2014:EDZ}, Definition~2.6). We only made the following modification: In part (iii) of Theorem~2.8 of~\cite{Polymath:2014:EDZ}, the sum was taken over $q \in \mathcal{D}_I^{(4)}(x^\delta)$, the set of $4$-tuply $x^\delta$-densely divisible integers, whereas we merely sum over $x^\delta$-smooth $q$. We made this change purely for convenience's sake and our proof could also be adapted to the case $q \in \mathcal{D}_I^{(3)}(x^\delta)$, the set of $3$-tuply $x^\delta$-densely divisible integers.

\vspace{3mm}
 Following the terminology of Polymath~\cite{Polymath:2014:EDZ}, we say that $\mbox{Type}_{I/II}^{(\infty)}[\omega,\delta,\sigma]$ holds if for all coefficient sequences $\alpha$ and $\beta$, located at scales $N$ and $M$, which satisfy the Siegel-Walfisz property and have $x^{1/2-\sigma} \leq N(x) \leq x^{1/2+\sigma}$, estimate (\ref{equ:statementI/II}) is true.  
  Proposition~\ref{prop:typeI/IIestimate} gives that $\mbox{Type}_{I/II}^{(\infty)}[\omega,\delta,\sigma]$ holds when 
\begin{align*}
72\omega + 24\delta <1 \qquad \mbox{ and } \qquad 
48 \omega +16\delta+4\sigma   < 1 \qquad \mbox{ and } \qquad 64\omega + 20\delta+2\sigma   <1.
\end{align*}
Using this new Type I estimate, we now prove Theorem~\ref{thm:equidistribution} as a corollary of Proposition~\ref{prop:typeI/IIestimate}.

\subsubsection{Proof of Theorem $1$ $($assuming Proposition $1)$} We  want to deduce that for any  $\varepsilon>0$ there exists $\delta >0$  such that for $x>1$,  $a \in \mathbb{Z}$ and $A>0$,
\begin{align}\label{equ:proofoftheorem1start}
\sum_{\substack{ q \leq x^{1/2 +1/40 -\varepsilon } \\ q \mid P(x^\delta) \\ (q,a)=1}} \Bigg| \sum_{\substack{ p \leq x \\ p \equiv a (q) }} 1 - \dfrac{1}{\phi(q)}\sum_{\substack{ p \leq x \\ (p,q)=1 }} 1  \,\,\,\, \Bigg| \ll_{A, \varepsilon} \dfrac{x}{\log(x)^A}.
\end{align} 
By partial summation and dyadic decomposition, the left-hand side of (\ref{equ:proofoftheorem1start}) can be bounded as follows:
\begin{align*}
\mbox{LHS} &\leq \sum_{\substack{ q \leq x^{1/2 +1/40 -\varepsilon } \\ q \mid P(x^\delta) \\ (q,a)=1}} \Bigg| \sum_{\substack{ p \leq x \\ p \equiv a (q) }} \!\frac{\log(p)}{\log(x)} - \dfrac{1}{\phi(q)}\!\sum_{\substack{ p \leq x \\ (p,q)=1 }} \!\frac{\log(p)}{\log(x)}   \Bigg| +\int_{2}^x\Bigg| \sum_{\substack{ p \leq t \\ p \equiv a (q) }}\! \frac{\log(p)}{t\log(t)^2} - \!\dfrac{1}{\phi(q)}\sum_{\substack{ p \leq t \\ (p,q)=1 }} \!\frac{\log(p)}{t\log(t)^2}   \Bigg|\,\mbox{d}t \\
&\leq  \sup_{t_0 \in [2,x]}\sum_{\substack{ q \leq x^{1/2 +1/40 -\varepsilon } \\ q \mid P(x^\delta) \\ (q,a)=1}}\Bigg| \sum_{\substack{ n \leq t_0 \\ n \equiv a (q) }}\! \Lambda(n)  - \!\dfrac{1}{\phi(q)}\sum_{\substack{ n \leq t_0 \\ (n,q)=1 }} \!\Lambda(n) \Bigg|\,\mbox{d}t + O_A\left( \dfrac{x}{\log(x)^A} \right) \\
&\leq \log(x) \sup_{x_0 \in [2,x]}\sum_{\substack{ q \leq x^{1/2 +1/40 -\varepsilon } \\ q \mid P(x^\delta) \\ (q,a)=1}}\Bigg| \sum_{\substack{ n \in [x_0,2x_0] \\ n \equiv a (q) }}\! \Lambda(n)  - \!\dfrac{1}{\phi(q)}\sum_{\substack{ n \in [x_0,2x_0]  \\ (n,q)=1 }} \!\Lambda(n) \Bigg|\,\mbox{d}t + O_A\left( \dfrac{x}{\log(x)^A} \right). 
\end{align*} 
For $x_0 \leq x^{1-\varepsilon}$, the above sum can  easily be bounded by $O(x/\log(x)^A)$. On the other hand, $x_0> x^{1-\varepsilon}$ gives $x^{1/2+1/40-\varepsilon} < x_0^{1/2+1/40-\varepsilon/4}$ and $x^\delta < x_0^{4\delta}$ when $\varepsilon<\frac{3}{4}$. The left-hand side of (\ref{equ:proofoftheorem1start}) has upper bound 
\begin{align*}
\mbox{LHS} \leq \log(x) \sup_{x_0 \in [2,x]}\sum_{\substack{ q \leq x_0^{1/2 +1/40 -\varepsilon/4 } \\ q \mid P(x_0^{4\delta}) \\ (q,a)=1}}\Bigg| \sum_{\substack{ n \in [x_0,2x_0] \\ n \equiv a (q) }}\! \Lambda(n)  - \!\dfrac{1}{\phi(q)}\sum_{\substack{ n \in [x_0,2x_0]  \\ (n,q)=1 }} \!\Lambda(n) \Bigg|\,\mbox{d}t + O_A\left( \dfrac{x}{\log(x)^A} \right).
\end{align*} 
To prove (\ref{equ:proofoftheorem1start}) it hence suffices to show that for sufficiently small $\delta >0$ and all $x>1$, $a \in \mathbb{N}$ and $A>0$, 
\begin{align} \label{equ:dyadic}
\sum_{\substack{ q \leq x^{1/2 +1/40 -\varepsilon/4 } \\ q \mid P(x^{4\delta}) \\ (q,a)=1}} \Bigg| \sum_{\substack{ n \in [x,2x] \\ n \equiv a (q) }} \Lambda(n) - \dfrac{1}{\phi(q)}\sum_{\substack{ n \in [x,2x] \\ (n,q)=1 }} \Lambda(n)   \Bigg| \ll_{A,  \delta} \dfrac{x}{\log(x)^{A+1}}.
\end{align}  
In the terminology of Polymath~\cite{Polymath:2014:VSS},  (\ref{equ:dyadic}) is referred to as a Motohashi-Pintz-Zhang estimate: We ask for $\mbox{MPZ}[\frac{1}{80} - \frac{\varepsilon}{8},4\delta]$ to hold.    Polymath~\cite{Polymath:2014:EDZ} gave a sufficient condition for when $\mbox{MPZ}[\omega,\delta]$ is true, dependent on the satisfaction of Type I/II/III estimates  $\mbox{Type}_I^{(i)}[\omega,\delta,\sigma]$, $\mbox{Type}_{II}^{(i)}[\omega,\delta]$ and  $\mbox{Type}_{III}^{(i)}[\omega,\delta,\sigma]$, described in Definition~2.6 of~\cite{Polymath:2014:EDZ}. This result allowed for  an extended range of summation, summing over $q \in \mathcal{D}_I^{(i)}(x^\delta)$. However, since we only use smooth moduli in this paper, we may replace $\mbox{Type}_I^{(i)}[\omega,\delta,\sigma]$ and $\mbox{Type}_{II}^{(i)}[\omega,\delta]$ by $\mbox{Type}_{I/II}^{(\infty)}[\omega,\delta, \sigma]$ and obtain the following criterion:

\begin{lemma}[Combinatorial Lemma]\label{lem:combinatoriallemma}
Let  $\omega \in (0,\frac{1}{4})$,  $\delta \in (0,\frac{1}{4}+\omega)$ and $\sigma \in (\frac{1}{10}, \frac{1}{2})$  with $\sigma > 2\omega$. Suppose that $\mbox{\rm Type}_{I/II}^{(\infty)}[\omega,\delta,\sigma]$ and  $\mbox{\rm Type}_{III}^{(1)}[\omega,\delta,\sigma]$ hold. Then $\mbox{\rm MPZ}[\omega,\delta]$ holds.
\end{lemma}

\begin{proof}
This follows directly from the proof of Lemma~2.7 of~\cite{Polymath:2014:EDZ} by replacing $\mathcal{Q} = \{q \leq Q: q \in \mathcal{D}_I^{(i)}(x^\delta)\}$ with $\mathcal{Q} = \{q \leq Q: q \mid P(x^\delta)\}$ everywhere.
\end{proof}

Using Proposition~\ref{prop:typeI/IIestimate}, Lemma~\ref{lem:combinatoriallemma} and part (v) of Theorem~2.8 of~\cite{Polymath:2014:EDZ}, we deduce the following:

\begin{cor}\label{cor:goodomegaanddelta}
Let $\omega>0$ and $ \delta >0$ with $80 \omega +\frac{80}{3} \delta  <  1$. Then $\mbox{\rm MPZ}[\omega,\delta]$ holds.
\end{cor}

\begin{proof}
By Proposition~\ref{prop:typeI/IIestimate}, $\mbox{\rm Type}_{I/II}^{(\infty)}[\omega,\delta,\sigma]$ holds if $72\omega +24\delta <1$ and $48 \omega +16\delta+4\sigma  < 1$ and $64\omega + 20\delta+2\sigma  <1$. The second and third inequality rearrange to  $\sigma < \frac{1}{4}-12 \omega -4 \delta $ and $\sigma <\frac{1}{2}-32\omega - 10\delta$. By part~(v) of Theorem~2.8 of~\cite{Polymath:2014:EDZ}, $\mbox{\rm Type}_{III}^{(1)}[\omega,\delta,\sigma]$ holds for $\omega < \frac{1}{12}$ when $\sigma > \frac{1}{18} + \frac{28 \omega}{9} + \frac{2 \delta}{9}$.   By Lemma~\ref{lem:combinatoriallemma}, $\mbox{MPZ}[\omega,\delta]$ thus holds provided that $72\omega +24\delta <1$ and provided the following  inequality is satisfied:
\begin{align*}
&\max\left\{\frac{1}{18} + \frac{28 \omega}{9} + \frac{2 \delta}{9}, \frac{1}{10}\right\} < \min\left\{ \frac{1}{4}-12 \omega -4\delta, \frac{1}{2}-32\omega - 10\delta \right\}.
\end{align*} 
 After some careful computations, this condition rearranges to  $80 \omega +\frac{80}{3}\delta <  1$. 
\end{proof}

 \vspace{3mm}

\emph{Proof of Theorem $1$ $($assuming Prop. $1)$ -- Final Step:} Recall that (\ref{equ:dyadic}) corresponds to $\mbox{MPZ}[\frac{1}{80}-\frac{\varepsilon}{8},4\delta]$. Hence we have the desired bound (\ref{equ:dyadic}) whenever  $\varepsilon$ and $\delta$ satisfy $80(\frac{1}{80}-\frac{\varepsilon}{8}) + \frac{320}{3}\delta = 1 -10\varepsilon + \frac{320}{3}\delta<1$. Choosing $\delta = \frac{\varepsilon}{100}$, this inequality holds. Theorem~\ref{thm:equidistribution} follows. \qed

\vspace{3mm}
 In~\cite{Polymath:2014:EDZ} and~\cite{Polymath:2014:VSS} it was shown and used that $\mbox{MPZ}[\omega,\delta]$ holds when $\frac{600 \omega}{7} + \frac{180 \delta}{7} <1$. Our alternative criterion $80 \omega +\frac{80}{3}\delta <  1$, stated in Corollary~\ref{cor:goodomegaanddelta}, is  an improvement over Polymath's result when $ \omega  >   \frac{\delta}{6}$.

\subsubsection{A first improvement of bounds on $H_m$} \label{ssec:hb} 
 We recall now from~\cite{Polymath:2014:VSS} some of the main steps of Maynard's and Tao's method for bounding $H_m$. The weak Dickson-Hardy-Littlewood prime tuples conjecture $\mbox{DHL}[k,m+1]$ (described in Claim~3.1 of~\cite{Polymath:2014:VSS}) plays a key role: If $\mbox{DHL}[k,m+1]$ holds, then $H_m \leq H(k)$, where  $H(k)$ is the diameter of the narrowest admissible $k$-tuple.
 
 \vspace{3mm} By Theorem~3.10 of~\cite{Polymath:2014:VSS}, $\mbox{DHL}[k,m+1]$ holds if $\mbox{MPZ}[\omega,\delta]$ holds and $M_k^{[\frac{\delta}{1/4+\omega}]} > \frac{m}{1/4+\omega}$, where $M_k^{[\alpha]}$ is as defined in equation (34) of~\cite{Polymath:2014:VSS}.  By the proof of  Theorem~6.7 of~\cite{Polymath:2014:VSS}, $M_k^{[\alpha]} > \log(k)-O_\alpha(1)$. Finally, by Theorem~3.3 of~\cite{Polymath:2014:VSS},  $H(k) \leq k\log k + k\log\log k +o(k)$. Combining all this with Corollary~\ref{cor:goodomegaanddelta}, we have 
 \begin{align*}
 H_m \ll_\epsilon  \exp\left(\frac{(1+\epsilon)m}{(1/4 + 1/80)} \right).
 \end{align*}
 In particular, we have $H_m \ll \exp(3.81m)$, as claimed in the introduction.

\subsubsection{Proof of Corollary $\ref{cor:H2}$}
To get bounds on $H_2, \dots , H_5$, we again follow the arguments of Polymath~\cite{Polymath:2014:VSS}. 

\vspace{3mm} In particular, we use Theorem~6.7 of~\cite{Polymath:2014:VSS} with 
\begin{align*}
(k,c \log(k),T \log(k)) = \begin{cases} 
(35265,  0.99479, 0.85213) &\mbox{ for }  m=2, \\
(1624545, 1.00422, 0.80148)  &\mbox{ for }  m=3, \\
(73807570,  1.00712, 0.77003)  &\mbox{ for }  m=4, \\
(3340375663,  1.0079318,  0.7490925)  &\mbox{ for } m=5,
\end{cases}
\end{align*}
to deduce that  
$M_{k}^{[\frac{\delta}{1/4+\omega}]} > \frac{m}{1/4+\omega}$ holds for the following choices of $\omega$ and $\delta$:
\begin{align*}
(\omega,\delta) = \begin{cases}
(0.00556625,  0.0207987) &\mbox{ for } (m,k) = (2,35265), \\
(0.00768602,  0.0144419) &\mbox{ for } (m,k) = (3,1624545), \\
(0.00883292,  0.0110012 ) &\mbox{ for } (m,k) = (4,73807570), \\
(0.0095447064335,  0.0088658806979) &\mbox{ for } (m,k) = (5,3340375663).
\end{cases}
\end{align*}
Each of these choices of $\omega$ and $\delta$ satisfies $80 \omega +\frac{80}{3}\delta  <  1$, so that $\mbox{MPZ}[\omega,\delta]$ holds. Thus $\mbox{DHL}[35265,3]$, $\mbox{DHL}[1624545,4]$, $\mbox{DHL}[73807570,5]$ and $\mbox{DHL}[3340375663,6]$ are true.   
In particular, $H_2 \leq H(35265)$, $H_3 \leq H(1624545)$, $H_4 \leq H(73807570)$ and $H_5 \leq H(3340375663)$. Sutherland computed the following:
\begin{align*}
&H(35265) \leq 396504, \\
&H(1624545) \leq 24407016,  \\
&H(73807570) \leq 1391051532, \\
&H(3340375663) \leq 77510685234.
\end{align*} 
 For $k \in \{35265, 1624545,  73807570\}$, admissible $k$-tuples of lengths $396504,  24407016, 1391051532$ can be found \href{https://www.dropbox.com/scl/fo/9my1skgb657z2sc49tqls/h?rlkey=pq8snurxku4o1dbi1ta0ltqbf&dl=0}{{\color{blue}\underline{here}}}. To compute the bound $H(3340375663) \leq 77510685234$, the  shifted Schinzel sieve was used as described in Section~10.2.2 of Polymath~\cite{Polymath:2014:VSS}, with $ [s,s+x] =
[3343896484,80854581718]$ and $m=13090609$.
This concludes the proof of Corollary~\ref{cor:H2}. \qed

\subsection{Harman's sieve} 
We now use Harman's sieve to further improve our bounds on $H_m$. 
Section~\ref{sec:harman} is focused on the construction of a suitable minorant for $1_{\mathbb{P}}(n)$, using as input the Type I/II equidistribution estimate given in  Proposition~\ref{prop:typeI/IIestimate} and the Type~III equidistribution estimate given in part (v) of Theorem~2.8 of~\cite{Polymath:2014:EDZ}. Write 
\begin{align*}
\psi(n,y) = \begin{cases} 
\,\, 1  \quad \mbox{ if } p \mid n \rightarrow p \geq y, \\ \,\, 0 \quad \mbox{ otherwise. }
\end{cases}
\end{align*}
We will prove the following: 

\begin{prop} \label{prop:construction}
For any  $n \in \mathbb{N}$ and $x>1$, write $p_j = x^{\alpha_j}$ and set
\begin{align*}
\rho(n;x) = 1_{\mathbb{P}}(n) - \sum_{\substack{ n = p_1p_2p_3p_4n_5 \\ 0.19038\leq \alpha_4 < \alpha_3 < \alpha_2 < \alpha_1 < 0.40481 \\ \alpha_1+\alpha_2 < 0.40481 \\ \alpha_2+\alpha_3+\alpha_4 >0.59519}} \psi(n_5, p_4)- \sum_{\substack{ n = p_2p_3p_4p_5p_6 \\ 0.19038 \leq \alpha_2, \alpha_3, \alpha_4, \alpha_5, \alpha_6 \leq 0.23848  \\ \alpha_2, \alpha_4 > \alpha_3  \\ \alpha_2+\alpha_4 < 0.40481\\ \alpha_2 + \alpha_3 + \alpha_5 >0.59519 \\   \alpha_6 \geq \alpha_5 }} 1.
\end{align*}
Then $\rho(n;x)$ is a minorant for $1_{\mathbb{P}}(n)$, has exponent of distribution $0.5253 $ to smooth moduli, and
\begin{align*}
\sum_{n \in [x,2x]} \rho(n;x) \geq (1- 2 \cdot 10^{-5} +o(1)) \dfrac{x}{\log(x)}.
\end{align*}
\end{prop}

\subsubsection{Proof of Theorem $2$ $($assuming Proposition $2)$} The following lemma, proved by Baker and Irving in~\cite{Baker:2017:BIP}, combines results of~\cite{Polymath:2014:VSS} with Harman's sieve: 

\vspace{3mm}
\begin{lemma}[Baker, Irving]~\label{lem:maynardtaoharman}
Fix a small $\xi >0$ and some $c_1 \in (0,1)$. Suppose 
there exists a function $\rho:\mathbb{N} \times (1,\infty) \rightarrow \mathbb{R}$ with the following properties:
\begin{enumerate}[{\rm (i) }]
\item For every $x$, $\rho(n;x)$ is a minorant for the indicator function of primes, that is, $\rho(n;x) \leq 1_{\mathbb{P}}(n)$.
\item If $n \in [x,2x]$ and $\rho(n;x) \neq 0$, then all prime factors of $n$ exceed $x^\xi$. 
\item The function $\rho(n;x)$ has exponent of distribution $\theta$ to smooth moduli.
\item The sum of $\rho(n;x)$ over $[x,2x]$ satisfies
\begin{align*}
\sum_{n \in [x,2x]} \rho(n;x) = (1-c_1 +o(1)) \dfrac{x}{\log(x)}.
\end{align*}
\end{enumerate}
Then we have, for every $\epsilon >0$, 
\begin{align*}
H_m = \liminf_{n \rightarrow \infty} (p_{n+m}-p_n) \ll_\epsilon \exp\left(\left(\dfrac{2(1+\epsilon) }{\theta(1-c_1)}\right)m \right).
\end{align*}
\end{lemma}

\begin{proof}
This is Lemma~1 of~\cite{Baker:2017:BIP}. 
\end{proof}

  Lemma~\ref{lem:maynardtaoharman} and  Proposition~\ref{prop:construction} give $H_m \ll \exp(3.8075 m)$, concluding the proof of Theorem~\ref{thm:infinumboundedgaps}.

\section{Equidistribution estimates} \label{sec:equidistribution}

In this section we consider a coefficient sequence $\alpha$ at scale $M$ and a coefficient sequence $\beta$ at scale $N$ with $x \ll M(x)N(x) \ll x$. We  assume  $\beta$ has the Siegel--Walfisz property. We write $N(x)=x^{\gamma(x)}$. Our aim is to prove the following equidistribution estimate for certain  $\omega>0$, $\delta>0$ and $\gamma(x) \in (12\omega+6\delta,\frac{1}{2}-2\omega-8\varepsilon)$: For every $x>1$, $a \in \mathbb{Z}$ and $A >0$,
\begin{align} \label{finalresult}
\sum_{\substack{q \leq x^{1/2+2\omega} \\ q \mid P(x^\delta) \\ (q,a)=1}}\Bigg| \sum_{n \equiv a(q)} (\alpha \star \beta)(n;x) - \dfrac{1}{\phi(q)} \sum_{(n,q)=1} (\alpha \star \beta)(n;x)\Bigg| \ll_{A,\varepsilon}  \dfrac{x}{\log(x)^A}.
\end{align}
(Note: To make our notation less crowded,  we will usually follow the convention of Polymath~\cite{Polymath:2014:EDZ} whereby the dependence of expressions on $x$ is suppressed, for instance by writing $\alpha(n)$ instead of $\alpha(n;x)$ and $N$ instead of $N(x)$.)

\subsection{A summary of Polymath's results}

We begin by summarizing some of the results of Section~5 and Section~8A of Polymath~\cite{Polymath:2014:EDZ}. These results will provide a starting point for our proof. 

\begin{defn}\label{def:Z1}
For  given $\omega,\varepsilon, \delta, M$ and $N$,  the expression $Z_1(x)$ denotes the set of tuples $$(Q,R,U,V,H^*,q_0,\ell,a,b_1,b_2,\psi_M,\psi_N)$$ which satisfy the following seven properties: 
\begin{enumerate}[{\rm (i)}]
\item $Q$, $R$, $U$ and $V$ are positive real numbers.
\item $H^*$ is a non-zero real number with  $1 \ll |H^*| \ll H$ for  $H := x^\varepsilon RQ^2(q_0 M)^{-1}$.
\item $q_0$ is a positive integer with $q_0 \ll Q$.
\item  $\ell \in \mathbb{Z}$ satisfies $0 \neq |\ell| \ll \frac{N}{R}.$
\item $a, b_1, b_2 \in \mathbb{Z}$ with  $(q_0,ab_1b_2)=1.$
\item $\psi_M$ and $\psi_N$ are real, non-negative coefficient sequences, smooth at scales $M$ and $N$.
\item  The following 6 bounds are satisfied:
\begin{align} 
&H = x^\varepsilon RQ^2(q_0 M)^{-1}\geq 1,\\
&x^{-4\varepsilon-\delta}N \ll R \ll x^{-2\varepsilon}N, \label{equ:rangeQ}\\
&x^{1/2-\varepsilon} \ll QR \ll x^{1/2+2\omega+\varepsilon}, \label{equ:rangeR}\\
\label{def:UV1}
&q_0^{-1} x^{-\delta-5\varepsilon}Q/H \ll U \ll q_0^{-1}x^{-5\varepsilon} Q/H,\\ \label{def:UV2}
&x^{5\varepsilon}H \ll V \ll x^{\delta+5\varepsilon}H,\\ \label{def:UV3}
&UV \asymp Q/q_0.
\end{align} 
\end{enumerate}
\end{defn}
(The implied constants in this definition are   large and independent of $x$.)

\begin{defn}\label{def:Z1functions}
For a given $z=(Q,R,U,V,H^*,q_0,\ell,a,b_1,b_2,\psi_M,\psi_N)\in Z_1(x)$, we further define:
\begin{align*}
&C_0(n) = 1_{\substack{\frac{b_1}{n}  \equiv \frac{b_2}{n+\ell r} (q_0)}},
\\
&\varphi_H^{\star}(h)=  \frac{1}{M}
\sum_m \psi_M(m)   e\!\left(\frac{-mh}{rq_0u_1v_1q_2} \right) \mbox{ and } 
\varphi_H^{\star\star}(h)=  \frac{1}{M}
\sum_m \psi_M(m)   e\!\left(\frac{-mh}{rq_0u_1v_2q_2} \right),\\
&\Psi(n,h_1,h_2,r,u_1,v_1,v_2,q_2)=e_r\left(\dfrac{ah_1}{nq_0u_1v_1q_2}\right)e_{q_0u_1v_1}\left(\dfrac{b_1h_1}{nrq_2}\right)e_{q_2}\left(\dfrac{b_2h_1}{(n+\ell r)rq_0u_1v_1}\right)\\
&\qquad\qquad\qquad\qquad\qquad\qquad\,\,\, e_r\left(-\dfrac{ah_2}{nq_0u_1v_2q_2}\right)e_{q_0u_1v_2}\left(-\dfrac{b_1h_2}{nrq_2}\right)e_{q_2}\left(-\dfrac{b_2h_2}{(n+\ell r)rq_0u_1v_2}\right). 
\end{align*}

\end{defn}

\begin{lemma}[Reduction of Type I bounds to exponential sums~\cite{Polymath:2014:EDZ}]\label{lem:summaryofpolymath}
Let $\omega, \delta >0$. Let $\varepsilon \in (0, 10^{-100}\delta)$. 

 Let $\alpha$ and $\beta$ be coefficient sequences at scales $M$ and $N$ with $x \ll M(x)N(x) \ll x$.  Let $N= x^{\gamma(x)}$, where $\gamma(x)  \in (12\omega+6\delta,\frac{1}{2}-2\omega-8\varepsilon)$. Assume $\beta$ has the Siegel-Walfisz property.

\vspace{3mm}

Let $Z_1(x)$ be as given in Definition~{\rm\ref{def:Z1}} and let  $C_0(n)$, $\varphi^{\star}_H(h)$, $\varphi^{\star\star}_H(h)$ and $\Psi(n,h_1,h_2,r,u_1,v_1,v_2,q_2)$ be as given in Definition~{\rm\ref{def:Z1functions}}. For $x>1$ and $z=(Q,R,U,V,H^*,q_0,\ell,a,b_1,b_2,\psi_M,\psi_N)\in Z_1(x)$, we then define
\begin{align*}
&\Sigma_1(z;x)= \!\!\!\!\!\!\!\sum_{\substack{r, u_1, v_1, v_2, q_2 \\ r \asymp R \\ u_1 \asymp U\\ v_1, v_2 \asymp V  \\ q_2 \asymp Q/q_0 \\
rq_0u_1[v_1,v_2]q_2  \mid P(x^\delta) \\ (rq_0u_1v_1v_2q_2,ab_1b_2)=1}}\!\!\!\!\!\!\!\!\Big| \sum_{\substack{h_1, h_2 \\ h_1 \sim H^* \\ h_2 \sim H^* }} \sum_{\substack{n \\ (n,rq_0u_1v_1v_2 )=1 \\ (n+\ell r, q_0q_2)=1 }} \!\!\!\!\!\!\! C_0(n) \varphi^{\star}_H(h_1) \overline{\varphi^{\star\star}_H(h_2)}\psi_N(n) \Psi(n,h_1,h_2,r,u_1,v_1,v_2,q_2) \Big|.
\end{align*}
Suppose that  $\Sigma_1(z;x) = O_{\varepsilon}((q_0,\ell)RQN UV^2x^{-4\varepsilon}) $ for  $x>1$ and  $z \in Z_1$. Then for all $a_0 \in \mathbb{Z}$ and $A>0$,
\begin{align*} 
\sum_{\substack{q \leq x^{1/2+2\omega} \\ q \mid P(x^\delta) \\ (q,a_0)=1}}\Bigg| \sum_{n \equiv a_0(q)} (\alpha \star \beta)(n;x) - \dfrac{1}{\phi(q)} \sum_{(n,q)=1} (\alpha \star \beta)(n;x)\Bigg| \ll_{A, \varepsilon}  \dfrac{x}{\log(x)^A}.
\end{align*}
\end{lemma}

This lemma follows almost directly from the arguments of Section~5 and Section~8A of~\cite{Polymath:2014:EDZ}. Some small changes have to be made, so for completeness' sake a full proof of Lemma~\ref{lem:summaryofpolymath}, starting from the point at which we first diverge from~\cite{Polymath:2014:EDZ},  can be found in the appendix of this paper. 

\vspace{3mm} (Note: For clarity, we highlight here once more that the upper bound on the equidistribution estimate, the upper bound on $\Sigma_1$, and also the upper bounds on $\Sigma_i$, $i \in \{2, \dots, 8\}$, in later lemmas, do not depend on any specific coefficient sequences, but are  allowed to depend on the implied constants in Definition~\ref{def:sequence}.)

\subsection{The $q$-van der Corput method}\label{sssec:twoaims} 

We now prepare to carry out the variant of the $q$-van der Corput method which was described in Section~\ref{sssec:variant}.  

\subsubsection{Diagonal terms $+$ Splitting r}

Our first lemma has two main aims: On the one hand, we remove certain inconvenient terms from our sums. (To be precise, we remove combinations of $h_1, h_2, v_1, v_2$ with $h_1v_2=h_2v_1$.) On the other hand, we split up $r$, writing $r=d r_1$ for $d$ and $r_1$ of a certain good size. The factor $d$ will later  be incorporated into the $q$-van der Corput method to reduce the modulus. 

\begin{defn}\label{def:Z2}
For  given $ \omega, \varepsilon, \delta, M, N$ and given $x>1$ and  $z_1 \in Z_1(x)$, we set 
\begin{align} \label{def:DII}
D :=  \dfrac{N}{x^{50\varepsilon}H^2}. 
\end{align}
The expression $Z_2(z_1;x)$ $($with $x>1$ and $z_1 \in Z_1(x))$ denotes the set of tuples $$(\Delta,r_1,u_1, v_1,v_2,q_2,d_0,\psi_{\Delta_1})$$ which satisfy the following four properties:
 \begin{enumerate}[{\rm (i)}]
 \item $\Delta$ is a positive real number with $x^{-\delta}D \ll \Delta \ll  D $. 
 \item   $r_1, u_1, v_1, v_2, q_2, d_0  \in \mathbb{N}$ with   $r_1q_0u_1[v_1,v_2]q_2 \mid P(x^\delta)$ and $(r_1q_0u_1v_1v_2q_2,ab_1b_2)=1$. 
 \item   $r_1 \asymp \frac{R}{\Delta}$, $u_1 \asymp U$,  $v_1,v_2 \asymp V$,  $q_2 \asymp \frac{Q}{q_0}$ and $d_0 \asymp \Delta$. 
 \item $\psi_{\Delta_1}(d_1)$ is a real, non-negative coefficient sequence, smooth at scale $\Delta_1=\Delta/(x^{5\varepsilon})$.
 \end{enumerate}
 
 \smallskip
 Finally, for given $z_1 \in Z_1(x)$ and $z_2 \in Z_2(z_1;x)$, we set  $C(n) = 1_{\substack{\frac{b_1}{n}  \equiv \frac{b_2}{n+\ell d r_1} (q_0)}}$. 
\end{defn}

\begin{lemma}[Extraction of a suitable factor $d$]\label{lem:preparationII}
Let $\Sigma_1$, $Z_1$ and $Z_2$   be as described in Lemma~{\rm \ref{lem:summaryofpolymath}}, Definition~{\rm \ref{def:Z1}} and Definition~~{\rm \ref{def:Z2}}.
For given $x>1$, $z_1 \in Z_1(x)$ and $z_2\in Z_2(z_1;x)$,  we set 
\begin{align*}
&\Sigma_2(z_1,z_2;x)=\!\!\!\!\!\!\!\!\!\!\!\!\!\!\!\!\!\! \sum_{\substack{d  \\ (d,r_1q_0u_1[v_1,v_2]q_2)=1 \\ d \mbox{\begin{scriptsize}
squarefree
\end{scriptsize}}
}} \!\!\!\!\!\!\!\!\!\!\!\!\!\!\!\!\!\psi_{\Delta_1}(d-d_0)\, \Big|\!\!\!\! \sum_{\substack{h_1, h_2 \\ h_1 \sim H^* \\ h_2 \sim H^*  \\ h_1 v_2 \neq h_2v_1 }}  \!\!\!\sum_{\substack{n }}^{\,\,*_2} \,C(n)\varphi_H^\star(h_1) \overline{\varphi_H^{\star\star}(h_2)}\psi_N(n)  \Psi(n,h_1,h_2,dr_1,u_1,v_1,v_2,q_2) \Big|,
\end{align*}

\vspace{-3mm}
where $\sum\limits^{\,\,*_2}_n$ denotes summation over $n$ with $(n,dr_1q_0u_1v_1v_2  )=1$ and $(n+\ell dr_1, q_0q_2)=1$.

\vspace{1mm}
If $\Sigma_2(z_1,z_2;x) = O_{\varepsilon}\!\left(\dfrac{q_0(q_0,\ell)\Delta N (v_1,v_2)}{x^{10\varepsilon}}\right)$ for all $z_2 \in Z_2$, then $\Sigma_1(z_1;x) = O_{\varepsilon}\!\left(\dfrac{(q_0,\ell)RQNUV^2}{x^{4\varepsilon}}\right)$. 
\end{lemma}

\begin{proof}
Our aim is to bound $\Sigma_1(z_1)$, which is given by 
\begin{align*}
\Sigma_1(z_1)=\!\!\!\!\!\!\!\!\!\!\! \sum_{\substack{r, u_1, v_1, v_2, q_2 \\ r \asymp R, u_1 \asymp U\\ v_1, v_2 \asymp V, q_2 \asymp Q/q_0 \\
rq_0u_1[v_1,v_2]q_2  \mid P(x^\delta) \\ (rq_0u_1v_1v_2q_2,ab_1b_2)=1}}\!\!\!\Big| \sum_{\substack{h_1, h_2 \\ h_1 \sim H^* \\ h_2 \sim H^* }} \sum_{\substack{n \\ (n,rq_0u_1v_1v_2 )=1 \\ (n+\ell r, q_0q_2)=1 }} \!\!\!\!\!\!\! C_0(n) \varphi_H^\star(h_1) \overline{\varphi_H^{\star\star}(h_2)}\psi_N(n) \Psi(n,h_1,h_2,r,u_1,v_1,v_2,q_2) \Big|.
\end{align*}
We now factorise $r$. Recall that $R \gg x^{-4\varepsilon-\delta} N$.  Hence $D \leq x^{-50\varepsilon}N \ll x^{-46\varepsilon+\delta}R$. On the other hand,
\begin{align} \label{h-inequality}
 H = \frac{x^\varepsilon R Q^2 }{q_0 M } = \frac{ x^\varepsilon ( R Q)^2}{q_0  M R } \ll  \frac{(x^{1+4\omega+3\varepsilon}) (x^{4\varepsilon+\delta}) }{q_0 M N } \ll \frac{x^{4\omega+\delta+7\varepsilon}  }{q_0  }
\end{align}
and so $D = x^{-50\varepsilon} N/H^2 \gg x^{-8\omega-2\delta-64\varepsilon} N $. Using that  $\gamma > 12\omega+6\delta$, then $D\gg 1$. For $r \in \mathbb{N}$ with $r \mid P(x^\delta)$ and $r \asymp R$, we can thus find $d$ and $r_1$ such that $r=dr_1$ with $x^{-\delta}D \leq d \leq D$. 
Write
 \begin{align*}
\Upsilon_0(\Delta,r_1,u_1,q_2)=\sum^{\,\,*_1}_{\substack{d, v_1, v_2 \\d \asymp \Delta \\ v_1, v_2 \asymp V 
}} \!\!\Big| \sum_{\substack{h_1, h_2 \\ h_1 \sim H^* \\ h_2 \sim H^* }}  \sum^{\,\,*_2}_{\substack{n  }} \, C(n) \varphi_H^\star(h_1) \overline{\varphi_H^{\star\star}(h_2)}\psi_N(n) \Psi(n,h_1,h_2,dr_1,u_1,v_1,v_2,q_2) \Big|.
\end{align*} 

\vspace{-3mm} 
(Here $\sum\limits^{\,\,*_1}$ indicates that the outer sum is restricted to choices of $ d,  v_1, v_2$ with  $dr_1q_0u_1[v_1,v_2]q_2  \mid P(x^\delta)$ and $(dr_1q_0u_1v_1v_2q_2,ab_1b_2)=1$ and $\sum\limits^{\,\,*_2}$ indicates that the sum over $n$ is restricted to $(n,dr_1q_0u_1v_1v_2  )=1$ and $(n+\ell dr_1, q_0q_2)=1$.)
Then by dyadic decomposition, \begin{align*}
\Sigma_1(z_1) \ll_\varepsilon \left(\frac{x^{0.5\varepsilon} RUQ}{q_0 \Delta} \right) \sup \Upsilon_0 (\Delta,r_1,u_1,q_2),
\end{align*}
where the supremum is taken over all choices of $\Delta$ with $x^{-\delta}D \ll \Delta \ll D$ and all choices of $r_1, u_1, q_2$ with $r_1 \asymp R/\Delta$, $u_1 \asymp U$, $q_2 \asymp Q/q_0$ and $(r_1u_1q_2,ab_1b_2)=1$. To show that $\Sigma_1(z_1) \ll_{\varepsilon}x^{-4\varepsilon}(q_0,\ell)RQNUV^2$, it thus suffices to show that $\sup \Upsilon_0  (\Delta,r_1,u_1,q_2) \ll_{\varepsilon} x^{-4.5\varepsilon}q_0(q_0,\ell)\Delta N V^2 $.

\vspace{3mm}
Next we wish to exclude choices of $h_i$ and $v_i$ for which $h_1v_2-h_2v_1=0$. We have $\Upsilon_0 \leq \Upsilon_1 + \Upsilon_2$, where   $\Upsilon_1$ and $\Upsilon_2$ are given by 
  \begin{align*}
&\Upsilon_1=\!\!\!\sum^{\,\,*_1}_{\substack{d, v_1, v_2 \\d \asymp \Delta \\ v_1, v_2 \asymp V 
}} \!\!\Big| \sum_{\substack{h_1, h_2 \\ h_1, h_2 \sim H^* \\ h_1 v_2 - h_2 v_1 = 0 }}  \sum_{\substack{n  }}^{\,\,*_2}\, C(n) \varphi_H^\star(h_1) \overline{\varphi_H^{\star\star}(h_2)}\psi_N(n) \Psi(n,h_1,h_2,dr_1,u_1,v_1,v_2,q_2) \Big|,\\
&\Upsilon_2=\!\!\!\sum^{\,\,*_1}_{\substack{d, v_1, v_2 \\d \asymp \Delta \\ v_1, v_2 \asymp V 
}} \!\!\Big| \sum_{\substack{h_1, h_2 \\ h_1, h_2 \sim H^* \\ h_1 v_2 - h_2 v_1 \neq 0 }}  \sum_{\substack{n  }}^{\,\,*_2}\, C(n) \varphi_H^\star(h_1) \overline{\varphi_H^{\star\star}(h_2)}\psi_N(n) \Psi(n,h_1,h_2,dr_1,u_1,v_1,v_2,q_2) \Big|.
\end{align*}
Recall that $H^* \ll |H| \ll x^{-5\varepsilon}V$. A trivial bound on  $\Upsilon_1$ is thus given by 
$$
\Upsilon_1  \ll x^{o(1)} \Delta H^*V N   \ll x^{-5\varepsilon+o(1)} \Delta NV^2  
.$$ To get the desired bound on $\Sigma_1(z_1)$, it thus suffices to show that $\Upsilon_2 \ll_{\varepsilon} x^{-4.5\varepsilon} q_0(q_0, \ell)\Delta NV^2 $. Consider 
\begin{align*}
\Upsilon_2'(v_1,v_2)=\!\!\!\!\!\!\!\!\!\!\!\!\!\!\!\!\!\!\sum_{\substack{d \asymp \Delta \\ (d,r_1q_0u_1[v_1,v_2]q_2)=1  \\  d \mbox{\begin{scriptsize}
squarefree
\end{scriptsize}}
}}\!\!\! \Big|\! \sum_{\substack{h_1, h_2 \\ h_1, h_2 \sim H^*  \\ h_1 v_2-h_2v_1 \neq 0}}  \sum_{\substack{n }}^{\,\,*_2} \, C(n)\varphi_H^\star(h_1) \overline{\varphi_H^{\star\star}(h_2)}\psi_N(n)  \Psi(n,h_1,h_2,dr_1,u_1,v_1,v_2,q_2) \Big|
\end{align*}
for  $v_1$, $v_2$ with $r_1 q_0u_1[v_1,v_2]q_2 \mid P(x^\delta)$ and $(v_1v_2,ab_1b_2)=1$. (Note here that we removed the conditions $d \mid P(x^\delta)$ and $(d,ab_1b_2)=1$. This increases the sum.) We now have
\begin{align*}
\Upsilon_2 \ll \sum_{v_1 \asymp V} \sum_{v_2 \asymp V} \Upsilon_2'(v_1,v_2).
\end{align*}

For technical reasons, it will later be important to have a shorter summation range for $d$. Currently we are looking at $d \asymp \Delta$.  Recall that $\Delta_1 = \Delta/( x^{5\varepsilon})$. We partition the set of $d$ with $d \asymp \Delta$ into intervals  $\{d: \Delta_1 < d-d_0\leq 2\Delta_1\}$ with $d_0 \asymp \Delta$. We only need $O(x^{5\varepsilon})$ choices of $d_0$ to cover the  range $d \asymp \Delta$.  Set 
 \begin{align*}
\Upsilon_3(v_1,v_2,d_0) = \!\!\!\!\!\!\!\!\!\!\!\!\!\!\!\!\!\!&\sum_{\substack{d  \\ (d,r_1q_0u_1[v_1,v_2]q_2)=1 \\ d \mbox{\begin{scriptsize}
squarefree
\end{scriptsize}}}} \!\!\!\!\!\!\!\!\!\!\!\!\!\!\!\!\!\psi_{\Delta_1}(d-d_0)\, \Big|\!\!\!\! \sum_{\substack{h_1, h_2 \\ h_1, h_2 \sim H^*  \\ h_1 v_2 \neq h_2v_1 }}  \!\sum_{\substack{n }}^{\,\,*_2}\,C(n)\varphi_H^\star(h_1) \overline{\varphi_H^{\star\star}(h_2)}\psi_N(n)  \Psi(n,h_1,h_2,dr_1,u_1,v_1,v_2,q_2) \Big|.
\end{align*}
Here  $\psi_{\Delta_1}(d_1)$ is a coefficient sequence which is smooth at scale $\Delta_1$ and has $\psi_{\Delta_1}(d_1) \geq 0$ for all $d_1 \in \mathbb{R}$ and $\psi_{\Delta_1}(d_1) \geq 1$ for all $d_1 \in [\Delta_1,2\Delta_1]$. We then have 
\begin{align*}
\Upsilon_2'(v_1,v_2) \ll x^{5\varepsilon} \sup \Upsilon_3(v_1,v_2,d_0).
\end{align*}
Now observe that $\Upsilon_3(v_1,v_2,d_0)$ is of the form $\Sigma_2(z_1,z_2)$ for some $z_2 \in Z_2(z_1)$. Hence, by the assumptions of Lemma~\ref{lem:preparationII}, we have  $ \Upsilon_3(v_1,v_2,d_0) \ll_{\varepsilon} x^{-10\varepsilon}q_0(q_0,\ell)\Delta N  (v_1,v_2)$. This in turn implies  $\Upsilon_2'(v_1,v_2) \ll_{\varepsilon} x^{-5\varepsilon}q_0(q_0,\ell)\Delta N (v_1,v_2) $. Summing over  $v_1, v_2 \asymp V$, we then obtain the desired bound
\begin{align*}
\Upsilon_2 \ll_{\varepsilon}  x^{-5\varepsilon}q_0(q_0,\ell)\Delta N  \sum_{v_1 \asymp V} \sum_{v_2 \asymp V}(v_1,v_2)\ll x^{-4.5\varepsilon}q_0(q_0,\ell)\Delta N  V^2.
\end{align*} 
Hence $\Sigma_2(z_1,z_2) = O_{\varepsilon}(\frac{q_0(q_0,\ell)\Delta N (v_1,v_2)}{x^{10\varepsilon}})$ indeed implies $\Sigma_1(z_1) = O_{\varepsilon}(\frac{(q_0,\ell)RQNUV^2}{x^{4\varepsilon}})$, as proposed.
\end{proof}

\subsubsection{Simplifying exponentials} We now inspect $\Psi(n,h_1,h_2,dr_1,u_1,v_1,v_2,q_2)$ more closely. 

\begin{lemma}[Factorising $\Psi$] \label{lem:simplifyingexponentialsII}
Consider $n, d, h_1,h_2,\ell,r_1, q_0, u_1,v_1,v_2,q_2,a,b_1,b_2 \in \mathbb{Z}$. Assume that 
\begin{align*}
r_1q_0u_1[v_1,v_2]q_2 \mid P(x) \quad \mbox{ and } \quad (r_1q_0u_1v_1v_2q_2,dab_1b_2)= (n,dr_1q_0u_1v_1v_2)=(n+\ell d r_1,q_0q_2)=1.
\end{align*} Let $\Psi(n,h_1,h_2,dr_1,u_1,v_1,v_2,q_2)$ be as given in Definition~{\rm\ref{def:Z1functions}}.   

\vspace{3mm}
Write $v^*=(v_1,v_2)$. Set $v_1^* = v_1/v^*$ and $v_2^* = v_2/ v^* $.  Then
\begin{align*}
\Psi(n,h_1,h_2,dr_1,u_1,v_1,v_2,q_2)=e_{d}\left(\dfrac{a(h_1v_2^*-h_2v_1^*)}{nr_1q_0u_1v^*v_1^*v_2^*q_2}\right) e_{m}\left(\dfrac{A (h_1v_2^*-h_2v_1^*) }{d(n+B d)}\right), 
\end{align*}
 where  $m= r_1q_0 u_1 v^*v_1^*v_2^*q_2= r_1q_0 u_1 [v_1,v_2]q_2$   and where $A$ and $B$ are residue classes mod $m$ which are  independent of $n$, $d$, $h_1$ and $h_2$. Furthermore, $A$ satisfies $(A,m)=1$. 
\end{lemma}

\begin{proof}
Using the Chinese Remainder Theorem as stated in (\ref{CRT}), we find the following: 
  \begin{align*}
 \Psi(n,h_1,h_2, dr_1, \dots)  
&=e_{dr_1}\left(\dfrac{ah_1}{nq_0u_1v_1q_2}-\dfrac{ah_2}{nq_0u_1v_2q_2}\right)e_{q_0u_1v_1}\left(\dfrac{b_1h_1}{ndr_1q_2}\right)e_{q_0u_1v_2}\left(-\dfrac{b_1h_2}{ndr_1q_2}\right) \\
  &e_{q_2}\left(\dfrac{b_2h_1}{(n+\ell dr_1)dr_1q_0u_1v_1}-\dfrac{b_2h_2}{(n+\ell dr_1)dr_1q_0u_1v_2}\right)\\
 &= e_{d}\left(\dfrac{a(h_1v_2-h_2v_1)}{nr_1q_0u_1v_1v_2q_2}\right)e_{r_1}\left(\dfrac{a(h_1v_2-h_2v_1)}{ndq_0u_1v_1v_2q_2}\right)e_{v^*}\left(\dfrac{b_1(h_1v_2^*-h_2v_1^*)}{ndr_1q_0u_1v_1^*v_2^*q_2}\right)\\
 &e_{q_0u_1v_1^*}\left(\dfrac{b_1h_1}{ndr_1v^*q_2}\right)e_{q_0u_1v_2^*}\left(-\dfrac{b_1h_2}{ndr_1v^*q_2}\right)e_{q_2}\left(\dfrac{b_2(h_1v_2-h_2v_1)}{(n+\ell dr_1)dr_1q_0u_1v_1v_2}\right)\\ 
 &= e_{d}\left(\dfrac{a(h_1v_2-h_2v_1)}{nr_1q_0u_1v_1v_2q_2}\right)e_{r_1}\left(\dfrac{a(h_1v_2^*-h_2v_1^*)}{ndq_0u_1v^*v_1^*v_2^*q_2}\right)e_{v^*}\left(\dfrac{b_1(h_1v_2^*-h_2v_1^*)}{ndr_1q_0u_1v_1^*v_2^*q_2}\right)\\
 &e_{q_0u_1}\left(\dfrac{b_1(h_1v_2^*-h_2v_1^*)}{ndr_1v^*v_1^*v_2^*q_2}\right)e_{v_1^*}\left(\dfrac{b_1(h_1v_2^*-h_2v_1^*)}{ndr_1q_0u_1v^*v_2^*q_2}\right)e_{v_2^*}\left(\dfrac{b_1(h_1v_2^*-h_2v_1^*)}{ndr_1q_0u_1v^*v_1^*q_2}\right) \\
 &e_{q_2}\left(\dfrac{b_2(h_1v_2^*-h_2v_1^*)}{(n+\ell dr_1)dr_1q_0u_1v^*v_1^*v_2^*}\right)\\ 
  &= e_{d}\left(\dfrac{a(h_1v_2^*-h_2v_1^*)}{nr_1q_0u_1v^*v_1^*v_2^*q_2}\right) e_{m}\left(\dfrac{A (h_1v_2^*-h_2v_1^*) }{d(n+B d)}\right) 
 \end{align*}
 where  $m= r_1q_0u_1v^*v_1^*v_2^*q_2 = r_1q_0 u_1 [v_1,v_2]q_2$   and where $A, B \in \{1,\dots,m\}$ are integers which satisfy:
\begin{align*}
&A \equiv a \mbox{ mod } r_1,   &&B \equiv 0 \mbox{ mod } r_1, \\
&A \equiv b_1 \mbox{ mod } q_0 u_1 v^*v_1^*v_2^*,   &&B \equiv 0 \mbox{ mod } q_0 u_1 v^*v_1^*v_2^*, \\
&A \equiv b_2 \mbox{ mod } q_2, 
&&B \equiv \ell r_1 \mbox{ mod } q_2. 
\end{align*}
Such $A$ and $B$ exist by the Chinese Remainder Theorem. From the conditions above, we see that $A$ and $B$  are independent of $n$, $d$, $h_1$ and $h_2$. Further, since $(ab_1b_2,r_1q_0u_1v_1v_2q_2)=1$, we have that $a$, $b_1$ and $b_2$ are coprime to $m$ and   $A$ satisfies $ (A,m)=1$.
\end{proof}

\subsubsection{$q$-van der Corput process}

We now apply the $q$-van der Corput method. This is where our arguments diverge from the work of Polymath~\cite{Polymath:2014:EDZ}: On page 2185 of~\cite{Polymath:2014:EDZ}, sums were partitioned according to the value of $n$ mod $d$. Instead, we now partition according to the value of $(h_1v_2^*-h_2v_1^*)/n$ mod $d$. By including both $h$ and $d$ in the innermost sum of the $q$-van der Corput method, the contribution of diagonal terms is decreased. This allows for a larger choice of $D$: In~\cite{Polymath:2014:EDZ}, the proof used $D \approx N/H^4$, but now we may work with $D \approx N/H^2$, which leads to a smaller modulus $m$ in the exponential sums. This change  is what allows for an improvement of the equidistribution estimates.

\begin{defn}\label{def:Z3}
For given $\omega, \varepsilon, \delta, M, N$ and given $x>1$, $z_1 \in Z_1(x)$ and $z_2 \in Z_2(z_1;x)$, the expression $Z_3(z_1,z_2;x)$ denotes the set of tuples  
$$(W_1,Y, w_0,w_1,w_2,A,B, C(n))$$
which satisfy the following five requirements:
\begin{enumerate}[{\rm (i)}]
\item  $W_1 \in (0,\infty)$ with $1 \ll W_1 \ll \Delta$.
\item $Y \in \mathbb{R}$  with $1 \ll |Y| \ll \frac{H^*V}{(v_1,v_2)}$.
 \item $w_0, w_1, w_2$ are positive integers. For $m = r_1q_0u_1[v_1,v_2]q_2$, they satisfy:
\begin{enumerate}[{\rm (a)}]
\item $w_1 \asymp W_1$, $w_1$ is  squarefree and $(w_1,m)=1$.
\item  $w_0 \mid w_1$ and $w_2 \mid m^*$, where  $ m^*=m^{\lfloor 2\log(x) \rfloor}$.
\end{enumerate} 
\item $A$ and $B$ are residue classes mod $m$ with $(A,m)=1$, where $m = r_1q_0u_1[v_1,v_2]q_2$.
\item $C(n)$ is the indicator function of a set $E_d$ mod $q_0$. \\
This set $E_d$  only depends on  the value of $d$ mod $q_0$ and has at most $(q_0,\ell)$ elements.
\end{enumerate}

\smallskip
Additionally, for a given $z_2 \in Z_2(z_1;x)$, we set 
$$\mathcal{L} = \{h_1 v_2^* - h_2 v_1^*: h_1, h_2 \sim H^*\},$$ where $v_1^*= v_1/(v_1,v_2)$ and $v_2^* = v_2/(v_1,v_2)$.
\end{defn}

\begin{lemma}[Alternative $q$-van der Corput method]\label{lem:qvanderCorputII}
Let $\Sigma_2$, $Z_1$, $Z_2$ and  $Z_3$ be as described in Lemma~{\rm \ref{lem:preparationII}}, Definition~{\rm \ref{def:Z1}}, Definition~~{\rm \ref{def:Z2}} and Definition~~{\rm \ref{def:Z3}}.
For given $x>1$, $z_1 \in Z_1(x)$, $z_2\in Z_2(z_1;x)$ and $z_3 \in Z_3(z_1,z_2;x)$,  we define  $\Sigma_3(z_1,z_2,z_3;x)$ as follows:
\begin{align*}
\Sigma_3(\, \cdot \,)= \!\!\!\!\!\! \sum_{\substack{y, \widetilde{y} \in \mathcal{L} \\ y, \widetilde{y} \sim Y \\ w_1 \mid (y, \widetilde{y}) \\ (y, m^*)=w_2 \\ (\widetilde{y}, m^*)=w_2}}\!\!\!\!\Bigg|\sum_{\substack{d \\ w_0w_1 \mid d  \\  (\frac{d}{w_1}, \frac{my\widetilde{y}}{w_1^2})=1}} \!\!\!\!\!\!\!\psi_{\Delta_1}(d-d_0) \!\!\!\!\!\!\!\!\!\!\sum_{\substack{n, \widetilde{n} \\ y\widetilde{n} \equiv 
\widetilde{y}n  \mbox{ \scriptsize mod } d}}^{\,\,*_3} \!\!\!\!\!\!\!\!\!\! C(n)  C(\widetilde{n})\psi_N(n)\psi_N(\widetilde{n}) e_{m}\left(\dfrac{A (y(\widetilde{n}+B d)-  \widetilde{y}(n+B d))}{d(n+B d)(\widetilde{n}+B d)}\right)\Bigg|,
\end{align*}
where $\sum\limits^{\,\,*_3}$ denotes summation over $n$ and $\widetilde{n}$ with $ (n\widetilde{n},w_1 r_1 q_0 u_1 v_1 v_2)=((n+\ell dr_1)(\widetilde{n}+\ell d r_1), q_0q_2)=1  $.

\medskip
We then have the bound
$$\Sigma_2(z_1,z_2;x) \ll_\varepsilon \sup_{z_3 \in Z_3(z_1,z_2)} (x^{3.5\varepsilon}(v_1,v_2)  \Delta) \, \Sigma_3(z_1,z_2,z_3;x)^{1/2}.$$
\end{lemma}

Hence if $\Sigma_3(z_1,z_2,z_3;x) = O_{\varepsilon}\!\left(\dfrac{q_0^2(q_0,\ell)^2 N^2}{x^{27\varepsilon}}\right) $ for $z_3 \in Z_3$,   then   $\Sigma_2(z_1,z_2;x) = O_{\varepsilon}\!\left(\dfrac{q_0(q_0,\ell)\Delta N (v_1,v_2)}{x^{10\varepsilon}}\right)$.

\begin{proof}
Our goal is to bound $\Sigma_2(z_1,z_2)$, which is given by 
\begin{align*}
 \sum_{\substack{d  \\ (d,r_1q_0u_1[v_1,v_2]q_2)=1 \\ d \mbox{\begin{scriptsize}
squarefree
\end{scriptsize}}
}} \!\!\!\!\!\!\!\!\!\!\!\!\!\!\!\!\!\psi_{\Delta_1}(d-d_0)\, \Big|\!\!\!\! \sum_{\substack{h_1, h_2 \\ h_1, h_2 \sim H^*  \\ h_1 v_2 \neq h_2v_1 }}  \sum_{\substack{n }}^{\,\,*_2}\,C(n)\varphi^{\star}_H(h_1) \overline{\varphi^{\star\star}_H(h_2)}\psi_N(n)  \Psi(n,h_1,h_2,dr_1,u_1,v_1,v_2,q_2) \Big|,
\end{align*}
where $\sum\limits^{\,\,*_2}$ indicates that the sum is restricted to $n$ with  $(n,dr_1q_0u_1v_1v_2  )=(n+\ell dr_1, q_0q_2)=1$. 
 Here $C(n) = 1_{\substack{\frac{b_1}{n}  \equiv \frac{b_2}{n+\ell d r_1} (q_0)}}
$, and this is the indicator function of a set $E_d$ mod $q_0$ (dependent on  $d$ mod $q_0$). 

\vspace{3mm}
 We begin our proof by splitting up the sums inside the absolute value signs of $\Sigma_2(z_1,z_2)$  according to the value of $(h_1v_2^*-h_2v_1^*)/n$ mod $d$, where  $v^*=(v_1,v_2)$, $v_1^* = v_1/v^*$ and $v_2^* = v_2/v^*$. We recall Lemma~\ref{lem:simplifyingexponentialsII}.
If $(h_1v_2^*-h_2v_1^*)/n \equiv \gamma$ mod $d$, then 
\begin{align*}
\Psi(n,h_1,h_2,dr_1,u_1,v_1,v_2,q_2)&=e_{d}\left(\dfrac{a(h_1v_2^*-h_2v_1^*)}{nr_1q_0u_1v^*v_1^*v_2^*q_2}\right) e_{m}\left(\dfrac{A (h_1v_2^*-h_2v_1^*) }{d(n+B d)}\right) \\ 
&=e_{d}\left(\dfrac{a\gamma}{r_1q_0u_1v^*v_1^*v_2^*q_2}\right) e_{m}\left(\dfrac{A (h_1v_2^*-h_2v_1^*) }{d(n+B d)}\right),
\end{align*}
 where  $m=  r_1q_0 u_1 [v_1,v_2]q_2$ and $A$ and $B$ are residue classes mod $m$,   independent of $n$, $d$, $h_1$ and $h_2$,  with $(A,m)=1$. 
The $d$-factor no longer depends on $h_1$, $h_2$ and $n$. So we find that: 
 \begin{align*}
\Sigma_2(z_1,z_2)&\ll\!\!\!\!\!\!\!\!
 \sum_{\substack{d  \\ (d,m)=1 \\ d \mbox{\begin{scriptsize}
squarefree
\end{scriptsize}}
}} \!\!\!\!\!\!\!\!\psi_{\Delta_1}(d-d_0)\sum_{\substack{\gamma (d) }}\, \Big|\!\!\!\! \sum_{\substack{h_1, h_2 \\ h_1, h_2 \sim H^* \\ h_1 v_2 \neq h_2v_1 }}  \!\sum_{\substack{n  \\\frac{(h_1v_2^*-h_2v_1^*)}{n} \equiv \gamma(d) }}^{\,\,*_2}\!\!\! \!\!\!\!\!\!\!\!\!\!\!\!C(n)\varphi^\star_H(h_1) \overline{\varphi^{\star\star}_H(h_2)}\psi_N(n)   e_{m}\left(\frac{A (h_1v_2^*-h_2v_1^*) }{d(n+B d)}\right) \Big|.
\end{align*}
The condition $(n,d)=1$, which is part of summation restriction $*_2$, could become bothersome at a later stage. To ease its removal, we now partition according to the value of $w_1=(h_1v_2^*-h_2v_1^*,d)$ and will later only insist that $(n,w_1)=1$. Similarly,  choices of  $h_1v_2^*-h_2v_1^*$ which share prime factors with $m$ may cause issues, and so we partition according to the value of $w_2 = (h_1v_2^*-h_2v_1^*, m^{\lfloor 2\log(x) \rfloor})$.  (Since $ 2^{\lfloor 2\log(x) \rfloor} > x$, this means that we partition according to the largest factor of $h_1v_2^*-h_2v_1^*$ which consists solely of powers of primes which divide $m$.) Writing $m^* = m^{\lfloor 2\log(x) \rfloor}$, we obtain:
\begin{align*}
\Sigma_2
&\ll\!\!\!\!\!\!\!\!\!\!
 \sum_{\substack{d \\ (d,m)=1 \\ d \mbox{\begin{scriptsize}
squarefree
\end{scriptsize}}
}} \!\!\!\!\!\!\!\!\psi_{\Delta_1}(d-d_0)\!\!\sum_{\substack{w_1, w_2 \\ w_1 \mid d \\ w_2 \mid m^*}} \sum_{\substack{\gamma (d) }} \Big| \!\!\!\!\!\!\!\!\!\!\!\!\! \!\!\sum_{\substack{h_1, h_2 \\ h_1, h_2 \sim H^*  \\ h_1v_2 \neq h_2 v_1 \\  (h_1v_2^*-h_2v_1^*,d)=w_1 \\  (h_1v_2^*-h_2v_1^*,m^*)=w_2 }}\!\!\!\!\!\!\!\!  \sum_{\substack{n  \\\frac{(h_1v_2^*-h_2v_1^*)}{n} \equiv \gamma(d) }}^{\,\,*_2} \!\!\!\!\!\!\!\!\!\!\!\!\!\!C(n)\varphi_H^\star(h_1) \overline{\varphi_H^{\star\star}(h_2)}\psi_N(n)   e_{m}\left(\dfrac{A (h_1v_2^*-h_2v_1^*) }{d(n+B d)}\right) \Big|.
\end{align*}
Now we count the number of $w_2 \mid m^*$ which have $w_2 \ll x$. Since $|h_1v_2^*-h_2v_1^*| \ll x$, these are the only $w_2$ which contribute a non-zero amount to $\Sigma_2(z_1,z_2)$. Recall that $m$ is squarefree and denote the number of prime factors of $m$ by $k$. Here $k \leq \lfloor 2\log(x) \rfloor$. There is an injection from  the set $\{w_2 \in \mathbb{N}: w_2 \mid m^*, w_2 \ll x \}$ to the set  $I_k=\{ n \in \mathbb{N}: n=p_1^{\alpha_1}\dots p_{k}^{\alpha_{k}}, n\ll x\}$, where $p_1, \dots, p_k$ are the $k$ smallest primes. 
(Denoting the prime factors of $m$ in ascending order by $q_1, \dots, q_k$, we send $w_2 = q_1^{\alpha_1} \dots q_k^{\alpha_k}$ to $p_1^{\alpha_1} \dots p_k^{\alpha_k}$, which is less than or equal to $w_2$.) But $I_k$ is the set of $p_k$-smooth numbers bounded by $x$. Here $p_k = O(k \log(k)) = O_\varepsilon(\log(x)^{1+\varepsilon/2})$ since $k \leq \lfloor 2\log(x) \rfloor$. It is known that  the number of $\log(x)^{1+\varepsilon/2}$-smooth numbers less than $x$ is bounded by $O_\varepsilon(x^{1-1/(1+\varepsilon/2)+o(1)})=O_\varepsilon(x^\varepsilon)$. (See for instance Corollary~1.3 and  bound (1.14) of~\cite{Hildebrand:1993:IWL}.) Hence we only need to consider $O_\varepsilon(x^\varepsilon)$  values of $w_2$. 
 
\vspace{3mm} 
Next we partition according to the size of $w_1$, restricting to $w_1 \asymp W_1$ for some $W_1$ with $1 \ll W_1 \ll \Delta$. We also fix $w_2$. Finally, we  restrict the value of $h_1v_2^*-h_2v_1^*$ to a dyadic interval, writing $h_1v_2^*-h_2v_1^* \sim Y$ for some $Y$ with $1 \ll |Y| \ll \frac{H^*V}{(v_1,v_2)}$. Here we use that $h_1 v_2^* -h_2 v_1^*\neq 0$. (Note that  $h_1v_2^*-h_2v_1^*$ and $Y$ may be negative, but not zero.) Using dyadic decomposition, we define $\Upsilon_0(W_1,Y,w_1,w_2)$ to equal
\begin{align*}
\Upsilon_0
&=\!\!\!\!\!\!\!
\sum_{\substack{d  \\(d,m)=1 \\ d \mbox{\begin{scriptsize}
squarefree
\end{scriptsize}} \\ w_1 \mid d
}} \!\!\!\!\!\!\!\!\!\psi_{\Delta_1}(d-d_0) \sum_{\substack{\gamma (d) }} \Big| \sum_{\substack{h_1, h_2  }}^*\sum_{\substack{n \\\frac{(h_1v_2^*-h_2v_1^*)}{n} \equiv \gamma(d) }}^{\,\,*_2} \!\!\!\!\!\!\!\!\!\!\!\!\!C(n) \varphi^\star_H(h_1) \overline{\varphi^{\star\star}_H(h_2)}\psi_N(n)  e_{m}\left(\dfrac{A (h_1v_2^*-h_2v_1^*) }{d(n+B d)}\right)\! \Big|,
\end{align*}
where  $\sum\limits^*$ denotes the sum over $h_1$ and $h_2$ with $h_1, h_2 \sim H^*$, $h_1v_2^* - h_2 v_1^* \sim Y$,   $(h_1v_2^*-h_2v_1^*,d)=w_1$ and $(h_1v_2^*-h_2v_1^*,m^*)=w_2$.  We now have
\begin{align*}
\Sigma_2(z_1,z_2) \ll_\varepsilon \sup_{W_1,Y,w_1,w_2} (x^{2\varepsilon} W_1) \Upsilon_0(W_1,Y,w_1,w_2).
\end{align*}
Here the supremum is taken over choices of $W_1$ and $Y$ with  $1 \ll W_1 \ll \Delta$ and $1 \ll |Y| \ll \frac{H^* V}{(v_1,v_2)}$ and over  all squarefree $w_1 \asymp W_1$ with $(w_1,m)=1$ and over all $w_2 \mid m^*$.

\vspace{3mm} We retained the condition that $d$  is squarefree up to this point to ensure that we only need to consider squarefree $w_1$. Now we remove the condition, replacing it by $(d/w_1,w_1)=1$. Observe that $\Upsilon_0 \ll \Upsilon_1$ for 
\begin{align*}
\Upsilon_1
&=\!\!\!\!\!\!\!\!\!\!
\sum_{\substack{d  \\(d,m)=1 \\ (d/w_1,w_1)=1 \\ w_1 \mid d
}} \!\!\!\!\!\!\!\psi_{\Delta_1}(d-d_0) \sum_{\substack{\gamma (d) }} \Big| \sum_{\substack{h_1, h_2 }}^* \sum_{\substack{n \\\frac{(h_1v_2^*-h_2v_1^*)}{n} \equiv \gamma(d) }}^{\,\,*_2} \!\!\!\!\!\!\!\!\!\!\!\!\!C(n)\varphi^\star_H(h_1) \overline{\varphi^{\star\star}_H(h_2)}\psi_N(n)   e_{m}\left(\dfrac{A (h_1v_2^*-h_2v_1^*) }{d(n+B d)}\right) \Big|.
\end{align*}
(Here we used that $\psi_{\Delta_1}$ is non-negative and that $w_1 \mid d$ and $d$ squarefree imply $(d/w_1,w_1)=1$.)

\vspace{3mm}
Since we only sum over  $h_1$ and $h_2$ with  $(h_1v_2^*-h_2v_1^*,d)=w_1$ and $n$ with $(n,d)=1$, residue classes $\gamma (d)$ with $(\gamma,w_1) \neq w_1$ do not contribute to $\Upsilon_1$. Hence, for each residue class mod $(d/w_1)$, only one lifted residue class mod $d$ contributes to $\Upsilon_1$. Furthermore, $\frac{h_1v_2^*-h_2v_1^*}{n}$ is invertible mod $(d/w_1)$. We may replace summation over $\gamma (d)$ by summation over $\gamma_1 (d/w_1)$ and invert $\frac{h_1v_2^*-h_2v_1^*}{n}$, as follows:
\begin{align*}
\Upsilon_1
&=\!\!\!\!\!\!\!\!\!\!\!
\sum_{\substack{d  \\(d,m)=1 \\ (d/w_1,w_1)=1 \\ w_1 \mid d 
}} \!\!\!\!\!\!\!\psi_{\Delta_1}(d-d_0)\!\!\!\! \sum_{\substack{\gamma_1 (d/w_1) }}\!\Big|\!  \sum_{\substack{h_1, h_2 }}^* \! \sum_{\substack{n \\ (n,dr_1q_0u_1v_1v_2 )=1 \\ (n+\ell dr_1, q_0q_2)=1 \\\frac{n}{(h_1v_2^*-h_2v_1^*)} \equiv \gamma_1(\frac{d}{w_1}) }} \!\!\!\!\!\!\!\!\!\!\!\!\!C(n)\varphi^\star_H(h_1) \overline{\varphi^{\star\star}_H(h_2)}\psi_N(n)   e_{m}\!\left(\!\dfrac{A (h_1v_2^*-h_2v_1^*) }{d(n+B d)}\!\right)\! \Big|.
\end{align*}
The condition $(n,d)=1$ ensures that only congruence classes $\gamma_1(\frac{d}{w_1})$ with $(\gamma_1,\frac{d}{w_1})\!=\! 1$ contribute to $\Upsilon_1$. If we replace $(n,d)=1$ by $(n,w_1)=1$, we get new summands corresponding to $(\gamma_1,\frac{d}{w_1}) \neq 1$ and this increases the sum. Hence 
 \begin{align*}
\Upsilon_1 &\ll\!\!\!\!\!
\sum_{\substack{d\\ (d,m)=1 \\ (d/w_1,w_1)=1 \\ w_1 \mid d }} \!\!\!\!\!\!\psi_{\Delta_1}(d-d_0) \sum_{\substack{\gamma_1 (d/w_1) }}\! \Big|  \sum_{\substack{h_1, h_2}}^* \!\sum_{\substack{n \\ (n,w_1r_1q_0,u_1v_1v_2 )=1 \\ (n+\ell dr_1, q_0q_2)=1  \\\frac{n}{(h_1v_2^*-h_2v_1^*)} \equiv \gamma_1(\frac{d}{w_1}) }} \!\!\!\!\!\!\!\!\!\!\!\!\!\!\!\!C(n) \varphi_H^\star(h_1) \overline{\varphi_H^{\star\star}(h_2)}\psi_N(n)  e_{m}\!\left(\!\dfrac{A (h_1v_2^*-h_2v_1^*) }{d(n+B d)}\!\right) \!\Big|.
\end{align*}

Next we apply the Cauchy-Schwarz inequality. We have $\Upsilon_0^2\ll\Upsilon_1^2 \ll_\varepsilon \Upsilon_2 \Upsilon_3$, where  
\begin{align*}
&\Upsilon_2 = x^{\varepsilon}(\Delta/W_1)^2,\\
&\Upsilon_3  =\!\!\!\!\!\! \!\!
\sum_{\substack{d \\ (d,m)=1 \\ (d/w_1,w_1)=1 \\ w_1 \mid d  }}\!\!\!\!\!\!\!\!\!\!\psi_{\Delta_1}(d-d_0)\sum_{\substack{ \gamma_1 (\frac{d}{w_1}) }} \sum_{\substack{h_1, h_2 }}^* \sum_{\substack{\widetilde{h}_1,  \widetilde{h}_2 }}^*  \!
\sum^{\,\,*_3}_{\substack{n, \widetilde{n}  \\ \frac{n} {(h_1v_2^*-h_2v_1^*)}\equiv  
 \gamma_1(\frac{d}{w_1}) 
 \\  
\frac{ \widetilde{n}}{(\widetilde{h}_1v_2^*-\widetilde{h}_2v_1^*)} \equiv \gamma_1(\frac{d}{w_1}) }}\!\!\!\!\!\!\!\!\!\!\!\varphi_H^{\star}(h_1) \overline{\varphi^{\star\star}_H(h_2)}\overline{\varphi_H^{\star}(\widetilde{h}_1)} \varphi_H^{\star\star}(\widetilde{h}_2) \,\,  (\star)   \nonumber \\
&\mbox{with } (\star) = C(n)  C(\widetilde{n})\psi_N(n)\psi_N(\widetilde{n})e_{m}\left(\dfrac{A (h_1v_2^*-h_2v_1^*) }{d(n+B d)}\right)e_{m}\left(-\dfrac{A (\widetilde{h}_1v_2^*-\widetilde{h}_2v_1^*) }{d(\widetilde{n}+B d)}\right).
\end{align*}
In particular, we now have the upper bound
\begin{align*}
\Sigma_2(z_1,z_2) \ll_\varepsilon \sup_{W_1,Y,w_1,w_2} (x^{2.5\varepsilon}  \Delta)  \Upsilon_3(W_1,Y,w_1,w_2)^{1/2}.
\end{align*}

Since  $(h_1v_2^*-h_2v_1^*,d)=w_1$ and $(\widetilde{h}_1v_2^*-\widetilde{h}_2v_1^*,d)=w_1$ and $(w_1,\frac{d}{w_1})=1$, the conditions $\frac{ n}{(h_1v_2^*-h_2v_1^*)} \equiv \frac{ \widetilde{n}}{(\widetilde{h}_1v_2^*-\widetilde{h}_2v_1^*)}  \equiv \gamma_1(\frac{d}{w_1})$ together with the sum over  $\gamma_1 (\frac{d}{w_1})$ reduce to 
\begin{align*}
(h_1v_2^*-h_2v_1^*)\widetilde{n} \equiv (\widetilde{h}_1v_2^*-\widetilde{h}_2v_1^*)n \mbox{ mod } d.
\end{align*}

\vspace{3mm}
Note also that  $|\varphi_H^\star(h)|, |\varphi_H^{\star\star}(h)| \ll  \log(M)^{O(1)}$ for $|h| \ll H$. 
Hence we have the following upper bound: 
\begin{align*}
&\Upsilon_3 \ll x^{o(1)} \!\!\!\!\!\!\!\!\!\sum_{\substack{h_1, \widetilde{h}_1, h_2, \widetilde{h}_2 \\ h_1,\widetilde{h}_1, h_2,\widetilde{h}_2 \sim H^* \\ h_1v_2^*-h_2v_1^*,\, \widetilde{h}_1v_2^*-\widetilde{h}_2v_1^* \sim Y  \\ w_1 \mid h_1v_2^*-h_2v_1^*, \, w_1 \mid \widetilde{h}_1v_2^*-\widetilde{h}_2v_1^* \\ (h_1v_2^*-h_2v_1^*,m^*) =w_2\\ (\widetilde{h}_1v_2^*-\widetilde{h}_2v_1^*,m^*)=w_2  }}  \Bigg|\sum_{\substack{d  \\  (d,m)=1 \\ (d/w_1,w_1)=1 \\ w_1 \mid d \\
(d/w_1, (h_1v_2^*-h_2v_1^*)/w_1)=1 \\ 
(d/w_1, (\widetilde{h}_1v_2^*-\widetilde{h}_2v_1^*)/w_1)=1
}} \!\!\!\!\!\!\!\!\!\!\!\!\psi_{\Delta_1}(d-d_0) \!\!\!\!\!\!\!\!\!\!\sum_{\substack{n, \widetilde{n} \\ (h_1v_2^*-h_2v_1^*)\widetilde{n} \equiv 
(\widetilde{h}_1v_2^*-\widetilde{h}_2v_1^*)n  \mbox{ \scriptsize mod } d }}^{\,\,*_3} \!\!\!\!\!\!\!\!\!\!\!\!\!\!\!\!\!\! (\star) \quad \,\,\quad\Bigg| \quad   \nonumber \\
&\mbox{with } (\star) =C(n)  C(\widetilde{n})\psi_N(n)\psi_N(\widetilde{n}) e_{m}\left(\dfrac{A (h_1v_2^*-h_2v_1^*) }{d(n+B d)}\right)e_{m}\left(-\dfrac{A (\widetilde{h}_1v_2^*-\widetilde{h}_2v_1^*) }{d(\widetilde{n}+B d)}\right).
\end{align*} 
Now consider $y \sim Y$ and count the number of choices of $(h_1,h_2)$ with $h_1 \sim H^*$ and $h_2 \sim H^*$ such that $y = h_1 v_2^* - h_2 v_1^*$. (Recall here that $v_1$ and $v_2$ are already fixed.) Suppose that $(h_1^*,h_2^*)$ is a solution. Since $(v_1^*,v_2^*)=1$, the other possible solutions to $h_1v_2^*-h_2v_1^*=y$ are of the form $(h_1,h_2)=(h_1^*+kv_1^*,h_2^*+kv_2^*)$. However, since $v_1^* \asymp V/(v_1,v_2)$ and since we also require $h_1 \sim H^*$, the number of possible $k$ is bounded by $O(\max\{1, (v_1,v_2)H^*/V\}) = O((v_1,v_2))$. Now write $\mathcal{L}(v_1^*,v_2^*, H^*) = \{h_1 v_2^* - h_2 v_1^*: h_1, h_2 \sim H^*\}$. Set
\begin{align*}
&\Upsilon_4 =\!\!\!\!\!\sum_{\substack{y, \widetilde{y} \in \mathcal{L}  \\ y, \widetilde{y} \sim Y \\ w_1 \mid (y, \widetilde{y}) \\ (y,m^*)=w_2 \\ (\widetilde{y},m^*)=w_2}}\!\Bigg|\sum_{\substack{d \\ w_1 \mid d  \\  (\frac{d}{w_1}, \frac{my\widetilde{y}}{w_1^2})=1 \\ (d/w_1,w_1)=1  
}} \!\!\!\!\!\!\psi_{\Delta_1}(d-d_0)\!\!\!\!\!\sum_{\substack{n, \widetilde{n}\\ y\widetilde{n} \equiv 
\widetilde{y}n  \mbox{ \scriptsize mod } d }}^{\,\,*_3} \!\!\!\!\!\!\! C(n)  C(\widetilde{n})\psi_N(n)\psi_N(\widetilde{n}) e_{m}\left(\dfrac{A (y(\widetilde{n}+B d)-  \widetilde{y}(n+B d))}{d(n+B d)(\widetilde{n}+B d)}\right)\Bigg|.
\end{align*} 
We then have $\Upsilon_3 \ll_\varepsilon (x^\varepsilon(v_1,v_2)^2 )\Upsilon_4$. In particular,
\begin{align*}
\Sigma_2(z_1,z_2) \ll \sup_{W_1,Y,w_1,w_2} (x^{3\varepsilon}(v_1,v_2)\Delta) \, \Upsilon_4(W_1,Y,w_1,w_2)^{1/2}.
\end{align*}
 As a final step, we  remove the condition that $(d/w_1,w_1)=1$: By Möbius inversion,  $\Upsilon_4$ is bounded by 
\begin{align*}
\sum_{\substack{w_0 \\ w_0 \mid w_1 }} \!\!|\mu(w_0)|\!\!\!\!\!\!\!\sum_{\substack{y, \widetilde{y} \in \mathcal{L} \\ y, \widetilde{y} \sim Y \\ w_1 \mid (y, \widetilde{y}) \\ (y,m^*)=w_2 \\ (\widetilde{y},m^*)=w_2}}\!\! \Bigg|\!\!\sum_{\substack{d \\ w_0w_1 \mid d  \\  (\frac{d}{w_1}, \frac{my\widetilde{y}}{w_1^2})=1}}   \!\!\!\!\!\!\!\!\!\!\psi_{\Delta_1}(d-d_0)\!\!\!\!\!\!\sum_{\substack{n, \widetilde{n} \\ y\widetilde{n} \equiv 
\widetilde{y}n  \mbox{ \scriptsize mod } d }}^{\,\,*_3} \!\!\!\!\!\!\! C(n)  C(\widetilde{n}) \psi_N(n)\psi_N(\widetilde{n})e_{m}\left(\dfrac{A (y(\widetilde{n}+B d)-  \widetilde{y}(n+B d))}{d(n+B d)(\widetilde{n}+B d)}\right)\Bigg|.
\end{align*}
But $w_1$ has only $O_\varepsilon(x^\varepsilon)$ divisors. For a given $w_0 \mid w_1$, define
\begin{align*}
\Upsilon_5(w_0)=\!\!\!\!\sum_{\substack{y, \widetilde{y} \in \mathcal{L} \\ y, \widetilde{y} \sim Y \\ w_1 \mid (y, \widetilde{y}) \\ (y,m^*)=w_2 \\ (\widetilde{y},m^*)=w_2}}\Bigg|\sum_{\substack{d \\ w_0w_1 \mid d  \\  (\frac{d}{w_1}, \frac{my\widetilde{y}}{w_1^2})=1}} \!\!\!\!\!\!\!\psi_{\Delta_1}(d-d_0) \!\!\!\!\!\!\!\!\sum_{\substack{n, \widetilde{n} \\ y\widetilde{n} \equiv 
\widetilde{y}n  \mbox{ \scriptsize mod } d }}^{\,\,*_3} \!\!\!\!\!\!\!\!\! C(n)  C(\widetilde{n})\psi_N(n)\psi_N(\widetilde{n}) e_{m}\left(\dfrac{A (y(\widetilde{n}+B d)-  \widetilde{y}(n+B d))}{d(n+B d)(\widetilde{n}+B d)}\right)\Bigg|.
\end{align*}
Observe that we now have the bound \begin{align*}
\Sigma_2(z_1,z_2) \ll \sup_{W_1,Y,w_1,w_2,w_0} (x^{3.5\varepsilon}(v_1,v_2) \Delta) \, \Upsilon_5(W_1,Y,w_1,w_2,w_0)^{1/2},
\end{align*}
where the supremum is taken over choices of $w_0$ with $w_0 \mid w_1$. But $\Upsilon_5(W_1,Y,w_1,w_2,w_0)$ is of the same form as $\Sigma_3(z_1,z_2,z_3)$ for some  $z_3 \in Z_3(z_1,z_2)$. This concludes the proof.
\end{proof}

  \subsubsection{Large divisors}
 
We have reduced the problem of obtaining equidistribution estimates to estimating exponential sums of the form $\Sigma_3$ given in Lemma~\ref{lem:qvanderCorputII}. Tuples $(n,\widetilde{n}, y, \widetilde{y})$ with $y(\widetilde{n}+B d)-  \widetilde{y}(n+B d)=0$ contribute a large amount to $\Sigma_3$. We also have no non-trivial estimates for the exponential sums when $(y(\widetilde{n}+B d)-  \widetilde{y}(n+B d))/d$ shares a very large divisor with $m$. In the next lemma we thus remove such bad choices of $(n,\widetilde{n}, y, \widetilde{y})$. Here it is of great importance that $D$ has been chosen sufficiently small. (Due to the involvement of $y$, the diagonal/bad terms here are different from the ones in~\cite{Polymath:2014:EDZ}, and this is what allows us to choose a larger $D$.)

\begin{lemma}[Removing diagonal terms]\label{lem:diagonalII}
Let $\Sigma_3$, $Z_1$, $Z_2$ and  $Z_3$ be as described in Lemma~{\rm \ref{lem:qvanderCorputII}}, Definition~{\rm \ref{def:Z1}}, Definition~~{\rm \ref{def:Z2}} and Definition~~{\rm \ref{def:Z3}}.
 For given $x>1$, $z_1 \in Z_1(x)$, $z_2 \in Z_2(z_1;x)$ and $z_3 \in Z_3(z_1,z_2;x)$, we denote by $Z_4(z_1,z_2,z_3;x)$ the set of tuples  
$
(y, \widetilde{y})
$ which satisfy the following three properties: 
\begin{enumerate}[{\rm(i)}]
\item $y$ and $\widetilde{y}$ are non-zero integers with $y, \widetilde{y} \sim Y$.
\item $w_1 \mid y$ and $w_1 \mid \widetilde{y}$.
\item$(y,m^*)=(\widetilde{y},m^*)=w_2$, where 
 $m = r_1q_0u_1[v_1,v_2]q_2$ and $m^* = m^{\lfloor 2 \log(x) \rfloor}$.
\end{enumerate}
Further, we set $T(z_1,z_2,z_3;x)= \max\{\tfrac{1}{(w_2,m)}, \tfrac{1}{H}\} \tfrac{x^{\delta+100\varepsilon}H^2 N}{(v_1,v_2)\Delta_1}$ and define $\Sigma_4(z_1,z_2,z_3,z_4;x)$ as follows:
 \begin{align*}
&\Sigma_4(\, \cdot \,)= \Bigg|\!\!\!\!\sum_{\substack{d \\ w_0w_1 \mid d  \\  (\frac{d}{w_1}, \frac{my\widetilde{y}}{w_1^2})=1}}\!\!\!\!\!\!\!\psi_{\Delta_1}(d-d_0) \!\!\!\!\!\!\!\!\!\!\!\!\!\!\!\!\!\!\!\!\!\!\sum_{\substack{n, \widetilde{n}   \\ y\widetilde{n} \equiv 
\widetilde{y}n  \mbox{ \scriptsize mod } d \\ (\frac{y (\widetilde{n} +Bd)-\widetilde{y}(n+Bd)}{d},m ) \leq T(z_1,z_2,z_3)}}^{\,\,*_3}  \!\!\!\!\!\!\!\!\!\!\!\!\!\!\!\!\!\!\!\!\!\!\!\!\!\! C(n)  C(\widetilde{n})\psi_N(n)\psi_N(\widetilde{n}) e_{m}\left(\dfrac{A (y(\widetilde{n}+B d)-  \widetilde{y}(n+B d))}{d(n+B d)(\widetilde{n}+B d)}\right)\Bigg|,
\end{align*}
where $\sum\limits^{\,\,*_3}$ denotes summation over $n$ and $\widetilde{n}$ with $ (n\widetilde{n},w_1 r_1 q_0 u_1 v_1 v_2)=((n+\ell dr_1)(\widetilde{n}+\ell d r_1), q_0q_2)=1  $.

\vspace{3mm}
We then have the bound $$\Sigma_3(z_1,z_2,z_3)  \ll_\varepsilon  \sup_{z_4 \in Z_4(z_1,z_2,z_3)}  \left(\max\left\{\dfrac{x^{\delta+10\varepsilon}H^3}{(w_2,m)}, \dfrac{H^4}{ (w_2,m)}\right\}\right) \Sigma_4(z_1,z_2,z_3,z_4) + (x^{-47\varepsilon}   N^2).$$
\end{lemma}
In particular, if $\Sigma_4(z_1,z_2,z_3,z_4;x) = O_{\varepsilon}\!\left(\min\left\{\dfrac{ (w_2,m)}{x^{\delta+10\varepsilon}H^3}, \dfrac{ (w_2,m)}{H^4}\right\}\dfrac{q_0^2(q_0,\ell)^2 N^2}{x^{27\varepsilon}}\right)$ for $z_4\in Z_4$, then $\Sigma_3(z_1,z_2,z_3;x) = O_{\varepsilon}\!\left(\dfrac{q_0^2(q_0,\ell)^2 N^2}{x^{27\varepsilon}}\right)$.

\begin{proof}
Our goal is to bound $\Sigma_3(z_1,z_2,z_3)$, which is given by
\begin{align*}
\Sigma_3= \!\!\!\!\!\! \sum_{\substack{y, \widetilde{y} \in \mathcal{L} \\ y, \widetilde{y} \sim Y \\ w_1 \mid (y, \widetilde{y}) \\ (y, m^*)=w_2 \\ (\widetilde{y}, m^*)=w_2}}\!\!\!\!\Bigg|\sum_{\substack{d \\ w_0w_1 \mid d  \\  (\frac{d}{w_1}, \frac{my\widetilde{y}}{w_1^2})=1}} \!\!\!\!\!\!\!\psi_{\Delta_1}(d-d_0) \!\!\!\!\!\!\!\!\!\!\sum_{\substack{n, \widetilde{n} \\ y\widetilde{n} \equiv 
\widetilde{y}n  \mbox{ \scriptsize mod } d}}^{\,\,*_3} \!\!\!\!\!\!\!\!\!\! C(n)  C(\widetilde{n})\psi_N(n)\psi_N(\widetilde{n}) e_{m}\left(\dfrac{A (y(\widetilde{n}+B d)-  \widetilde{y}(n+B d))}{d(n+B d)(\widetilde{n}+B d)}\right)\Bigg|.
\end{align*}
 Notice that $d$ divides $(y(\widetilde{n}+B d)-  \widetilde{y}(n+B d))$ since $y\widetilde{n} \equiv 
\widetilde{y}n $ mod $d$. Write 
\begin{align*}
&S(q;y,\widetilde{y},d) = \left\{ (n,\widetilde{n}): \mbox{gcd}\left( \dfrac{y(\widetilde{n}+B d)-  \widetilde{y}(n+B d)}{d}, m \right) = q\right\}, \\
&T^\star(y,\widetilde{y},d) = \bigcup_{\substack{  q \leq \max\left\{\tfrac{1}{(w_2,m)}, \tfrac{1}{H}\right\} \tfrac{x^{\delta+100\varepsilon}H^2 N}{(v_1,v_2)\Delta_1}  }} S(q;y,\widetilde{y},d) .
\end{align*} 
We  partition the set of $(n, \widetilde{n})$, using $T^\star(y,\widetilde{y},d)$ and $S(q;y,\widetilde{y},d)$ with $q >\max\{\tfrac{1}{(w_2,m)}, \tfrac{1}{H}\} \tfrac{x^{\delta+100\varepsilon}H^2 N}{(v_1,v_2)\Delta_1}$. Set 
\begin{align*}
&\Upsilon_T=\sum_{\substack{y, \widetilde{y}\in \mathcal{L} \\ y, \widetilde{y} \sim Y \\ w_1 \mid (y, \widetilde{y}) \\ (y,m^*)=w_2 \\ (\widetilde{y},m^*)=w_2}}\Bigg|\sum_{\substack{d \\ w_0w_1 \mid d  \\  (\frac{d}{w_1}, \frac{my\widetilde{y}}{w_1^2})=1}} \!\!\!\!\!\!\!\psi_{\Delta_1}(d-d_0) \!\!\!\!\!\!\!\!\sum_{\substack{n, \widetilde{n} \\  y\widetilde{n} \equiv 
\widetilde{y}n  \mbox{ \scriptsize mod } d \\ (n,\widetilde{n}) \in T^\star(y,\widetilde{y},d) }}^{\,\,*_3} \!\!\!\!\!\!\!\!\! C(n)  C(\widetilde{n})\psi_N(n)\psi_N(\widetilde{n}) e_{m}\left(\dfrac{A (y(\widetilde{n}+B d)-  \widetilde{y}(n+B d))}{d(n+B d)(\widetilde{n}+B d)}\right)\Bigg|,\\
&\Upsilon_S(q)=\sum_{\substack{y, \widetilde{y} \in \mathcal{L} \\ y, \widetilde{y} \sim Y \\ w_1 \mid (y, \widetilde{y}) \\ (y,m^*)=w_2 \\ (\widetilde{y},m^*)=w_2}}\!\!\Bigg|\!\!\sum_{\substack{d \\ w_0w_1 \mid d  \\  (\frac{d}{w_1}, \frac{my\widetilde{y}}{w_1^2})=1}} \!\!\!\!\!\!\!\!\!\psi_{\Delta_1}(d-d_0)\!\!\!\!\!\!\!\!\sum_{\substack{n, \widetilde{n}  \\ y\widetilde{n} \equiv 
\widetilde{y}n  \mbox{ \scriptsize mod } d \\ (n,\widetilde{n}) \in S(q;y,\widetilde{y},d)}}^{\,\,*_3} \!\!\!\!\!\!\!\!\!\!\! C(n)  C(\widetilde{n})\psi_N(n)\psi_N(\widetilde{n}) e_{m}\left(\dfrac{A (y(\widetilde{n}+B d)-  \widetilde{y}(n+B d))}{d(n+B d)(\widetilde{n}+B d)}\right)\Bigg|.
\end{align*}
The number of divisors of $m$ is bounded by $O_\varepsilon(x^\varepsilon)$, so $S(q; y, \widetilde{y},d)$ is empty for all but $O_\varepsilon(x^\varepsilon)$ choices of $q$. Therefore we have the upper bound 
$$\Sigma_3(z_1,z_2,z_3) \ll_\varepsilon \Upsilon_T + \sup_{\substack{q>T(z_1,z_2,z_3)}} x^\varepsilon\,\Upsilon_S(q).$$
  Initially, we focus on $\Upsilon_S(q)$. 
 Fix $q >T(z_1,z_2,z_3)=\max\{\tfrac{1}{(w_2,m)}, \tfrac{1}{H}\} \tfrac{x^{\delta+100\varepsilon}H^2 N}{(v_1,v_2)\Delta_1}$.  We wish to count the number of choices of $(y,\widetilde{y},n,\widetilde{n},d)$ which appear in $\Upsilon_S(q)$ with  $(n,\widetilde{n}) \in S(q;y,\widetilde{y},d)$. For use in the computations below, recall  that $H \ll x^{4\omega+\delta+7\varepsilon}$, $\Delta \ll D = \frac{N}{x^{50\varepsilon}H^2}$ and $Y \ll \frac{H^*V}{(v_1,v_2)} \ll \frac{x^{\delta+5\varepsilon}H^2}{(v_1,v_2)}$. Hence
 $$ \frac{Y}{q} \ll  \frac{\min\{(w_2,m),H\}}{x^{150\varepsilon}H^2}.$$

 \vspace{3mm}
 We first count the number of $y \in \mathcal{L}$ with $(w_2,m) \mid y$. Note that $(w_2,m)$ is squarefree since $m$ is squarefree. Write $w_2^\circ = ((w_2,m),v_1^*)$. Suppose $y = h_1v_2^*-h_2v_1^*$ with $h_1, h_2 \sim H^*$. Since $(v_1^*,v_2^*)=1$, we have $(w_2,m) \mid y$ if and only if $( h_1  \equiv 0$ mod $w_2^\circ$ and $h_2  \equiv (v_2^*/v_1^*) h_1$ mod $(w_2,m)/w_2^\circ)$.  Thus there are $O(H^*/w_2^\circ)$ choices for $h_1$, and for each $h_1$, there are $O(\max\{w_2^\circ H^*/(w_2,m),1\})$ choices for $h_2$. In total, the number of $y \in \mathcal{L}$ with $w_2 \mid y$ is thus bounded by $O(\max\{H^2/(w_2,m),H\})$. 

\vspace{3mm}  To count the number of choices of $(y,\widetilde{y},n,\widetilde{n},d)$ which appear in $\Upsilon_S(q)$ with  $(n,\widetilde{n}) \in S(q;y,\widetilde{y},d)$, we now first  consider the case ($y = \widetilde{y}$ and $n = \widetilde{n}$). There are $O(\Delta)$ choices for $d$ with $d \asymp \Delta$, there are $O((H^*)^2)$ choices for $y \sim Y$ with $y \in \mathcal{L} = \{h_1 v_2^*-h_2 v_1^*: h_1, h_2 \sim H^*\}$ and there are $O(N)$ choices for $n$. The total contribution of the case   ($y = \widetilde{y}$ and $n = \widetilde{n}$) is thus bounded by $O( \Delta H^2 N) = O(x^{-50\varepsilon}  N^2)$. 

\vspace{3mm} Next we consider the case $(y = \widetilde{y}$ and $n \neq \widetilde{n})$. 
Then $(n,\widetilde{n}) \in S(q;y,\widetilde{y},d)$ implies $q \mid \frac{1}{d}(y(\widetilde{n}-n))$. Given $d \asymp \Delta$,  $y \sim Y$ with $y \in \mathcal{L}$ and $w_2 \mid y$, and $n \asymp N$, there are  only $O(NY/(\Delta q))  $ suitable choices for $\widetilde{n}-n$, each of which gives only one choice for $\widetilde{n}$. There are $O(\max\{H^2/(w_2,m),H\})$ choices for $y \in \mathcal{L}$ with $w_2 \mid y$,    there are $O(\Delta)$ choices for $d \asymp \Delta$ and there are $O(N)$ choices for $n \asymp N$. 
Overall, this bounds the contribution of the case $y = \widetilde{y}$ by $O(\max\{H^2/(w_2,m),H\}\Delta  N^2  Y/(\Delta q) ) = O(x^{-50\varepsilon}  N^2 ).   $ From now on, we thus assume $y \neq \widetilde{y}$. 

\vspace{3mm} Next we treat the case $n = \widetilde{n}$.  Fix a value of $d \asymp \Delta$. We count the number of  $(n, y,\widetilde{y})$ for which $q \mid \frac{1}{d}(y-\widetilde{y})(n+Bd)$. First we note that $dq$ has $O_\varepsilon(x^\varepsilon)$ divisors and fix some $\gamma \mid dq$. We now assume that $(y-\widetilde{y},dq)=\gamma$.  Then we must have $\frac{dq}{\gamma} \mid (n+Bd)$. There are $O(\frac{Y}{\gamma})$ choices for $y-\widetilde{y}$, $O(\frac{Y}{w_2})$ choices for $y \sim Y$ with $w_2 \mid y$, and $O(\gamma N/(\Delta q)+1)$ choices for $n$. Note that $(\frac{Y}{\gamma})(\frac{\gamma N}{\Delta q}+1) = O(\frac{YN}{\Delta q} + Y )$. Overall, this bounds the contribution of the case $n = \widetilde{n}$ by $O_\varepsilon(x^{\varepsilon} \Delta Y^2 \frac{1}{w_2} (\frac{ N }{\Delta q} +1)) = O_\varepsilon(x^{\delta-50\varepsilon}   N  +x^{2\delta} H^2 N) = O(x^{-50\varepsilon} N^2)$  since $H^2 \ll x^{8\omega+3\delta} \ll x^{-3\delta}N$.  We thus now also assume $n \neq \widetilde{n}$.

\vspace{3mm} Fix some $d$ with $d-d_0 \asymp \Delta_1$, some $\widetilde{n} \asymp N$ and some $y_d \in \mathcal{L}-\mathcal{L}$ with $|y_d| \ll Y$, $y_d \neq 0$  and $w_2 \mid y_d$. Since $\mathcal{L}-\mathcal{L} \subseteq \{h_1 v_2^*-h_2 v_1^*: h_1, h_2 \ll H^*\}$, there are $O(\max\{H^2/(w_2,m),H\})$ choices for $y_d$ and so in total there are $O(\max\{H^2/(w_2,m),H\} \Delta N )$ choices for $(d,\widetilde{n},y_d)$. Consider $y$, $\widetilde{y}$ and  $n$ with $y-\widetilde{y}=y_d$ and $(n,\widetilde{n}) \in S(q;y,\widetilde{y},d)$.  Then $q$ divides $\frac{1}{d}(y(\widetilde{n}+B d)-  \widetilde{y}(n+B d))$ and $(y\widetilde{n}-\widetilde{y}n)+B d y_d = Cdq$ for some $C \in \mathbb{Z}$.  Observe that  $|y\widetilde{n}-\widetilde{y}n| \ll  Y N \ll \frac{1}{(v_1,v_2)}x^{\delta+5\varepsilon}H^2 N $, while $|dq| \gg  \max\{\frac{1}{(w_2,m)}, \frac{1}{H}\}\frac{x^{\delta+100\varepsilon}  H^2 N}{(v_1,v_2)}$.  Hence $C = By_d/q+ O\left(  \max\{\frac{1}{(w_2,m)}, \frac{1}{H}\}^{-1}\right)$ and this bounds the number of possible choices for $C$. Choose one such $C$. We have $ \widetilde{y}(\widetilde{n}-n) = Cdq -(Bd+\widetilde{n})y_d$. All values on the RHS have been fixed, whereas the LHS is assumed to be non-zero. This gives $O_\varepsilon(x^\varepsilon)$ choices for the factors on  the LHS.  The values of $\widetilde{n}$ and $\widetilde{n}-n$ also determine the value of $n$. Overall, this bounds the number of   $(y, \widetilde{y}, n, \widetilde{n},d)$  which appear in $\Upsilon_S(q)$ with $(n,\widetilde{n}) \in S(q;y,\widetilde{y},d)$ by 
$O(x^{\varepsilon}H^2\Delta N)= O(
 x^{-49\varepsilon}   N^2).$   We obtain $$
\Upsilon_S(q) \ll_\varepsilon x^{-48\varepsilon}   N^2. $$

Now we look at $\Upsilon_T$. Consider $y$ and $\widetilde{y}$ with $y, \widetilde{y} \in \mathcal{L}$, $y, \widetilde{y} \sim Y$ and  $w_2 \mid (y,\widetilde{y})$.  As observed earlier, there are $O(\max\{H^4/(w_2,m)^2, H^2\})$ pairs $y, \widetilde{y}$ with $y, \widetilde{y} \in \mathcal{L}$ and  $w_2 \mid (y, \widetilde{y})$. 
On the other hand, the number of pairs $y, \widetilde{y} \sim Y$ with $w_2 \mid (y,\widetilde{y})$ is clearly also bounded by $O(Y^2/ (w_2,m)^2)$. Now observe that $Y^2/(w_2,m)^2 = O(x^{\delta+10\varepsilon}H^3/ (w_2,m))$ when $ (w_2,m) > x^\delta H$. But when $(w_2,m) \leq x^{\delta}H$, we also have $H^2 = O(x^{\delta}H^3/(w_2,m))$.  Overall, the number of relevant $y$ and $\widetilde{y}$ is thus bounded by $O(\max\{\frac{x^{\delta+10\varepsilon}H^3}{ (w_2,m)}, \frac{H^4}{ (w_2,m)}\})$. 
Recall that $\Sigma_4(z_1,z_2,z_3,(y, \widetilde{y}))$ is given by
\begin{align*}
\Bigg|\!\!\sum_{\substack{d \\ w_0w_1 \mid d  \\  (\frac{d}{w_1}, \frac{my\widetilde{y}}{w_1^2})=1}}\!\!\!\!\!\!\!\psi_{\Delta_1}(d-d_0) \!\!\!\!\!\!\!\!\!\!\!\!\!\!\!\!\!\!\!\!\!\!\sum_{\substack{n, \widetilde{n}   \\ y\widetilde{n} \equiv 
\widetilde{y}n  \mbox{ \scriptsize mod } d \\ (\frac{y (\widetilde{n} +Bd)-\widetilde{y}(n+Bd)}{d},m ) \leq T(z_1,z_2,z_3)}}^{\,\,*_3}  \!\!\!\!\!\!\!\!\!\!\!\!\!\!\!\!\!\!\!\!\!\!\!\!\!\! C(n)  C(\widetilde{n})\psi_N(n)\psi_N(\widetilde{n}) e_{m}\left(\dfrac{A (y(\widetilde{n}+B d)-  \widetilde{y}(n+B d))}{d(n+B d)(\widetilde{n}+B d)}\right)\Bigg|.
\end{align*}
We then have $\Upsilon_T \ll (\max\{\frac{x^{\delta+10\varepsilon}H^3}{ (w_2,m)}, \frac{H^4}{ (w_2,m)}\}) \sup_{(y,\widetilde{y})} \Sigma_4(z_1,z_2,z_3,(y, \widetilde{y}))$, where the supremum is taken over all choices of $y$ and $\widetilde{y}$ with $y, \widetilde{y} \in \mathcal{L}$, $y, \widetilde{y} \sim Y$, $w_1\mid (y, \widetilde{y})$ and $(y,m^*)=(\widetilde{y},m^*)=w_2$. Thus
$$\Sigma_3(z_1,z_2,z_3)  \ll_\varepsilon  \sup_{z_4 \in Z_4(z_1,z_2,z_3)}  \left(\max\left\{\frac{x^{\delta+10\varepsilon}H^3}{ (w_2,m)}, \frac{H^4}{ (w_2,m)}\right\}\right) \Sigma_4(z_1,z_2,z_3,z_4) + (x^{-47\varepsilon}  N^2).$$
This completes the proof. 
\end{proof}

\subsection{Summary - Results so far}

By now, we have accumulated a lot of notation. For this reason we now summarize our results up to this point and simplify their presentation. 

\begin{defn} \label{def:Z5}
For  given $\omega,\varepsilon, \delta, M$ and $N$,  the expression $Z_5(x)$ denotes the set of tuples $$(Q,R,\Delta_1,\Lambda, m,q_0,z_1,d_0,v_1,v_2,w_1,w_2,c_1,c_2,\ell,l,\lambda,\widetilde{\lambda},A,B,\psi_N,\psi_{\Delta_1}, C(n;d))$$ which satisfy the following  properties: 
\begin{enumerate}[{\rm (i)}]
\item $Q$, $R$ and $\Delta_1$ are positive real numbers.
\item $\Lambda$ is a non-zero real number.
\item $m,q_0,z_1,d_0,v_1,v_2,w_1,w_2, c_1,c_2$ are positive integers which satisfy:
\begin{enumerate}[{\rm (a)}]
\item $m \mid P(x^\delta)$.
\item $q_0 \ll Q$  and $q_0 \mid m$.
\item $z_1 \ll x^{5\varepsilon}\Delta_1$ and  $(z_1,m)=1$.
\item $d_0 \asymp x^{5\varepsilon}\Delta_1$. 
\item $w_1$ is  squarefree and $w_1 \mid z_1$. 
\item  $w_2 \mid m^*$ for $m^* = m^{\lfloor 2\log(x)\rfloor} $.
\item $c_1 \mid m$ and $c_2 \mid m$.
\end{enumerate}
\item  $\ell, l, \lambda, \widetilde{\lambda} \in \mathbb{Z}\setminus \{0\}$  with $|\ell|\ll \frac{N}{R}$, $|l| \ll N$ and $\lambda, \widetilde{\lambda} \sim \Lambda$ and $(\lambda,m^*) = (\widetilde{\lambda},m^*)=w_2$.
\item $A$ and $B$ are residue classes mod $m$ with $(A,m) =1 $. 
\item $\psi_N$ and $\psi_{\Delta_1}$  are non-negative, real coefficient sequences, smooth at scales $N$ and $\Delta_1$.
\item $C(n;d)=1_{n \in E_d (q_0)}$ is the indicator function of a set $E_d$ mod $q_0$ $($dependent on $d$ mod $q_0)$ which has at most $(q_0,\ell)$ elements. 
\item  Finally, the following 7 bounds are satisfied:
\begin{align*} 
&1 \leq H = x^\varepsilon RQ^2(q_0 M)^{-1}\ll \tfrac{x^{4\omega+\delta+7\varepsilon}}{q_0},\\
&x^{-4\varepsilon-\delta}N \ll R \ll x^{-2\varepsilon}N, \\
&x^{1/2-\varepsilon} \ll QR \ll x^{1/2+2\omega+\varepsilon}, \\ 
&x^{5\varepsilon}H \ll v_1, v_2 \ll x^{\delta+5\varepsilon}H,\\
 &\tfrac{N}{x^{\delta+55\varepsilon}H^2}\ll \Delta_1 \ll \tfrac{N}{x^{55\varepsilon}H^2}, \\
 &\tfrac{RQ^2H}{q_0 (v_1,v_2)\Delta_1} \ll m \ll \tfrac{x^{\delta}RQ^2H}{q_0 (v_1,v_2)\Delta_1}, \\
&1 \ll |\Lambda| \ll \tfrac{1}{w_1(v_1,v_2)}x^{\delta+5\varepsilon}H^2.
\end{align*}
\end{enumerate}
\end{defn}

\begin{lemma}[Summary]\label{lem:summary}
Let $\omega, \delta >0$. Let $\varepsilon \in (0, 10^{-100}\delta)$. 
 Let $\alpha$ and $\beta$ be coefficient sequences at scales $M$ and $N$ with $x \ll M(x)N(x) \ll x$.  Let $N= x^{\gamma(x)}$, where $\gamma(x)  \in (12\omega+6\delta,\frac{1}{2}-2\omega-8\varepsilon)$. Assume $\beta$ has the Siegel-Walfisz property.
 
\vspace{3mm} Let $Z_5(x)$ be as given in Definition~{\rm\ref{def:Z5}}. 
For $z \in Z_5(x)$, we set $T(z;x) = \max\{\tfrac{1}{(w_2,m)}, \tfrac{1}{H}\} \tfrac{x^{\delta+100\varepsilon}H^2 N}{(v_1,v_2)\Delta_1}$ and
 \begin{align*}
\Sigma_5(z;x)= \Bigg|\!\!\!\!\sum_{\substack{d \\   (d,m \lambda \widetilde{\lambda})=1}}\!\!\!\!\!\!\!\psi_{\Delta_1}(z_1d-d_0) \!\!\!\!\!\!\!\!\!\!\!\!\!\!\!\!\!\!\!\!\!\!\!\!\!\!\!\sum_{\substack{n, \widetilde{n} \\ (n\widetilde{n},w_1 c_1)=1 \\ ((n+l d)(\widetilde{n}+l d ), c_2)=1  \\ \lambda \widetilde{n} \equiv 
\widetilde{\lambda}n  \mbox{ \scriptsize mod } (z_1 d/w_1) \\ (\frac{\lambda (\widetilde{n} +Bd)-\widetilde{\lambda}(n+Bd)}{d},m ) \leq T(z;x)}} \!\!\!\!\!\! \!\!\!\!\!\!\!\!\!\!\!\!\!\!\!\!\!\!\!\!\! C(n;d)  C(\widetilde{n};d)\psi_N(n)\psi_N(\widetilde{n}) e_{m}\left(\dfrac{A (\lambda(\widetilde{n}+B d)-  \widetilde{\lambda}(n+B d))}{d(n+B d)(\widetilde{n}+B d)}\right)\Bigg|. 
\end{align*}
If $\sup_{z \in Z_5}\Sigma_5(z;x) = O_\varepsilon\!\left(\min\left\{\dfrac{(w_2,m)}{x^{\delta+10\varepsilon}H^3}, \dfrac{ (w_2,m)}{H^4}\right\}\dfrac{q_0^2(q_0,\ell)^2 N^2 }{x^{27\varepsilon}}\right) $ for all $x>1$ and $z \in Z_5(x)$,  then the following equidistribution estimate holds for  $x>1$, $C >0$ and  $a \in \mathbb{Z}$:
\begin{align*} 
\sum_{\substack{d \leq x^{1/2+2\omega} \\ q \mid P(x^\delta) \\ (q,a)=1}}\Bigg| \sum_{n \equiv a(q)} (\alpha \star \beta)(n;x) - \dfrac{1}{\phi(q)} \sum_{(n,q)=1} (\alpha \star \beta)(n;x)\Bigg| \ll_{C,\varepsilon}  \dfrac{x}{\log(x)^C}.
\end{align*}
\end{lemma}

\begin{proof}
This is simply a summary of Lemma~\ref{lem:summaryofpolymath}, Lemma~\ref{lem:preparationII}, Lemma~\ref{lem:qvanderCorputII} and Lemma~\ref{lem:diagonalII}: 
Recall that in Lemma~\ref{lem:diagonalII}, for given $\alpha_1 \in Z_1$, $\alpha_2 \in Z_2(\alpha_1)$, $\alpha_3 \in Z_3(\alpha_1,\alpha_2)$ and $\alpha_4 \in Z_4(\alpha_1,\alpha_2,\alpha_3)$, we defined the quantity 
\begin{align*}
&\Sigma_4(\alpha_1,\alpha_2,\alpha_3,\alpha_4)= \Bigg|\!\!\!\!\!\!\sum_{\substack{d \\ w_0w_1 \mid d  \\  (\frac{d}{w_1}, \frac{my\widetilde{y}}{w_1^2})=1}}\!\!\!\!\!\!\!\psi_{\Delta_1}(d-d_0) \!\!\!\!\!\!\!\!\!\!\!\!\!\!\!\!\!\!\!\!\!\!\!\!\!\!\!\!\!\!\sum_{\substack{n, \widetilde{n} \\ (n\widetilde{n},w_1 r_1 q_0 u_1 v_1 v_2)=1 \\ ((n+\ell dr_1)(\widetilde{n}+\ell d r_1), q_0q_2)=1  \\ y\widetilde{n} \equiv 
\widetilde{y}n  \mbox{ \scriptsize mod } d \\ (\frac{y (\widetilde{n} +Bd)-\widetilde{y}(n+Bd)}{d},m ) \leq T(\alpha_1,\alpha_2,\alpha_3)}}  \!\!\!\!\!\!\!\!\!\!\!\!\!\!\!\!\!\!\!\!\!\!\!\!\!\!\!\!\!\!\!\! C(n)  C(\widetilde{n})\psi_N(n)\psi_N(\widetilde{n}) e_{m}\!\left(\dfrac{A (y(\widetilde{n}+B d)-  \widetilde{y}(n+B d))}{d(n+B d)(\widetilde{n}+B d)}\right)\!\Bigg|.
\end{align*}
Here $m=r_1q_0u_1[v_1,v_2]q_2=r_1q_0u_1v_1v_2q_2/(v_1,v_2)$ with $r_1 \asymp R/\Delta=R/(x^{5\varepsilon}\Delta_1)$, $u_1v_1, u_1v_2 \asymp Q/q_0$, $x^{5\varepsilon}H \ll v_1,v_2 \ll x^{\delta+5\varepsilon}H$ and $q_2 \asymp Q/q_0$. Hence $\frac{1}{q_0 (v_1,v_2)}RQ^2H/\Delta_1 \ll m \ll \frac{1}{q_0 (v_1,v_2)}x^{\delta}RQ^2H/\Delta_1$. 

\vspace{3mm}
We now   write  $z_1 = w_0w_1$.  We  observe that $w_1 \mid z_1$, $(z_1,m)=1$ and $z_1 \ll x^{5\varepsilon}\Delta_1$ when $\Sigma_4(\alpha_1,\alpha_2,\alpha_3,\alpha_4)$ is non-zero. Recall that $w_1 \mid y$ and $w_1 \mid \widetilde{y}$ and write $y = w_1 \lambda$ and $\widetilde{y} = w_1 \widetilde{\lambda}$. Then $\lambda, \widetilde{\lambda} \sim \Lambda = Y/w_1$, where $|\Lambda| \ll  \tfrac{1}{w_1(v_1,v_2)}x^{\delta+5\varepsilon}H^2$.  Write $d = z_1 e$. 
If $(z_1/w_1,m\lambda\widetilde{\lambda}) \neq 1$, then $\Sigma_4(\alpha_1,\alpha_2,\alpha_3,\alpha_4)$ is zero. Hence we  assume $(z_1/w_1,m\lambda\widetilde{\lambda}) =1$ and replace $(z_1e/w_1, m\lambda\widetilde{\lambda})=1$ by $(e, m\lambda\widetilde{\lambda})=1$ in the first sum in $\Sigma_4(\alpha_1,\alpha_2,\alpha_3,\alpha_4)$. Further, $w_1 \lambda \widetilde{n} \equiv 
w_1 \widetilde{\lambda} n $ mod $z_1 e$ is equivalent to $\lambda \widetilde{n} \equiv 
 \widetilde{\lambda} n $ mod $ (z_1/w_1) e$ since $w_1 \mid z_1$. Finally, since $(z_1,m)=1$ and $w_1 \mid z_1$, we have $(\frac{w_1}{z_1}\frac{ \lambda (\widetilde{n} +Bz_1 e)- \widetilde{\lambda}(n+B z_1 e)}{e},m ) =  (\frac{ \lambda (\widetilde{n} +Bz_1 e)- \widetilde{\lambda}(n+B z_1 e)}{e},m ) $. 

\vspace{3mm}
 So whenever $\Sigma_4(\alpha_1,\alpha_2,\alpha_3,\alpha_4)\neq 0$, there exists some $z \in Z_5$ with $\Sigma_4(\alpha_1,\alpha_2,\alpha_3,\alpha_4) = \Sigma_5(z)$. This choice of $\Sigma_5(z)$ has   $\ell z_1 r_1 $ in the place of $l$, $r_1q_0u_1[v_1,v_2]$ in the place of $c_1$, $q_0q_2$ in the place of $c_2$, $Aw_1/z_1$ mod $m$ in the place of $A$ mod $m$, $Bz_1$ in the place of $B$ and $Y/w_1$  in the place of $\Lambda$. 
 
 \vspace{3mm} We now assume $\sup_{z \in Z_5}\Sigma_5(z) = O_\varepsilon(\min\{\tfrac{(w_2,m)}{x^{\delta+10\varepsilon}H^3}, \tfrac{ (w_2,m)}{H^4}\}\tfrac{q_0^2(q_0,\ell)^2 N^2 }{x^{27\varepsilon}}) $. Then $$\Sigma_4(\alpha_1,\alpha_2,\alpha_3,\alpha_4)\ll\sup_{z \in Z_5} \Sigma_5(z)=O_\varepsilon(\min\{\tfrac{(w_2,m)}{x^{\delta+10\varepsilon}H^3}, \tfrac{ (w_2,m)}{H^4}\}\tfrac{q_0^2(q_0,\ell)^2 N^2 }{x^{27\varepsilon}})$$ for $\alpha_1 \in Z_1$, $\alpha_2 \in Z_2(\alpha_1)$, $\alpha_3 \in Z_3(\alpha_1,\alpha_2)$ and $\alpha_4 \in Z_4(\alpha_1,\alpha_2,\alpha_3)$.  By Lemma~\ref{lem:diagonalII}, this implies  $\Sigma_3(\alpha_1,\alpha_2,\alpha_3) = O_{\varepsilon}(\frac{q_0^2(q_0,\ell)^2 N^2}{x^{27\varepsilon}})$. By Lemma~\ref{lem:qvanderCorputII}, then  $\Sigma_2(\alpha_1,\alpha_2) = O_{\varepsilon}(\frac{q_0(q_0,\ell)\Delta N (v_1,v_2)}{x^{10\varepsilon}})$. By Lemma~\ref{lem:preparationII}, then $\Sigma_1(\alpha_1)=O_{\varepsilon}(\frac{(q_0,\ell)RQNUV^2}{x^{4\varepsilon}})$. Finally, by Lemma~\ref{lem:summaryofpolymath},  
$$
\sum_{\substack{d \leq x^{1/2+2\omega} \\ q \mid P(x^\delta) \\ (q,a)=1}}\Bigg| \sum_{n \equiv a(q)} (\alpha \star \beta)(n;x) - \dfrac{1}{\phi(q)} \sum_{(n,q)=1} (\alpha \star \beta)(n;x)\Bigg| \ll_{C,\varepsilon}  \dfrac{x}{\log(x)^C}. \eqno \qedhere
$$
\end{proof}

\subsection{Further adjustments}\label{ssec:adjustments} 

In Section~\ref{ssec:adjustments} we simplify $\Sigma_5$ further, aiming to replace it by an expression to which the exponential sum bounds of~\cite{Polymath:2014:EDZ}, and in particular a version of Polymath's Proposition~8.4, can be applied. 
 
\subsubsection{Removal of $\widetilde{n}$} As a first step, we remove $\widetilde{n}$. We use $\lambda \widetilde{n} \equiv \widetilde{\lambda} n$ mod $(z_1 d/w_1)$ to replace $\widetilde{n}$ by a linear combination of $d$ and $n$, and then use Taylor expansions to separate variables $n$ and $d$ in $\psi_N(\widetilde{n})$. 

\begin{lemma}[Removing a linear congruence condition] \label{lem:adjustments1} Let $Z_5$ and $\Sigma_5$ be as given in Lemma~{\rm\ref{lem:summary}}. 

For $z \in Z_5(x)$,  $W_6(z;x)$ denotes the set of tuples $ ( \psi_N^*,\psi_{\Delta_1}^*, d_\star, n_\star )$ which satisfy the following conditions:
\begin{enumerate}[{\rm (i)}]
\item $\psi_N^*(n)$ and $\psi_{\Delta_1}^*(d_1)$ are coefficient sequences, smooth at scales $N$ and $\Delta_1$. 
\item $d_\star$ and $n_\star$  are residue classes mod $q_0$.
\end{enumerate}
Further, we denote by $\Xi(z;x)$ the set of sequences $(\xi_k)$ with $\xi_k \in [0, \infty)$ and 
$\!\!\!\!\!\!\sum_{\substack{|k| \ll \frac{w_1\Lambda N}{x^{5\varepsilon}\Delta_1}  
\\ ((\frac{z_1}{w_1})k+(\lambda-\widetilde{\lambda})B,m) \leq T(z;x)\\ w_2 \mid k}}\limits \!\!\!\!\!\!\!\!\!\!\!\!\!\! \xi_k \,\,\,\,\,\ll_\varepsilon 1.$

 For  given $x>1$, $z \in Z_5(x)$, $w \in W_6(z;x)$ and $k \in \mathbb{Z}$ with $|k| \ll  \frac{w_1\Lambda N}{x^{5\varepsilon}\Delta_1}$,  we set
\begin{align*}
\Sigma_6(z,w,k;x)&=  \Bigg| \sum_{\substack{d \\ d \equiv d_\star (q_0) \\  (d, m \lambda \widetilde{\lambda})=1 }} \psi_{\Delta_1}^*(z_1 d-d_0)\sum^{\,\,*_6}_{\substack{n \\  n \equiv n_\star (q_0)  }} \psi_N^*(n)e_{m}\!\!\left(\dfrac{A ((\frac{z_1}{w_1}) k   +  (\lambda-\widetilde{\lambda})B)}{(n+B d)(\frac{1}{\lambda}(\widetilde{\lambda}n+(\frac{z_1}{w_1})kd)+B d)}\right)\Bigg|,
\end{align*}

\vspace{-2mm}
where $\sum\limits^{\,\,*_6}$ denotes summation over $n$ with $ \lambda \mid \widetilde{\lambda} n +k(\frac{z_1}{w_1})d$, $ (n,w_1 c_1)=1,$  $(\frac{1}{\lambda}(\widetilde{\lambda}n+k(\frac{z_1}{w_1})d),w_1 c_1)=1$, $(n+l d, c_2)=1$ and $(\frac{1}{\lambda}(\widetilde{\lambda}n+k(\frac{z_1}{w_1})d)+l d, c_2)=1$.

\vspace{3mm} Suppose that for every $x>1$ and $z \in Z_5(x)$, there exists $(\xi_k) \in \Xi(z;x)$ such that for all $k$ with   $|k| \ll  \frac{w_1\Lambda N}{x^{5\varepsilon}\Delta_1}$,  
 $w_2 \mid k$  and $((\frac{z_1}{w_1})k+(\lambda-\widetilde{\lambda})B,m) \leq T(z;x)$, 
\begin{align*}
\sup_{w \in W_6(z)} \Sigma_6(z,w,k)  \ll_\varepsilon \xi_k  \left(\min\left\{\frac{(w_2,m)}{x^{\delta+10\varepsilon}H^3},\frac{(w_2,m)}{H^4}\right\}\frac{q_0(q_0,\ell)  N^2}{ x^{27\varepsilon}} \right).
\end{align*}
Then   also  $\sup_{z \in Z_5}\Sigma_5(z;x) \ll_\varepsilon  \min\left\{\dfrac{(w_2,m)}{x^{\delta+10\varepsilon}H^3}, \dfrac{ (w_2,m)}{H^4}\right\}\dfrac{q_0^2(q_0,\ell)^2 N^2 }{x^{27\varepsilon}} $.
\end{lemma}

\begin{proof}
Our goal is to show that $\Sigma_5(z) \ll_\varepsilon \min\{\tfrac{(w_2,m)}{x^{\delta+10\varepsilon}H^3}, \tfrac{ (w_2,m)}{H^4}\}\tfrac{q_0^2(q_0,\ell)^2 N^2 }{x^{27\varepsilon}} $.  Recall that 
\begin{align*}
\Sigma_5(z)=\Bigg|\sum_{\substack{d \\   (d,m \lambda \widetilde{\lambda})=1}}\!\!\!\!\!\!\!\psi_{\Delta_1}(z_1d-d_0) \!\!\!\!\!\!\!\!\!\!\!\!\!\!\!\!\!\!\!\!\!\!\!\!\!\!\!\!\!\!\!\sum_{\substack{n, \widetilde{n} \\ (n\widetilde{n},w_1 c_1)=1 \\ ((n+ld)(\widetilde{n}+l d), c_2)=1  \\ \lambda \widetilde{n} \equiv 
\widetilde{\lambda}n  \mbox{ \scriptsize mod } (z_1 d/w_1) \\ ((\lambda (\widetilde{n} +Bd)-\widetilde{\lambda}(n+Bd))/d,m ) \leq T(z)}} \!\!\!\!\!\! \!\!\!\!\!\!\!\!\!\!\!\!\!\!\!\!\!\!\!\!\!\!\!\!\! C(n;d)  C(\widetilde{n};d)\psi_N(n)\psi_N(\widetilde{n}) e_{m}\left(\dfrac{A (\lambda(\widetilde{n}+B d)-  \widetilde{\lambda}(n+B d))}{d(n+B d)(\widetilde{n}+B d)}\right)\Bigg|.
\end{align*}

Consider the condition $\lambda \widetilde{n} \equiv \widetilde{\lambda} n$ mod $(z_1 d/w_1)$. Rearranging, we get 
\begin{align*}
 \widetilde{n} = \dfrac{\widetilde{\lambda} n +k(\frac{z_1}{w_1}) d}{\lambda}
\end{align*}
for some $k \in \mathbb{Z}$ with $\lambda \mid \widetilde{\lambda} n +k(\frac{z_1}{w_1}) d$. Since $\psi_N(n)$ is a coefficient sequence at scale $N$, the contribution of $n$ and $\widetilde{n}$ to $\Sigma_5(z)$ is only non-zero when $1 \leq n, \widetilde{n} \ll N$.   Since $\lambda, \widetilde{\lambda} \sim \Lambda$, then $|\lambda \widetilde{n} - \widetilde{\lambda} n| \ll \Lambda N$. Further, $\psi_{\Delta_1}(d_1)$ is a coefficient sequence at scale $\Delta_1$ and $d_0 \asymp x^{5\varepsilon}\Delta_1$. Hence the contribution of $d$ to $\Sigma_5(z)$ is only  non-zero when $ z_1 d  \asymp x^{5\varepsilon}\Delta_1$. However, $k = (\frac{w_1}{z_1 d} ) (\lambda \widetilde{n}  - \widetilde{\lambda} n)
$ and so $|k| \ll \frac{w_1 \Lambda N}{x^{5\varepsilon}\Delta_1} $ whenever the contribution of $k$ to $\Sigma_5(z)$ is non-zero.
Observe  that now  $(\lambda(\widetilde{n}+B d)-  \widetilde{\lambda}(n+B d)) = ((\frac{z_1}{w_1}) kd   +  (\lambda-\widetilde{\lambda})Bd)$ and 
\begin{align*}
e_{m}\left(\dfrac{A (\lambda(\widetilde{n}+B d)-  \widetilde{\lambda}(n+B d))}{d(n+B d)(\widetilde{n}+B d)}\right) = e_{m}\left(\dfrac{A ((\frac{z_1}{w_1}) k   +  (\lambda-\widetilde{\lambda})B)}{(n+B d)(\frac{1}{\lambda}(\widetilde{\lambda}n+(\frac{z_1}{w_1})kd)+B d)}\right).
\end{align*}
Hence we may write
\begin{align*}
&\Sigma_5(z) =\, \Bigg|\!\!\!\!\!\!\!\!\!\!\!\!\!\!\!\!\!\sum_{\substack{k \\ |k| \ll \frac{w_1\Lambda N}{x^{5\varepsilon}\Delta_1} \\ ((\frac{z_1}{w_1}) k   +  (\lambda-\widetilde{\lambda})B,m) \leq T(z) }}\!\!\!\!\!\!\!\!\!\!\!\!\!\!\sum_{\substack{d \\  (d, m \lambda \widetilde{\lambda})=1}} \!\!\!\!\!\!\!\psi_{\Delta_1}(z_1 d-d_0)\sum_{\substack{n }}^{\,\,*_6} C(n;d)  C\!\left(\dfrac{\widetilde{\lambda} n +k(\frac{z_1}{w_1}) d}{\lambda};d\right)\psi_N(n)\psi_N\!\left(\frac{\widetilde{\lambda} n +k(\frac{z_1}{w_1}) d}{\lambda}\right) (\star)\,\Bigg|\\
&\mbox{with }  (\star ) = e_{m}\left(\dfrac{A ((\frac{z_1}{w_1}) k   +  (\lambda-\widetilde{\lambda})B)}{(n+B d)(\frac{1}{\lambda}(\widetilde{\lambda}n+(\frac{z_1}{w_1})kd)+B d)}\right).
\end{align*}
Since $\psi_N(n)$ is a smooth coefficient sequence, there exists a smooth function $\psi:\mathbb{R} \rightarrow \mathbb{C}$, supported on $[c,C]$, with $|\psi^{(j)}(t)|\ll_j \log(x)^{O_j(1)}$ and $\psi_N(n)=\psi(n/N)$. We now focus on 
\begin{align*}
\psi_N\!\left(\frac{\widetilde{\lambda} n +k(\frac{z_1}{w_1}) d}{\lambda}\right)= \psi\!\left(\frac{\widetilde{\lambda} n +(\frac{k}{w_1}) d_0}{\lambda N} + \frac{(\frac{k}{w_1})(z_1d-d_0)}{\lambda N}\right).
\end{align*}
Since $\psi_{\Delta_1}(d_1)$ is smooth at scale $\Delta_1$, the contribution of $d$ to $\Sigma_5(z)$ is zero unless $|z_1d-d_0|\ll \Delta_1$. Then
\begin{align*}
\left| \frac{(\frac{k}{w_1})(z_1d-d_0)}{\lambda N} \right| \ll \dfrac{\left(\frac{w_1 \Lambda N}{x^{5\varepsilon}\Delta_1}\right)\Delta_1}{w_1 \Lambda N} \ll \dfrac{1}{x^{5\varepsilon}}.
\end{align*}
Choosing $J= \lceil \frac{20}{\varepsilon} \rceil$ and using Taylor expansions, we may write 
\begin{align} \label{equ:taylorexpansion}
\psi_N\!\left(\frac{\widetilde{\lambda} n +k(\frac{z_1}{w_1}) d}{\lambda}\right)= \sum_{j=0}^J  \dfrac{1}{j!} \left( \frac{(\frac{k}{w_1})(z_1d-d_0)}{\lambda N}  \right)^j \psi^{(j)}\left(\frac{\widetilde{\lambda} n +(\frac{k}{w_1}) d_0}{\lambda N}  \right) + O_\varepsilon(x^{-100}).
\end{align}
Consider the function $\psi_\star(t) = \psi(t)\psi^{(j)}((\frac{\widetilde{\lambda}}{\lambda})t+(\frac{k d_0}{w_1 \lambda N} ))$.  Since $\lambda, \widetilde{\lambda} \sim \Lambda$, we have $|\frac{\widetilde{\lambda}}{\lambda}| \ll 1$.  So the Leibniz product rule and $|\psi^{(i)}(t)|\ll_i \log(x)^{O_i(1)}$ tell us that $\psi_\star(t)$ is also smooth with $|\psi_\star^{(i)}(t)|\ll_i \log(x)^{O_i(1)}$. Further, the support of $\psi_\star(t)$ is contained in the support of $\psi(t)$. Overall, we see that $\psi_\star(n/N)$ is a coefficient sequence which is smooth at scale $N$.

\vspace{3mm}
Since $\psi_{\Delta_1}(d_1)$ is a smooth coefficient sequence, there also exists a smooth function $\phi:\mathbb{R} \rightarrow \mathbb{C}$, supported on $[c,C]$, with $|\phi^{(j)}(t)|\ll_j \log(x)^{O_j(1)}$ and $\psi_{\Delta_1}(d_1)=\phi(d_1/\Delta_1)$.  Now set $\phi_\star(t) = \frac{1}{j!}((\frac{k}{w_1})(\frac{\Delta_1}{\lambda N}) t)^j \phi(t)$. We already know that $|(\frac{k}{w_1})(\frac{\Delta_1}{\lambda N}) | \ll \frac{1}{x^{5\varepsilon}} \ll 1$ and so by the product rule we find  that $|\phi_\star^{(i)}(t)|\ll_i \log(x)^{O_i(1)}$. Hence the coefficient sequence $\phi_\star(d_1/\Delta_1)$ is smooth at scale $\Delta_1$.

\vspace{3mm}
 We substitute the RHS of (\ref{equ:taylorexpansion}) for $\psi_N( (\widetilde{\lambda}n+\frac{z_1}{w_1}kd)/\lambda)$ in $\Sigma_5(z)$. The contribution  of $1_{[cN,CN]}(n)O(x^{-100})$  is trivially bounded by  $O_\varepsilon(x^{-90})$, while the substitution of $\frac{1}{j!}((\frac{k}{w_1 \lambda N}) (z_1 d -d_0))^j\psi^{(j)}((\widetilde{\lambda}n+k d_0/w_1)/(\lambda N))$ creates a new smooth coefficient sequence of the form $\psi_\star(n/N)$ and a new shifted smooth coefficient sequence of the form $\phi_\star((z_1 d- d_0)/\Delta_1)$. Note that $J = \lceil \frac{20}{\varepsilon} \rceil = O_\varepsilon(1)$.  We define 
 \begin{align*}
\Upsilon_1(\psi^\star_N,\psi^\star_{\Delta_1})=\Bigg|\!\!\!\!\!\!\!\sum_{\substack{k \\ |k| \ll \frac{w_1\Lambda N}{x^{5\varepsilon}\Delta_1}\\ ((\frac{z_1}{w_1}) k   +  (\lambda-\widetilde{\lambda})B,m) \leq T(z)}}\!\!\!\!\sum_{\substack{d \\  (d, m \lambda \widetilde{\lambda})=1 }} \psi^*_{\Delta_1}(z_1 d-d_0)  \sum_{\substack{n }}^{\,\,*_6} C(n;d)  C\!\left(\dfrac{\widetilde{\lambda} n +k(\frac{z_1}{w_1}) d}{\lambda};d\right)\psi_N^*(n) \, (\star)\Bigg|.
\end{align*}
Taking the supremum over all coefficient sequences   $\psi_N^*(n)$ and $\psi_{\Delta_1}^*(d_1)$,   smooth at scales $N$ and $\Delta_1$, then
 $$\Sigma_5(z) \ll_\varepsilon     \sup_{(\psi_N^*,\psi_{\Delta_1}^*)} \Upsilon_1(\psi^\star_N,\psi^\star_{\Delta_1}) + O_\varepsilon(x^{-90}).$$
Hence we have the desired bound on $\Sigma_5(z)$ if we can show that $\Upsilon_1 \ll_\varepsilon \min\{\tfrac{(w_2,m)}{x^{\delta+10\varepsilon}H^3}, \tfrac{ (w_2,m)}{H^4}\}\tfrac{q_0^2(q_0,\ell)^2 N^2 }{x^{27\varepsilon}} $. 
 
\vspace{3mm }
 
Next we partition according to the value of $k$ and write  \begin{align*}
\Upsilon_2^{(k)}&= \Bigg| \!\!\!\!\!\sum_{\substack{d \\  (d, m \lambda \widetilde{\lambda})=1 }} \!\!\!\!\!\!\!\!\psi^*_{\Delta_1}(z_1 d-d_0) \sum_{\substack{n }}^{\,\,*_6} \! C(n;d)  C\!\left(\dfrac{\widetilde{\lambda} n +k(\frac{z_1}{w_1}) d}{\lambda};d\right)\psi_N^*(n)e_{m}\!\!\left(\dfrac{A ((\frac{z_1}{w_1}) k   +  (\lambda-\widetilde{\lambda})B)}{(n\!+\!B d)(\frac{1}{\lambda}(\widetilde{\lambda}n+(\frac{z_1}{w_1})kd)\!+\!B d)}\right)\!\Bigg|. 
\end{align*}
Suppose that $(\xi_k)$ is a sequence with  $\sum_{|k| \ll \frac{w_1\Lambda N}{x^{5\varepsilon}\Delta_1}} \xi_k 1_{((\frac{z_1}{w_1}) k   +  (\lambda-\widetilde{\lambda})B,m) \leq T(z)} = O_\varepsilon(1)$. If $\Upsilon_2^{(k)} =O_\varepsilon( \xi_k I)$ for each $k \in \mathbb{Z}$ with $|k| \ll \frac{w_1\Lambda N}{x^{5\varepsilon}\Delta_1}$ and $((\frac{z_1}{w_1}) k   +  (\lambda-\widetilde{\lambda})B,m) \leq T(z)$, then also $\Upsilon_1 = O_\varepsilon(I)$.  

\vspace{3mm} Finally, we also remove $C(n;d)$. Recall that $C(n;d)$ is the indicator function of a set $E_d$ mod $q_0$, which depends on the value of $d$ mod $q_0$. In particular, if we fix the values of $n$ and $d$ mod $q_0$, then $C(n;d)$ is constant. Since $\lambda, \widetilde{\lambda}$ and $k$ are all already fixed in $\Upsilon_2^{(k)}$, the function $C((\widetilde{\lambda}n+k(\frac{w_1}{z_1})d)/\lambda;d)$ is also constant when $n$ and $d$ are fixed mod $q_0$. But $E_d$ has at most $(q_0,\ell)$ elements for each $d$ mod $q_0$ and thus $C(n;d)$ is only non-zero for at most $q_0(q_0,\ell)$ choices of $n$ and $d$ mod $q_0$. We recall
\begin{align*}
\Sigma_6(z,(\psi_N^*,\psi_{\Delta_1}^*,d_\star,n_\star),k)=  \Bigg| \sum_{\substack{d \\ d \equiv d_\star (q_0) \\  (d, m \lambda \widetilde{\lambda})=1 }} \!\!\!\!\!\!\!\psi^*_{\Delta_1}(z_1 d-d_0) \!\!\!\!\sum_{\substack{n \\  n \equiv n_\star (q_0)    }}^{\,\,*_6} \!\!\!\!\!\psi_N^*(n)e_{m}\!\!\left(\dfrac{A ((\frac{z_1}{w_1}) k   +  (\lambda-\widetilde{\lambda})B)}{(n+B d)(\frac{1}{\lambda}(\widetilde{\lambda}n+(\frac{z_1}{w_1})kd)+B d)}\right)\!\Bigg|
\end{align*}
and observe that $$\sup_{(\psi_N^*,\psi_{\Delta_1}^*)} \Upsilon_2^{(k)} \ll q_0 (q_0,\ell) \sup_{(\psi_N^*,\psi_{\Delta_1}^*,d_\star,n_\star)} \Sigma_6(z,(\psi_N^*,\psi_{\Delta_1}^*,d_\star,n_\star),k).$$ The sum $\sum\limits^{\,\,*_6}$ is restricted to $n$ which satisfy $\lambda \mid \widetilde{\lambda} n + k (\frac{z_1}{w_1}) d$. But $w_2 \mid \lambda$ and $w_2 \mid \widetilde{\lambda}$, while $(z_1, m) = (d,m)=1$, so that $((\frac{z_1}{w_1}) d, w_2)=1$. Hence if $w_2 \nmid k$, we immediately get $\Sigma_6(z,w,k)= 0$. Additionally, we only need to consider $k$ which satisfy $((\frac{z_1}{w_1})k+(\lambda-\widetilde{\lambda})B,m) \leq T(z)$. 
But by the assumptions of Lemma~\ref{lem:adjustments1}, there exists $(\xi_k) \in \Xi(z)$ so that for all $|k| \ll \frac{w_1\Lambda N}{x^{5\varepsilon}\Delta_1}$ 
with  $w_2 \mid k$ and $((\frac{z_1}{w_1})k+(\lambda-\widetilde{\lambda})B,m) \leq T(z)$, 
\begin{align*}
\sup_{w \in W_6(z)} \Sigma_6(z,w,k)  \ll_\varepsilon \xi_k  \left(\min\left\{\frac{(w_2,m)}{x^{\delta+10\varepsilon}H^3},\frac{(w_2,m)}{H^4}\right\}\frac{q_0(q_0,\ell)  N^2}{ x^{27\varepsilon}} \right).
\end{align*}
 Then also  $\Upsilon_2^{(k)} \ll_\varepsilon \xi_k \min\{\tfrac{(w_2,m)}{x^{\delta+10\varepsilon}H^3}, \tfrac{ (w_2,m)}{H^4}\}\tfrac{q_0^2(q_0,\ell)^2 N^2 }{x^{27\varepsilon}} $  and $\Upsilon_1 \ll_\varepsilon \min\{\tfrac{(w_2,m)}{x^{\delta+10\varepsilon}H^3}, \tfrac{ (w_2,m)}{H^4}\}\tfrac{q_0^2(q_0,\ell)^2 N^2 }{x^{27\varepsilon}} $. It follows that  $\Sigma_5(z) \ll_\varepsilon \min\{\tfrac{(w_2,m)}{x^{\delta+10\varepsilon}H^3}, \tfrac{ (w_2,m)}{H^4}\}\tfrac{q_0^2(q_0,\ell)^2 N^2 }{x^{27\varepsilon}} $,  as proposed.
\end{proof}

\subsubsection{Divisibility by $\lambda$} 
Next we remove the condition $\lambda \mid \widetilde{\lambda} n +(\frac{z_1}{w_1})kd$.

\begin{lemma}[Removing a divisibility condition]\label{lem:adjustments2a} Let $Z_5$  be as given in Lemma~{\rm\ref{lem:summary}}  and let $W_6$ and  $\Sigma_6$ be as given in   Lemma~{\rm\ref{lem:adjustments1}}. For $z \in Z_5(x)$ and $k \in \mathbb{Z}$, denote by $W_7(z,k;x)$ the set of tuples
$
(\psi_{\frac{\Delta^*}{z_1}}, \psi_{(\frac{\lambda^*}{\lambda})N},  F_k, G_k)
$
which satisfy the following conditions: 
\begin{enumerate}[{\rm (i)}]
\item $\psi_{\frac{\Delta^*}{z_1}}(d)$  and $ \psi_{(\frac{\lambda^*}{\lambda})N}(n)$ are shifted smooth coefficient sequences at scales $(\frac{\Delta^*}{z_1}) $ and   $(\frac{\lambda^*}{\lambda})N$, with corresponding shifts $x_0 \ll x$.  Here  $\Delta^* = \min\{\frac{N}{\Lambda x^{5\varepsilon}}, \Delta_1\}$  and $\lambda^* = (\lambda, \widetilde{\lambda})$. 
\item $F_k$ and $G_k$ are residue classes mod $m$ with: 
\end{enumerate}
\begin{align} \label{equ:relationshipFG}
 \left( \frac{\lambda^*}{\widetilde{\lambda}}\right)(G_k+B) -\left(\frac{\lambda^*}{\lambda}\right) (F_k +B ) \equiv \left\{\left(\frac{\lambda^*}{\lambda }\right)\left(\frac{w_2}{ \widetilde{\lambda}}\right)\right\}\left\{\dfrac{1}{w_2}\left(\left(\frac{z_1}{w_1}\right)k + (\lambda-\widetilde{\lambda})B\right)\right\} \mbox{\rm  mod } m.
\end{align}

\vspace{3mm}
For given $x>1$, $k \in \mathbb{Z}$,  $z \in Z_5(x)$, $w_6 \in W_6(z;x)$ and $w_7 \in W_7(z,k;x)$,  
  we set 
 \begin{align*}
&\Sigma_7(z,w_6,w_7,k;x) =  \Bigg|\!\! \sum_{\substack{d \\ d \equiv d_\star (q_0) \\  (d, m \lambda \widetilde{\lambda})=1 }} \!\!\!\!\psi_{\frac{\Delta^*}{z_1}}(d) \sum_{\substack{n_1  }}^{\,\,*_7} \psi_{(\frac{\lambda^*}{\lambda})N}(n_1)
e_{m}\!\!\left(\dfrac{(\frac{\lambda^*}{\lambda})(\frac{\lambda^*}{\widetilde{\lambda}})A ((\frac{z_1}{w_1}) k   +  (\lambda-\widetilde{\lambda})B)}{(n_1 + (\frac{\lambda^*}{\lambda}) (F_k +B )d)( n_1 + ( \frac{\lambda^*}{\widetilde{\lambda}})(G_k+B) d)}\right) \Bigg|,
\end{align*} 

\vspace{-2mm}
 where $\sum\limits^{\,\,*_7}$ denotes summation over $n_1$ which satisfy  $(\frac{\lambda}{\lambda^*})n_1 + F_k d \equiv n_\star (q_0)$, $ ((\frac{\lambda}{\lambda^*})n_1 + F_k d,w_1 c_1)=1$,  $(( \frac{\widetilde{\lambda}}{\lambda^*}) n_1 + G_kd,w_1 c_1)=1$,  $((\frac{\lambda}{\lambda^*})n_1 + F_k d+ld, c_2)=1$ and $(( \frac{\widetilde{\lambda}}{\lambda^*}) n_1 + G_kd+l d, c_2)=1$.

\vspace{3mm}
For  $x>1$,  $z \in Z_5(x)$, $w_6 \in W_6(z;x)$ and $k $  
 with $|k| \ll  \frac{w_1\Lambda N}{x^{5\varepsilon}\Delta_1}$, $w_2 \mid k$ and  $((\frac{z_1}{w_1})k+(\lambda-\widetilde{\lambda})B,m) \leq T(z;x)$,
 \begin{align*}
 \Sigma_6(z,w_6,k;x) \ll_\varepsilon \sup_{w_7 \in W_7(z,k;x)} \left( \frac{\Delta_1}{\Delta^*}\right) (\Sigma_7(z,w_6,w_7,k;x) + O_\varepsilon(x^{-90})).
\end{align*}   
\end{lemma}
In particular, in order to show that $\Sigma_6(z,w_6,k) \ll_\varepsilon \xi_k  \left(\min\left\{\frac{(w_2,m)}{x^{\delta+10\varepsilon}H^3},\frac{(w_2,m)}{H^4}\right\}\frac{q_0(q_0,\ell)  N^2}{ x^{27\varepsilon}} \right)$, it suffices to prove  $\Sigma_7(z,w_6,w_7,k) + O_\varepsilon(x^{-90}) \ll_\varepsilon \xi_k  \left( \frac{\Delta^*}{\Delta_1}\right)\left(\min\left\{\frac{(w_2,m)}{x^{\delta+10\varepsilon}H^3},\frac{(w_2,m)}{H^4}\right\}\frac{q_0(q_0,\ell)  N^2}{ x^{27\varepsilon}} \right)$ for all $w_7 \in W_7(z,k)$.

\begin{proof} Our goal is to bound $\Sigma_6(z,w_6,k)$, which is given by
\begin{align*}
\Sigma_6(z,w_6,k)&=  \Bigg| \sum_{\substack{d \\ d \equiv d_\star (q_0) \\  (d, m \lambda \widetilde{\lambda})=1 }} \psi_{\Delta_1}^*(z_1 d-d_0) \!\!\!\!\!\!\!\!\!  \!\!\!\!\!\!\!\!\! \sum_{\substack{n \\  n \equiv n_\star (q_0) \\  \lambda \mid \widetilde{\lambda} n +k(\frac{z_1}{w_1})d \\
 (n,w_1 c_1)= (n+ld, c_2)=1 \\ 
 (\frac{1}{\lambda}(\widetilde{\lambda}n+k(\frac{z_1}{w_1})d),w_1 c_1)=1 \\ 
 \\ (\frac{1}{\lambda}(\widetilde{\lambda}n+k(\frac{z_1}{w_1})d)+ld, c_2)=1  }} \!\!\!\!\!\!\!\!\! \!\!\!\!\!\!\!\!\!  \psi_N^*(n)e_{m}\!\!\left(\dfrac{A ((\frac{z_1}{w_1}) k   +  (\lambda-\widetilde{\lambda})B)}{(n+B d)(\frac{1}{\lambda}(\widetilde{\lambda}n+(\frac{z_1}{w_1})kd)+B d)}\right)\Bigg|.
\end{align*}
Here we assume that $k \in \mathbb{Z}$  
 with $|k| \ll  \frac{w_1\Lambda N}{x^{5\varepsilon}\Delta_1}$, $w_2 \mid k$ and  $((\frac{z_1}{w_1})k+(\lambda-\widetilde{\lambda})B,m) \leq T(z)$.
 
 \vspace{3mm}
We now consider the condition $\lambda \mid \widetilde{\lambda} n +(\frac{z_1}{w_1})kd$. Write $\lambda^* = (\lambda, \widetilde{\lambda})$. The divisibility condition implies that $\lambda^* \mid (\frac{z_1}{w_1})kd$. But $(d,m \lambda \widetilde{\lambda})=1$ in $\Sigma_6$. So if $\lambda^* \nmid (\frac{z_1}{w_1})k$, then  $\Sigma_6(z,w_6,k) =0$. Hence we may assume that $\lambda^* \mid (\frac{z_1}{w_1})k$.  Note that $\lambda \mid \widetilde{\lambda} n +(\frac{z_1}{w_1})kd$ is equivalent to $\frac{\lambda}{\lambda^*} \mid \frac{\widetilde{\lambda}}{\lambda^*} n +(\frac{z_1k}{w_1 \lambda^*})d$. But since $(\lambda/\lambda^*, \widetilde{\lambda}/\lambda^*)=1$, this is further equivalent to $ n  \equiv -(\frac{\widetilde{\lambda}}{\lambda^*})^{-1}(\frac{z_1k}{w_1\lambda^* }) d  $ mod $\frac{\lambda}{\lambda^*}$.  Choose $F_k \in \{1, \dots, \frac{\lambda}{\lambda^*}\}$ with $F_k \equiv  -(\frac{\widetilde{\lambda}}{\lambda^*})^{-1}(\frac{z_1k}{w_1\lambda^* }) $ mod $\frac{\lambda}{\lambda^*}$.
 Then $n = F_k d + (\frac{\lambda}{\lambda^*})n_1$ for some $n_1 \in \mathbb{Z}$. Further, $\frac{1}{\lambda}(\widetilde{\lambda}n+(\frac{z_1}{w_1})kd) =\frac{1}{\lambda}(\widetilde{\lambda}(F_k d + (\frac{\lambda}{\lambda^*})n_1)+(\frac{z_1}{w_1})kd)$. But there exists $G_k \in \mathbb{Z}$ such that $(\frac{\widetilde{\lambda}}{\lambda^*}) F_k = - \frac{1}{\lambda^*}(\frac{z_1}{w_1})k  + G_k(\frac{\lambda}{\lambda^*})$ and so 
 $\frac{1}{\lambda} (\widetilde{\lambda}F_kd + (\frac{z_1}{w_1})kd)=  G_kd$. Hence  
 \begin{align*}
 n =\left(\frac{\lambda}{\lambda^*}\right)n_1 + F_k d   \qquad \mbox{ and } \qquad \frac{1}{\lambda}\left(\widetilde{\lambda}n+\left(\frac{z_1}{w_1}\right)kd\right) = \left( \dfrac{\widetilde{\lambda}}{\lambda^*}\right) n_1 + G_kd. 
 \end{align*}
 Therefore, $\Sigma_6(z,w_6,k)$ transforms as follows: 
 \begin{align*} 
\Sigma_6&=\Bigg|\!\sum_{\substack{d \\ d \equiv d_\star (q_0) \\  (d, m \lambda \widetilde{\lambda})=1 }} \!\!\!\!\!\!\!\psi_{\Delta_1}^*(z_1 d-d_0) \!\!\!\!\!\!\!\!\!\!\!\!\!\!\!\!\!\!\!\!\!\sum_{\substack{n_1 \\  (\frac{\lambda}{\lambda^*})n_1 + F_k d \equiv n_\star (q_0) \\ ((\frac{\lambda}{\lambda^*})n_1 + F_k d,w_1 c_1)=1 \\ (( \frac{\widetilde{\lambda}}{\lambda^*}) n_1 + G_kd,w_1 c_1)=1 \\ ((\frac{\lambda}{\lambda^*})n_1 + F_k d+l d, c_2)=1  \\ (( \frac{\widetilde{\lambda}}{\lambda^*}) n_1 + G_kd+l d, c_2)=1    }} \!\!\!\!\!\!\!\!\!\!\!\!\!\!\!\!\!\!\!\!\psi_N^*((\tfrac{\lambda}{\lambda^*})n_1 + F_k d)e_{m}\!\!\left(\!\dfrac{A ((\frac{z_1}{w_1}) k   +  (\lambda-\widetilde{\lambda})B)}{((\frac{\lambda}{\lambda^*})n_1 + F_k d+B d)(( \frac{\widetilde{\lambda}}{\lambda^*}) n_1 + G_kd+B d)}\!\right)\!\Bigg|.
\end{align*}
Recall that $(\lambda, m^*) =(\widetilde{\lambda},m^*)=w_2$. Here $w_2$ is the largest factor of $\lambda$ (and $\widetilde{\lambda}$) which only consists of powers of primes that divide $m$. Since $w_2 \mid \lambda^*$, this further implies that $(\frac{\lambda}{\lambda^*})$ and $(\frac{\widetilde{\lambda}}{\lambda^*})$ are coprime to $m$.
  \begin{align*}
\Sigma_6&=  \Bigg|\sum_{\substack{d \\ d \equiv d_\star (q_0) \\  (d, m \lambda \widetilde{\lambda})=1 }} \!\!\!\!\!\!\!\psi_{\Delta_1}^*(z_1 d-d_0)\sum_{\substack{n_1  }}^{\,\,*_7}\psi_N^*((\tfrac{\lambda}{\lambda^*})n_1 + F_k d)e_{m}\!\!\left(\dfrac{(\frac{\lambda^*}{\lambda})(\frac{\lambda^*}{\widetilde{\lambda}})A ((\frac{z_1}{w_1}) k   +  (\lambda-\widetilde{\lambda})B)}{(n_1 + (\frac{\lambda^*}{\lambda}) (F_k +B )d)( n_1 + ( \frac{\lambda^*}{\widetilde{\lambda}})(G_k+B) d)}\right)\Bigg|.
\end{align*}

 \vspace{3mm}
 Observe now  that  $( \frac{\lambda^*}{\widetilde{\lambda}})(G_k+B) -(\frac{\lambda^*}{\lambda}) (F_k +B ) = (\frac{\lambda^*}{\lambda \widetilde{\lambda}})(\frac{z_1}{w_1})k + (\frac{\lambda^*}{\lambda \widetilde{\lambda}})(\lambda-\widetilde{\lambda})B$. We assumed that  $w_2 \mid k$. Hence  we may write 
 \begin{align} \label{equ:congruencyrelations}
 \left( \frac{\lambda^*}{\widetilde{\lambda}}\right)(G_k+B) -\left(\frac{\lambda^*}{\lambda}\right) (F_k +B ) \equiv \left\{\left(\frac{\lambda^*}{\lambda }\right)\left(\frac{w_2}{ \widetilde{\lambda}}\right)\right\}\left\{\dfrac{1}{w_2}\left(\left(\frac{z_1}{w_1}\right)k + (\lambda-\widetilde{\lambda})B\right)\right\} \mbox{ mod } m,
 \end{align}
 where both of  the two factors surrounded by $\{ \}$--brackets are well-defined residue classes mod $m$.

  \vspace{3mm}   
  Next we use Taylor expansions to treat $\psi_N^*((\frac{\lambda}{\lambda^*})n_1 + F_k d)$, which depends on both $n_1$ and $d$. Currently, the summation range of $d$ is still a bit too long, so we first use a partition of unity to shorten it: Choose a smooth function $\phi: \mathbb{R} \rightarrow [0,\infty)$ with $\mbox{supp}(\phi) \subseteq [-2,2]$ and  $\phi(t) \geq 1$ on $[-1,1]$. Define $\phi_m: \mathbb{R} \rightarrow [0,\infty)$ by setting $\phi_m(t) = (\sum_{j=-4}^4 \phi(t-m+j))^{-1} \phi(t-m)$ on $[m-3,m+3]$ and $\phi_m(t)=0$ elsewhere. The function $\phi_m$ is well-defined and smooth since $ (\sum_{j=-4}^4 \phi(t-m+j))^{-1} > 0$ on $[m-3,m+3]$. Further, for every $t \in \mathbb{R}$, 
  $\sum_{m=-\infty}^\infty \phi_m(t)=1$. Now recall that  $\psi^*_{\Delta_1}(d_1)$ is a smooth coefficient sequence and there is a smooth function $\psi:\mathbb{R} \rightarrow \mathbb{C}$, supported on $[c,C]$, with $|\psi^{(j)}(t)|\ll_j \log(x)^{O_j(1)}$ and $\psi^*_{\Delta_1}(d_1)=\psi(\frac{d_1}{\Delta_1})$. Set 
  \begin{align*}
  \Delta^* = \min\left\{\dfrac{N}{\Lambda x^{5\varepsilon}}, \Delta_1\right\} 
  \end{align*}
  and observe that
  \begin{align*}
  \psi^*_{\Delta_1}(z_1 d -d_0) &= \psi\left(\left(\frac{\Delta^*}{\Delta_1}\right)\dfrac{(z_1 d -d_0)}{\Delta^*}\right) \sum_{m=-\infty}^\infty \phi_m\left(\dfrac{z_1 d -d_0}{\Delta^*}\right)\\
  &=  \sum_{m=-\infty}^\infty\psi\left(\left(\frac{\Delta^*}{\Delta_1}\right) \frac{\left(d -d_0/z_1\right)}{(\Delta^*/z_1)}\right) \phi_5\left(\dfrac{(d -d_0/z_1) }{(\Delta^*/z_1)}- (m-5)\right) \\
   &=  \sum_{m \ll \Delta_1/\Delta^*} \tau_m\left(\frac{d -x_m}{(\Delta^*/z_1)} \right),
  \end{align*}
 where $x_m = d_0/z_1 + (\Delta^*/z_1)(m-5) \in \mathbb{R}$ and $\tau_m(t) = \psi((\frac{\Delta^*}{\Delta_1})(t+m-5)) \phi_5(t)$. In the last line we used that $\tau_m$ is identical to zero unless  $m \ll \Delta_1/\Delta^*$. The support of $\tau_m$ is contained in $[3,7]$ and since $\Delta^*/\Delta_1 \ll 1$, we also have $|\tau^{(j)}(t)| \ll_j \log(x)^{O_j(1)}$.  Thus $\tau_m((d-x_m)/(\Delta^*/z_1))$ is a coefficient sequence which is shifted smooth at scale $\Delta^*/z_1$. We then define 
 \begin{align*}
\Upsilon_0(\psi^\circ_{\frac{\Delta^*}{z_1}})&=  \Bigg| \sum_{\substack{d \\ d \equiv d_\star (q_0) \\  (d, m \lambda \widetilde{\lambda})=1 }} \!\!\!\!\!\!\!\psi_{\frac{\Delta^*}{z_1}}^\circ(d)\sum_{\substack{n_1  }}^{\,\,*_7}\psi_N^*((\tfrac{\lambda}{\lambda^*})n_1 + F_k d)e_{m}\!\!\left(\dfrac{(\frac{\lambda^*}{\lambda})(\frac{\lambda^*}{\widetilde{\lambda}})A ((\frac{z_1}{w_1}) k   +  (\lambda-\widetilde{\lambda})B)}{(n_1 + (\frac{\lambda^*}{\lambda}) (F_k +B )d)( n_1 + ( \frac{\lambda^*}{\widetilde{\lambda}})(G_k+B) d)}\right)\Bigg|.
\end{align*}
We now substitute $\sum \tau_m(\frac{d-x_m}{(\Delta^*/z_1)})$ for $\psi^*_{\Delta_1}(z_1d-d_0)$ in $\Sigma_6$. Taking a supremum over all   coefficient sequences $\psi^\circ_{\frac{\Delta^*}{z_1}}(d)$ which are shifted smooth at scale $\frac{\Delta^*}{z_1}$ and have shift $x_0 \ll x^{5\varepsilon}\Delta_1$, we find that
 \begin{align*}
 \Sigma_6 \ll_\varepsilon \left(\dfrac{\Delta_1}{\Delta^*}\right) \sup_{\psi^\circ_{\frac{\Delta^*}{z_1}}} \Upsilon_0(\psi^\circ_{\frac{\Delta^*}{z_1}}).
 \end{align*}
We have thus shortened the summation range of $d$. In particular, for a given $\Upsilon_0$, there exists $x_0 \in \mathbb{R}$  with $x_0 \ll x^{5\varepsilon}\Delta_1$ such that the contribution of $d$ to $\Upsilon_0$ is zero unless $|d-x_0| \ll \Delta^*/z_1$.  Note that 
  \begin{align*}
 \dfrac{ |F_k (d-x_0)|}{N} \ll \left|\frac{\lambda}{\lambda^*}\right|\dfrac{\Delta^*}{z_1 N} \ll \dfrac{1}{x^{5\varepsilon}}
  \end{align*}
  whenever $d$ contributes to $\Upsilon_0$.   We set $\psi^*_N(n) = \phi(\frac{n}{N})$. 
Choosing $J= \lceil \frac{20}{\varepsilon} \rceil$ and using Taylor expansions,  
\begin{align} \label{equ:taylorexpansion222}
\psi_N^*\!\left(\left(\tfrac{\lambda}{\lambda^*}\right)n_1+F_k d\right)= \sum_{j=0}^J  \dfrac{1}{j!} \left( \frac{F_k (d-x_0)}{N}  \right)^j \phi^{(j)}\left(\frac{n_1 + (\frac{\lambda^*}{\lambda})F_k x_0}{(\frac{\lambda^*}{\lambda}) N}  \right) + O_\varepsilon(x^{-100} 1_{n_1 \ll x}).
\end{align}
The sequence $\phi^{(j)}(\frac{n_1 + (\frac{\lambda^*}{\lambda})F_k x_0}{(\frac{\lambda^*}{\lambda}) N})$ is a shifted smooth coefficient sequence in $n_1$ at scale $(\frac{\lambda^*}{\lambda})N$, while the sequence 
$\frac{1}{j!}\left(\frac{F_k (\Delta^*/z_1)}{N}\frac{(d-x_0)}{(\Delta^*/z_1)}\right)^j\psi^\circ_{\frac{\Delta^*}{z_1}}(d)$ is a shifted smooth coefficient sequence in $d$ at scale $\frac{\Delta^*}{z_1}$.   We now substitute the RHS of (\ref{equ:taylorexpansion222}) for   $\psi_N^*((\frac{\lambda}{\lambda^*})n_1 + F_k d)$ in $\Upsilon_0$. 
 The contribution of $O_\varepsilon(x^{-100}1_{n_1 \ll x})$  is trivially bounded by  $O_\varepsilon(x^{-90})$. For the treatment of other terms, set 
 \begin{align*}
\Upsilon_1(\psi_{(\frac{\lambda^*}{\lambda})N},\psi_{\frac{\Delta^*}{z_1}}) &= \Bigg|\!\!\! \sum_{\substack{d \\ d \equiv d_\star (q_0) \\  (d, m \lambda \widetilde{\lambda})=1 }} \!\!\!\!\!\psi_{\frac{\Delta^*}{z_1}}(d) \sum_{\substack{n_1  }}^{\,\,*_7} \psi_{(\frac{\lambda^*}{\lambda})N}(n_1)
e_{m}\!\!\left(\dfrac{(\frac{\lambda^*}{\lambda})(\frac{\lambda^*}{\widetilde{\lambda}})A ((\frac{z_1}{w_1}) k   +  (\lambda-\widetilde{\lambda})B)}{(n_1 + (\frac{\lambda^*}{\lambda}) (F_k +B )d)( n_1 + ( \frac{\lambda^*}{\widetilde{\lambda}})(G_k+B) d)}\right)\Bigg|. 
\end{align*}
Here we consider  coefficient sequences $\psi_{(\frac{\lambda^*}{\lambda})N}(n_1)$  and $\psi_{\frac{\Delta^*}{z_1}}(d)$ which are shifted smooth at scales $(\frac{\lambda^*}{\lambda})N$ and $\frac{\Delta^*}{z_1}$, and for which the corresponding shifts $x_0$ are bounded by $x_0 \ll x$. By the discussion above, we find
$ \Upsilon_0(\psi^\circ_{\frac{\Delta^*}{z_1}}) \ll_\varepsilon
\sup_{(\psi_{(\frac{\lambda^*}{\lambda})N},\psi_{\frac{\Delta^*}{z_1}})}\Upsilon_1(\psi_{(\frac{\lambda^*}{\lambda})N},\psi_{\frac{\Delta^*}{z_1}}) + O_\varepsilon(x^{-90})
$. Hence, 
$$ \Sigma_6(z,w_6,k)\ll_\varepsilon  \sup_{(\psi_{(\frac{\lambda^*}{\lambda})N},\psi_{\frac{\Delta^*}{z_1}})}\left(\frac{\Delta_1}{\Delta^*}\right)(\Upsilon_1(\psi_{(\frac{\lambda^*}{\lambda})N},\psi_{\frac{\Delta^*}{z_1}}) + O_\varepsilon(x^{-90})).  $$
 Since $\Upsilon_1(\psi_{(\frac{\lambda^*}{\lambda})N},\psi_{\frac{\Delta^*}{z_1}}) = \Sigma_7(z,w_6,w_7,k)$ for some $w_7 \in W_7(z,k)$, this concludes the proof.
\end{proof}

\subsubsection{Coprimality conditions}  For a final simplification, we remove the various coprimality conditions by using Möbius inversion.

\begin{lemma}[Removal of coprimality conditions]\label{lem:adjustments2}  Let $Z_5$ and $\Sigma_7$ be as given in Lemma~{\rm\ref{lem:summary}} and  Lemma~{\rm\ref{lem:adjustments2a}}.
 For $z \in Z_5(x)$, denote by $W_8(z;x)$ the set of tuples
\begin{align*}
(\Delta_2, N_2, \psi_{\Delta_2/z_1}, \psi_{N_2}, q_3, d^\star, n^\star, A_1, A_2, B_1)
\end{align*}
which satisfy the following conditions: 
\begin{enumerate}[{\rm (i)}]
\item $\Delta_2, N_2 \in (0,\infty)$ with $\Delta_2 \leq \Delta^*= \min\{\frac{N}{\Lambda x^{5\varepsilon}}, \Delta_1\}$ and $N_2 \leq N$. 
\item $\psi_{\Delta_2/z_1}(d)$  and $\psi_{N_2}(n)$  are shifted smooth coefficient sequences at scales $\Delta_2/z_1 $ and $N_2 $. 
\item $q_3 \in \mathbb{N}$ with 
 $q_3 \mid (\frac{m}{q_0})$.
\item  $d^\star$ and $n^\star$ are residue classes mod $(q_0q_3)$.
\item $A_1$, $A_2$, $B_1$ are residue classes mod $m$ with $(A_1,m)=(A_2,m)=1$. 
\end{enumerate}
For given $x>1$, $z \in Z_5(x)$, $w_8 \in W_8(z;x)$ and $k \in \mathbb{Z}$,  set 
 \begin{align*}
&\Sigma_8(z,w_8,k;x) =  \Bigg|  \sum_{\substack{d \\  d \equiv d^\star  (q_0q_3) }} \sum_{\substack{n \\   n \equiv n^\star  (q_0q_3)    }}  \psi_{\Delta_2/z_1}(d)\psi_{N_2}(n)  e_{m}\left(\dfrac{w_2A_1L_1}{( n+B_1 d)(  n+(B_1+L_1)d)}\right) \Bigg| \\
&\mbox{where }  L_1 \equiv  A_2 \left\{\dfrac{1}{w_2}\left(\left(\frac{z_1}{w_1}\right)k + (\lambda-\widetilde{\lambda})B\right)\right\} \mbox{\rm mod } m.
\end{align*}

Then for all $k \in \mathbb{Z}$ 
 with $|k| \ll  \frac{w_1\Lambda N}{x^{5\varepsilon}\Delta_1}$, $w_2 \mid k$ and  $((\frac{z_1}{w_1})k+(\lambda-\widetilde{\lambda})B,m) \leq T(z;x)$, and for all $x>1$, $z \in Z_5(x)$, $w_6 \in W_6(z;x)$ and $w_7 \in W_7(z,k;x)$, we have
 \begin{align*}
 \Sigma_7(z,w_6,w_7, k;x) \ll_\varepsilon \sup_{w_8 \in W_8(z;x)} \left( x^{2\varepsilon} q_3 \right) \Sigma_8(z,w_8,k;x) + O_\varepsilon(x^{-89}).
 \end{align*}
\end{lemma}
In particular, in order to show that $\Sigma_7(z,w_6,w_7,k) \ll_\varepsilon \xi_k  \left( \frac{\Delta^*}{\Delta_1}\right)\left(\min\left\{\frac{(w_2,m)}{x^{\delta+10\varepsilon}H^3},\frac{(w_2,m)}{H^4}\right\}\frac{q_0(q_0,\ell)  N^2}{ x^{27\varepsilon}} \right)$, it suffices to prove  that $\Sigma_8(z,w_8,k) + O_\varepsilon(x^{-89}) \ll_\varepsilon \xi_k   \left( \frac{1}{x^{2\varepsilon} q_3} \right) \left(\min\left\{\frac{(w_2,m)}{x^{\delta+10\varepsilon}H^3},\frac{(w_2,m)}{H^4}\right\}\frac{q_0(q_0,\ell)  N^2}{ x^{27\varepsilon}} \right)$.

\begin{proof}
 Our goal is to bound $\Sigma_7(z,w_6,w_7,k)$, which is given by
\begin{align*}
\Sigma_7=\Bigg|\sum_{\substack{d \\ d \equiv d_\star (q_0) \\  (d, m \lambda \widetilde{\lambda})=1 }} \psi_{\frac{\Delta^*}{z_1}}(d) \!\!\!\!\!\!\!\!\!\!\sum_{\substack{n_1  \\  (\frac{\lambda}{\lambda^*})n_1 + F_k d \equiv n_\star (q_0) \\ ((\frac{\lambda}{\lambda^*})n_1 + F_k d,w_1 c_1)=1 \\ (( \frac{\widetilde{\lambda}}{\lambda^*}) n_1 + G_kd,w_1 c_1)=1 \\ ((\frac{\lambda}{\lambda^*})n_1 + F_k d+ld, c_2)=1  \\ (( \frac{\widetilde{\lambda}}{\lambda^*}) n_1 + G_kd+ld, c_2)=1   }} \!\!\!\!\!\!\!\!\!\! \psi_{(\frac{\lambda^*}{\lambda})N}(n_1)
e_{m}\!\!\left(\dfrac{(\frac{\lambda^*}{\lambda})(\frac{\lambda^*}{\widetilde{\lambda}})A ((\frac{z_1}{w_1}) k   +  (\lambda-\widetilde{\lambda})B)}{(n_1 + (\frac{\lambda^*}{\lambda}) (F_k +B )d)( n_1 + ( \frac{\lambda^*}{\widetilde{\lambda}})(G_k+B) d)}\right) \Bigg|.
\end{align*}
Here we assume that $k \in \mathbb{Z}$  
 with $|k| \ll  \frac{w_1\Lambda N}{x^{5\varepsilon}\Delta_1}$, $w_2 \mid k$ and  $((\frac{z_1}{w_1})k+(\lambda-\widetilde{\lambda})B,m) \leq T(z)$.

\vspace{3mm} \textbf{Step 1: Removal of coprimality conditions (Part 1).} We first remove coprimality conditions which do not involve factors of $m$. Recall that $(\lambda,m^*)=w_2$ and $(\widetilde{\lambda},m^*)=w_2$. Hence we may replace $(d, m \lambda \widetilde{\lambda})=1$ by $(d,m (\lambda/w_2)(\widetilde{\lambda}/w_2))=1$, where both $\lambda/w_2$ and $\widetilde{\lambda}/w_2$ are coprime to $m$. Recall also that $(w_1,m)=1$. For the moment, we wish to remove the requirements $(d, (\lambda/w_2)(\widetilde{\lambda}/w_2))=1$, $((\frac{\lambda}{\lambda^*})n_1 + F_k d,w_1 )=1$ and $((\frac{\widetilde{\lambda}}{\lambda^*})n_1 + G_k d,w_1 )=1$. Using Möbius inversion, we observe that $\Sigma_7$ equals 
 \begin{align*}
 \Bigg|\!\!\sum_{\substack{ f_1, f_2, f_3 \\ f_1 \mid \frac{\lambda}{w_2}\frac{\widetilde{\lambda}}{w_2} \\ f_2 \mid w_1 \\ f_3 \mid w_1}}\!\!\!\!\!\mu(f_1)\mu(f_2)\mu(f_3)\!\!\!\!\!\sum_{\substack{d \\ d \equiv d_\star (q_0) \\ f_1 \mid d \\ (d, m )=1 }}\!\!\!\!\psi_{\frac{\Delta^*}{z_1}}(d) \!\!\!\!\!\!\!\!\!\!\!\!\!\!\!\sum_{\substack{n_1 \\  n_1 + (\frac{\lambda^*}{\lambda})F_k d \equiv (\frac{\lambda^*}{\lambda})n_\star (q_0) \\ f_2 \mid (\frac{\lambda}{\lambda^*})n_1 + F_k d \\ f_3 \mid ( \frac{\widetilde{\lambda}}{\lambda^*}) n_1 + G_kd \\ ((\frac{\lambda}{\lambda^*})n_1 + F_k d, c_1)=1 \\ (( \frac{\widetilde{\lambda}}{\lambda^*}) n_1 + G_kd, c_1)=1 \\ ((\frac{\lambda}{\lambda^*})n_1 + F_k d+ld, c_2)=1  \\ (( \frac{\widetilde{\lambda}}{\lambda^*}) n_1 + G_kd+ld, c_2)=1    }} \!\!\!\!\!\!\!\!\!\!\!\!\!\!\!\!\!\!\!\! \psi_{(\frac{\lambda^*}{\lambda})N}(n_1)
e_{m}\!\!\left(\dfrac{(\frac{\lambda^*}{\lambda})(\frac{\lambda^*}{\widetilde{\lambda}})A ((\frac{z_1}{w_1}) k   +  (\lambda-\widetilde{\lambda})B)}{(n_1 \!+\! (\frac{\lambda^*}{\lambda}) (F_k +B )d)( n_1 \!+\! ( \frac{\lambda^*}{\widetilde{\lambda}})(G_k+B) d)}\right)\!\Bigg|.
\end{align*} 
For given $f_1,f_2,f_3 \in\mathbb{N}$ with  $f_1 \mid \frac{\lambda}{w_2}\frac{\widetilde{\lambda}}{w_2}$, $f_2 \mid w_1$ and  $f_3 \mid w_1$, we define $\Upsilon_1(f_1,f_2,f_3)$ to equal
 \begin{align*}
\Upsilon_1 =\Bigg|\! \sum_{\substack{d_1 \\ f_1d_1 \equiv d_\star (q_0) \\ (d_1, m )=1 }}\!\!\!\!\!\!\psi_{\frac{\Delta^*}{z_1}}(f_1 d_1)\!\! \!\!\!\!\!\!\!\!\!\!\sum_{\substack{n_1 \\  n_1 + (\frac{\lambda^*}{\lambda})F_k f_1d_1 \equiv (\frac{\lambda^*}{\lambda})n_\star (q_0) \\ f_2 \mid (\frac{\lambda}{\lambda^*})n_1 + F_k f_1d_1 \\ f_3 \mid ( \frac{\widetilde{\lambda}}{\lambda^*}) n_1 + G_kf_1d_1 \\ ((\frac{\lambda}{\lambda^*})n_1 + F_k f_1d_1, c_1)=1 \\ (( \frac{\widetilde{\lambda}}{\lambda^*}) n_1 + G_kf_1d_1, c_1)=1 \\ ((\frac{\lambda}{\lambda^*})n_1 + F_k f_1d_1+l f_1d_1,  c_2)=1  \\ (( \frac{\widetilde{\lambda}}{\lambda^*}) n_1 + G_kf_1d_1+l f_1d_1,  c_2)=1    }} \!\!\!\!\!\!\!\!\!\!\!\!\!\!\!\!\!\!\!\!\!\!\!\!\!\!\!\! \psi_{(\frac{\lambda^*}{\lambda})N}(n_1)
e_{m}\!\!\left(\dfrac{(\frac{\lambda^*}{\lambda})(\frac{\lambda^*}{\widetilde{\lambda}})A ((\frac{z_1}{w_1}) k   +  (\lambda-\widetilde{\lambda})B)}{(n_1 \!+\! (\frac{\lambda^*}{\lambda}) (F_k \!+\!B )f_1d_1)( n_1 \!+\! ( \frac{\lambda^*}{\widetilde{\lambda}})(G_k\!+\!B) f_1d_1)}\right)  \Bigg|.
\end{align*}
Since $ \frac{\lambda}{w_2}\frac{\widetilde{\lambda}}{w_2}$ and $w_1$ have only $O_\varepsilon(x^{\varepsilon/3})$ divisors, we then have
\begin{align}\label{chain1}
\Sigma_7(z,w_6,w_7,k) \ll_\varepsilon \sup_{f_1,f_2,f_3} (x^\varepsilon) \Upsilon_1(f_1,f_2,f_3),
\end{align}
where the supremum is taken over $f_1,f_2,f_3 \in\mathbb{N}$ with  $f_1 \mid \frac{\lambda}{w_2}\frac{\widetilde{\lambda}}{w_2}$, $f_2 \mid w_1$ and  $f_3 \mid w_1$.

\vspace{3mm} Next we consider the conditions $  (\frac{\lambda}{\lambda^*})n_1  \equiv - F_k f_1d_1 $ mod $(f_2)$ and $  ( \frac{\widetilde{\lambda}}{\lambda^*}) n_1  \equiv -G_kf_1d_1  $ mod $(f_3)$. Write $f_2^* = (f_2, \lambda/\lambda^*)$ and $f_3^* = (f_3, \widetilde{\lambda}/\lambda^*)$.  Write $f_2^\circ=f_2^*/(f_2^*, F_k f_1)$ and $f_3^\circ=f_3^*/(f_3^*, G_k f_1)$ and set $f_{23} = [f_2^\circ, f_3^\circ]$. For $  (\frac{\lambda}{\lambda^*})n_1  \equiv - F_k f_1d_1 $ mod $(f_2)$ and $  ( \frac{\widetilde{\lambda}}{\lambda^*}) n_1  \equiv -G_kf_1d_1  $  mod $(f_3)$ to simultaneously have solutions, we require that $f_{23} \mid d_1$. Write $d_1 = f_{23} d_2$. In that case, the conditions are satisfied if $  n_1  \equiv - (\frac{\lambda}{\lambda^* f_2^*})^{-1}(\frac{F_k f_1f_{23}}{f_2^*})d_2 $ mod $(f_2/f_2^*)$ and $  n_1  \equiv -( \frac{\widetilde{\lambda}}{\lambda^*f_3^*})^{-1}(\frac{G_kf_1f_{23}}{f_3^*})d_2  $ mod $(f_3/f_3^*)$. Now recall that $w_1$ is squarefree. Hence $f_2/f_2^*$ and $f_3/f_3^*$ are also squarefree. Suppose $p \mid (f_2/f_2^*,f_3/f_3^*)$. If $- (\frac{\lambda}{\lambda^* f_2^*})^{-1}(\frac{F_k f_1f_{23}}{f_2^*}) \not\equiv -( \frac{\widetilde{\lambda}}{\lambda^*f_3^*})^{-1}(\frac{G_kf_1f_{23}}{f_3^*})$ mod $p$, then $- (\frac{\lambda}{\lambda^* f_2^*})^{-1}(\frac{F_k f_1f_{23}}{f_2^*})d_2 \equiv -( \frac{\widetilde{\lambda}}{\lambda^*f_3^*})^{-1}(\frac{G_kf_1f_{23}}{f_3^*}) d_2$ mod $p$ only holds if $p \mid d_2$. Denote the product of $p\mid (f_2/f_2^*,f_3/f_3^*)$ with  $- (\frac{\lambda}{\lambda^* f_2^*})^{-1}(\frac{F_k f_1f_{23}}{f_2^*}) \not\equiv -( \frac{\widetilde{\lambda}}{\lambda^*f_3^*})^{-1}(\frac{G_kf_1f_{23}}{f_3^*})$ mod $p$ by $f_4$. Choose $H_k \in \{1, \dots, \frac{1}{f_4} [f_2/f_2^*, f_3/f_3^*]\}$ with $H_k \equiv - (\frac{\lambda}{\lambda^* f_2^*})^{-1}(\frac{F_k f_1f_{23}}{f_2^*})$ mod $p$ for every $p \mid f_2/(f_2^*f_4)$ and $H_k \equiv  -( \frac{\widetilde{\lambda}}{\lambda^*f_3^*})^{-1}(\frac{G_kf_1f_{23}}{f_3^*})$ mod $p$ for every $p \mid f_3/(f_3^*f_4)$. Overall, we then find that the conditions $  (\frac{\lambda}{\lambda^*})n_1  \equiv - F_k f_1d_1 $ mod $(f_2)$ and $  ( \frac{\widetilde{\lambda}}{\lambda^*}) n_1  \equiv -G_kf_1d_1  $ mod $(f_3)$ are satisfied  if and only if $d_1 = f_{23} f_4 d_3$ for some $d_3 \in \mathbb{Z}$ and $n_1 = [f_2/f_2^*, f_3/f_3^*] n_3+H_k f_4 d_3$ for some $n_3 \in \mathbb{Z}$.  To simplify notation, we now set $f^*=f_1f_{23}f_4  $, where $f^* \mid \frac{\lambda}{w_2}\frac{\widetilde{\lambda}}{w_2}w_1$ and $(f^*,m)=1$, and we set $w_3 = [f_2/f_2^*, f_3/f_3^*]$, where $w_3 \mid w_1$ and $(w_3,m)=1$. Then $\Upsilon_1(f_1,f_2,f_3)$ transforms as follows: 
 \begin{align*}
&\Upsilon_1(f_1,f_2,f_3) =\Bigg| \sum_{\substack{d_3 \\  d_3 \equiv (f^*)^{-1}d_\star  (q_0) \\ (d_3, m )=1 }}\!\!\!\!\!\!\!\psi_{\frac{\Delta^*}{z_1}}(f^*d_3) \!\!\!\!\!\!\!\!\!\!\sum_{\substack{n_3 \\   n_3\equiv - \frac{1}{w_3}(H_kf_4 + (\frac{\lambda^*}{\lambda})F_k f^* )(f^*)^{-1}d_\star + \frac{1}{w_3} (\frac{\lambda^*}{\lambda})n_\star (q_0) \\ ((\frac{\lambda}{\lambda^*})w_3n_3+ ( (\frac{\lambda}{\lambda^*})H_kf_4 + F_k f^* )d_3, c_1)=1 \\ (( \frac{\widetilde{\lambda}}{\lambda^*})  w_3n_3+ (( \frac{\widetilde{\lambda}}{\lambda^*})H_kf_4 + G_kf^* )d_3, c_1)=1 \\ ((\frac{\lambda}{\lambda^*}) w_3n_3+ ((\frac{\lambda}{\lambda^*})H_kf_4+ F_k f^* +l f^* )d_3, c_2)=1  \\ (( \frac{\widetilde{\lambda}}{\lambda^*})  w_3n_3+ ( ( \frac{\widetilde{\lambda}}{\lambda^*})H_kf_4 + G_kf^* +l f^* )d_3, c_2)=1      }} \!\!\!\!\!\!\!\!\!\!\!\!\!\!\!\!\!\!\!\!\!\!\!\!\!\!\!\! \psi_{(\frac{\lambda^*}{\lambda})N}( w_3n_3+H_kf_4d_3)\, (\star)\Bigg|
 \\
&\mbox{where } (\star) = e_{m}\left(\dfrac{(\frac{1}{w_3})^2(\frac{\lambda^*}{\lambda})(\frac{\lambda^*}{\widetilde{\lambda}})A ((\frac{z_1}{w_1}) k   +  (\lambda-\widetilde{\lambda})B)}{( n_3+\frac{1}{w_3}(H_kf_4 + (\frac{\lambda^*}{\lambda}) (F_k +B )f^*) d_3)(  n_3+\frac{1}{w_3}(H_kf_4 + ( \frac{\lambda^*}{\widetilde{\lambda}})(G_k+B) f^*) d_3)}\right).
\end{align*}
\textbf{Step 2: Taylor expansions.} Next we use Taylor expansions once more, this time to treat $\psi_{(\frac{\lambda^*}{\lambda})N}( w_3n_3+H_kf_4d_3)$.  Since $\psi_{(\lambda^*/\lambda)N}(n)$ is a shifted smooth coefficient sequence at scale $(\lambda^*/\lambda)N$, there exist a constant $x_0 \in \mathbb{R}$ and a smooth function $\phi:\mathbb{R} \rightarrow \mathbb{C}$, supported on $[c,C]$, with $|\phi^{(j)}(t)|\ll_j \log(x)^{O(1)}$ and $\psi_{(\lambda^*/\lambda)N}(n)=\phi(\frac{n-x_0}{(\lambda^*/\lambda)N})$.  Since $\psi_{\frac{\Delta^*}{z_1}}(d)$ is a shifted smooth coefficient sequence at scale $\Delta^*/z_1$, there also exist a constant $y_0 \in \mathbb{R}$ and a smooth function $\varphi:\mathbb{R} \rightarrow \mathbb{C}$, supported on $[c,C]$, with $|\varphi^{(j)}(t)|\ll_j \log(x)^{O(1)}$ and $\psi_{\frac{\Delta^*}{z_1}}(d) = \varphi(\frac{d-y_0}{\Delta^*/z_1})$. (Looking at the definition of $W_7(z,k)$, we also note the additional assumptions $x_0 \ll x $ and $y_0 \ll x$.) Recall that $\Delta^* \ll N/(\Lambda x^{5\varepsilon})$ and observe that 
\begin{align*}
\left|\dfrac{H_k f_4 (\Delta^*/(z_1 f^*))}{(\lambda^*/\lambda)N}\dfrac{(d_3 - \frac{y_0}{f^*})}{\Delta^*/(z_1 f^*)} \right| \ll \left|\dfrac{[f_2/f_2^*,f_3/f_3^*] \Delta^* \Lambda}{z_1 N} \right|\ll \left|\dfrac{w_1 \Delta^* \Lambda}{z_1 N} \right| \ll \dfrac{\Delta^* \Lambda}{ N} \ll \dfrac{1}{x^{5\varepsilon}}
\end{align*}
whenever $|f^*d_3 -y_0| \ll \Delta^*/z_1$.  Choosing $J= \lceil \frac{20}{\varepsilon} \rceil$, we may write 
\begin{align} \label{equ:taylorexpansion3}
\psi_{(\frac{\lambda^*}{\lambda})N}( w_3n_3+H_kf_4d_3)&= \sum_{j=0}^J  \dfrac{1}{j!} \left( \frac{H_k f_4 (d_3-\frac{y_0}{f^*})}{(\lambda^*/\lambda) N}  \right)^j \!\! \phi^{(j)}\!\left(\frac{w_3n_3-x_0+(\frac{H_k f_4}{f^*})y_0}{(\lambda^*/\lambda) N}  \right)\\ &+ O_\varepsilon(x^{-100}1_{n_3 \ll x}). \nonumber
\end{align} 
The sequence $\phi^{(j)}((n_3-\frac{x_0}{w_3}+(\frac{H_k f_4}{f^* w_3})y_0)/(\frac{\lambda^*N}{\lambda w_3}))$ is a shifted smooth coefficient sequence in $n_3$ at scale $(\frac{\lambda^*}{\lambda w_3})N$.  The sequence $\frac{1}{j!} \left( \frac{H_k f_4 (d_3-\frac{y_0}{f^*})}{(\lambda^*/\lambda) N}  \right)^j\varphi\left(\frac{ (d_3-\frac{y_0}{f^*})}{\Delta^*/( z_1f^*)}\right)$ is a shifted smooth coefficient sequence in $d_3$ at scale $\frac{\Delta^*}{z_1 f^*}$. We now substitute the RHS of (\ref{equ:taylorexpansion3}) for   $\psi_{(\frac{\lambda^*}{\lambda})N}( w_3n_3+H_kf_4d_3)$ in $\Upsilon_1(f_1,f_2,f_3)$. 
 Write 
  \begin{align*}
&\Upsilon_2(f_1,f_2,f_3,\psi_{(\frac{\lambda^*}{\lambda w_3})N},\psi_{\frac{\Delta^*}{z_1 f^*}}) =\Bigg| \sum_{\substack{d_3 \\  d_3 \equiv (f^*)^{-1}d_\star \!\!\!\!\!\mod (q_0) \\ (d_3, m )=1 }}\!\!\!\!\!\!\!\psi_{\frac{\Delta^*}{z_1 f^*}}(d_3) \!\!\!\!\!\!\!\!\!\!\!\!\!\sum_{\substack{n_3 \\   n_3\equiv - \frac{1}{w_3}(H_kf_4 + (\frac{\lambda^*}{\lambda})F_k f^* )(f^*)^{-1}d_\star + \frac{1}{w_3} (\frac{\lambda^*}{\lambda})n_\star \!\!\!\!\!\mod (q_0) \\ ((\frac{\lambda}{\lambda^*})w_3n_3+ ( (\frac{\lambda}{\lambda^*})H_kf_4 + F_k f^* )d_3, c_1)=1 \\ (( \frac{\widetilde{\lambda}}{\lambda^*})  w_3n_3+ (( \frac{\widetilde{\lambda}}{\lambda^*})H_kf_4 + G_kf^* )d_3, c_1)=1 \\ ((\frac{\lambda}{\lambda^*}) w_3n_3+ ((\frac{\lambda}{\lambda^*})H_kf_4+ F_k f^* +l f^* )d_3, c_2)=1  \\ (( \frac{\widetilde{\lambda}}{\lambda^*})  w_3n_3+ ( ( \frac{\widetilde{\lambda}}{\lambda^*})H_kf_4 + G_kf^* +l f^* )d_3, c_2)=1     }} \!\!\!\!\!\!\!\!\!\!\!\!\!\!\!\!\!\!\!\!\!\!\!\!\!\!\!\!\!\!\!\!\!\!\!\!\!\!\!\!\!\!\!\!\!\! \psi_{(\frac{\lambda^*}{\lambda w_3})N}( n_3)\, (\star)\,\,\,\,\,\Bigg|
 \\
&\mbox{where } (\star) = e_{m}\left(\dfrac{(\frac{1}{w_3})^2(\frac{\lambda^*}{\lambda})(\frac{\lambda^*}{\widetilde{\lambda}})A ((\frac{z_1}{w_1}) k   +  (\lambda-\widetilde{\lambda})B)}{( n_3+\frac{1}{w_3}(H_kf_4 + (\frac{\lambda^*}{\lambda}) (F_k +B )f^*) d_3)(  n_3+\frac{1}{w_3}(H_kf_4 + ( \frac{\lambda^*}{\widetilde{\lambda}})(G_k+B) f^*) d_3)}\right).
\end{align*}
Taking the supremum over all coefficient sequences $\psi_{(\frac{\lambda^*}{\lambda w_3})N}(n_3)$ and $\psi_{\frac{\Delta^*}{z_1 f^*}}(d_3)$ which are shifted smooth at scales $(\frac{\lambda^*}{\lambda w_3})N$ and $\frac{\Delta^*}{z_1 f^*}$, we then  have
\begin{align}\label{chain2}
\Upsilon_1(f_1,f_2,f_3) \ll_\varepsilon \sup_{(\psi_{(\frac{\lambda^*}{\lambda w_3})N},\psi_{\frac{\Delta^*}{z_1 f^*}})}\Upsilon_2(f_1,f_2,f_3,\psi_{(\frac{\lambda^*}{\lambda w_3})N},\psi_{\frac{\Delta^*}{z_1 f^*}}) + O_\varepsilon(x^{-90}).
\end{align}

\vspace{3mm} \textbf{Step 3: Removal of coprimality conditions (Part 2).} Finally, we now remove the remaining coprimality conditions. Again we use Möbius inversion. For some  $e_0 \mid m$, $e_1, e_2 \mid c_1$ and  $e_3, e_4 \mid c_2$, we define the quantity 
 \begin{align*}
&\Upsilon_3(\dots,e_0,e_1,e_2,e_3,e_4) =\Bigg| \sum_{\substack{d_3 \\  d_3 \equiv (f^*)^{-1}d_\star \!\!\!\!\!\mod (q_0) \\  d_3 \equiv 0 \!\!\!\mod (e_0) }}\!\!\!\!\!\!\!\psi_{\frac{\Delta^*}{z_1 f^*}}(d_3)\!\!\!\!\!\!\!\!\!\!\!\!\!\!\sum_{\substack{n_3 \\   n_3\equiv - \frac{1}{w_3}(H_kf_4 + (\frac{\lambda^*}{\lambda})F_k f^* )(f^*)^{-1}d_\star + \frac{1}{w_3} (\frac{\lambda^*}{\lambda})n_\star \!\!\!\!\!\mod (q_0) \\ (\frac{\lambda}{\lambda^*})w_3n_3+ ( (\frac{\lambda}{\lambda^*})H_kf_4 + F_k f^* )d_3 \equiv 0 \!\!\!\mod (e_1) \\ ( \frac{\widetilde{\lambda}}{\lambda^*})  w_3n_3+ (( \frac{\widetilde{\lambda}}{\lambda^*})H_kf_4 + G_kf^* )d_3 \equiv 0 \!\!\!\mod (e_2) \\ (\frac{\lambda}{\lambda^*}) w_3n_3+ ((\frac{\lambda}{\lambda^*})H_kf_4+ F_k f^* +l f^* )d_3 \equiv 0 \!\!\!\mod (e_3)  \\ ( \frac{\widetilde{\lambda}}{\lambda^*})  w_3n_3+ ( ( \frac{\widetilde{\lambda}}{\lambda^*})H_kf_4 + G_kf^* +l f^* )d_3 \equiv 0 \!\!\!\mod (e_4)     }} \!\!\!\!\!\!\!\!\!\!\!\!\!\!\!\!\!\!\!\!\!\!\!\!\!\!\!\! \psi_{(\frac{\lambda^*}{\lambda w_3})N}( n_3)\, (\star) \,\,\,\,\, \Bigg|
 \\
&\mbox{where } (\star) = e_{m}\left(\dfrac{(\frac{1}{w_3})^2(\frac{\lambda^*}{\lambda})(\frac{\lambda^*}{\widetilde{\lambda}})A ((\frac{z_1}{w_1}) k   +  (\lambda-\widetilde{\lambda})B)}{( n_3+\frac{1}{w_3}(H_kf_4 + (\frac{\lambda^*}{\lambda}) (F_k +B )f^*) d_3)(  n_3+\frac{1}{w_3}(H_kf_4 + ( \frac{\lambda^*}{\widetilde{\lambda}})(G_k+B) f^*) d_3)}\right).
\end{align*}
Using the same argument as in Step 1, we then have 
\begin{align}\label{chain3}
\Upsilon_2(f_1,f_2,f_3,\psi_{(\frac{\lambda^*}{\lambda w_3})N},\psi_{\frac{\Delta^*}{z_1 f^*}}) \ll_\varepsilon \sup_{e_0,e_1,e_2,e_3,e_4} (x^\varepsilon)\Upsilon_3(f_1,f_2,f_3,\psi_{(\frac{\lambda^*}{\lambda w_3})N},\psi_{\frac{\Delta^*}{z_1 f^*}},e_0,e_1,e_2,e_3,e_4),
\end{align}
where the supremum is taken over $e_0,e_1,e_2,e_3,e_4$ with $e_0 \mid m$, $e_1, e_2 \mid c_1$ and  $e_3, e_4 \mid c_2$.

\vspace{3mm}
Now recall that $c_1 \mid m$ and $c_2 \mid m$ and $(\lambda,m^*)=(\widetilde{\lambda},m^*)=w_2$. Since $(\lambda,\widetilde{\lambda})=\lambda^*$, then $w_2 \mid \lambda^*$ and $((\frac{\lambda}{\lambda^*})w_3,mc_1c_2)=((\frac{\widetilde{\lambda}}{\lambda^*})w_3,mc_1c_2)=1$. So the summation conditions of $\Upsilon_3$ can be rewritten as follows:
\begingroup
\allowdisplaybreaks
\begin{align} \label{equ:congruencyddd}
&d_3 \equiv \frac{d_\star }{f^*} \mbox{ mod } (q_0) \quad \mbox{ and } \quad  d_3 \equiv 0 \mbox{ mod } (e_0),  \\  
&n_3\equiv - \frac{1}{w_3}\left(H_kf_4 + \left(\frac{\lambda^*}{\lambda}\right)F_k f^* \right)\frac{d_\star}{f^*} + \frac{1}{w_3} \left(\frac{\lambda^*}{\lambda}\right)n_\star  \mbox{ mod } (q_0),\label{equ:congruency0} \\ 
\label{equ:congruency1}
&n_3 \equiv - \frac{1}{w_3}\left( H_kf_4  +\left(\frac{\lambda^*}{\lambda}\right) F_k f^*  \right)d_3\mbox{ mod } \left(e_1\right),\\
\label{equ:congruency2}
&n_3 \equiv - \frac{1}{w_3}\left( H_kf_4  +\left(\frac{\lambda^*}{\widetilde{\lambda}}\right) G_k f^*  \right)d_3\mbox{ mod } \left(e_2\right), \\
\label{equ:congruency3}
&n_3 \equiv - \frac{1}{w_3}\left( H_kf_4  +\left(\frac{\lambda^*}{\lambda}\right)(F_k f^* +l f^*) \right)d_3\mbox{ mod } \left(e_3\right),
\\
\label{equ:congruency4}
&n_3 \equiv - \frac{1}{w_3}\left( H_kf_4  +\left(\frac{\lambda^*}{\widetilde{\lambda}}\right)(G_k f^* +l f^*) \right)d_3\mbox{ mod } \left(e_4\right).
\end{align}
\endgroup
Now denote by $g_1$ the least common multiple of the six integers $q_0$, $e_0$, $e_1$, $e_2$, $e_3$ and $e_4$. Denote by $g_2$ the least common multiple of $q_0$ and $e_0$. Denote by $g_3$ the least common multiple of $q_0$, $e_1$, $e_2$, $e_3$ and $e_4$. Note that $g_1 = [g_2,g_3]$ and $q_0 \mid (g_2,g_3)$.

\vspace{3mm}
\textbf{Step 4: A final partition.} 
 We now partition $\Upsilon_3$ according to two new congruence conditions: 
\begin{align} \label{equ:congruency7}
d_3 \equiv d^\circ \mbox{ mod } \left(\dfrac{g_1}{g_2} \right) \quad \mbox{ and }  \quad n_3 \equiv n^\circ \mbox{ mod } \left(\dfrac{g_1}{g_3} \right). 
\end{align}
The conditions (\ref{equ:congruencyddd}), (\ref{equ:congruency0}), (\ref{equ:congruency1}), (\ref{equ:congruency2}), (\ref{equ:congruency3}) and (\ref{equ:congruency4}) may be incompatible. However, if they are compatible, then  (\ref{equ:congruencyddd}) and the LHS of  (\ref{equ:congruency7}) determine the congruence class of $d_3$ mod $g_1$. Replacing $d_3$ by this congruence class in (\ref{equ:congruency0}), (\ref{equ:congruency1}), (\ref{equ:congruency2}), (\ref{equ:congruency3}) and (\ref{equ:congruency4}), these expressions, together with the RHS of (\ref{equ:congruency7}), in turn determine the congruence class of $n_3$ mod $g_1$. We write $d_3 \equiv d_4$ mod $g_1$ and $n_3 \equiv n_4$ mod $g_1$. The only choices we made were the value of $d^\circ$ and the value of $n^\circ$, so there are at most $(g_1/g_2)(g_1/g_3)$ choices of $d_4$ and $n_4$ mod $g_1$ which contribute to $\Upsilon_3$. Set 
 \begin{align*}
&\Upsilon_4(\dots,d_4,n_4) = \Bigg|\sum_{\substack{d_3 \\  d_3 \equiv d_4  (g_1) }}\sum_{\substack{n_3 \\   n_3 \equiv n_4   (g_1)   }}  \psi_{\frac{\Delta^*}{z_1 f^*}}(d_3)\psi_{(\frac{\lambda^*}{\lambda w_3})N}( n_3)\, (\star) \,\,\,\,\,\Bigg|
 \\
&\mbox{where } (\star) = e_{m}\left(\dfrac{(\frac{1}{w_3})^2(\frac{\lambda^*}{\lambda})(\frac{\lambda^*}{\widetilde{\lambda}})A ((\frac{z_1}{w_1}) k   +  (\lambda-\widetilde{\lambda})B)}{( n_3+\frac{1}{w_3}(H_kf_4 + (\frac{\lambda^*}{\lambda}) (F_k +B )f^*) d_3)(  n_3+\frac{1}{w_3}(H_kf_4 + ( \frac{\lambda^*}{\widetilde{\lambda}})(G_k+B) f^*) d_3)}\right).
\end{align*}
For $\underline{w}=(f_1,f_2,f_3,\psi_{(\frac{\lambda^*}{\lambda w_3})N},\psi_{\frac{\Delta^*}{z_1 f^*}},e_0,e_1,e_2,e_3,e_4)$, we then have
\begin{align}\label{chain4}
\Upsilon_3(\underline{w}) \ll \left(\frac{g_1}{g_2}\right)\left(\frac{g_1}{g_3}\right) \sup_{(d_4,n_4)} \Upsilon_4(\underline{w},d_4,n_4) \ll \left(\frac{g_1}{q_0}\right) \sup_{(d_4,n_4)} \Upsilon_4(\underline{w},d_4,n_4),
\end{align}
where the supremum is taken over all residue classes $d_4$ and $n_4$ mod $g_1$. (In the second upper bound we also used that  $[g_2, g_3/q_0]= g_1$ and $g_2g_3/q_0 \geq g_1$.)

\vspace{3mm}
To conclude the proof, we now simplify the expression for $\Upsilon_4$ a bit further. Choose  $B_1, L \in \mathbb{Z}$ with  
\begin{align*}
&B_1 \equiv \frac{1}{w_3}\left(H_kf_4 + \left(\frac{\lambda^*}{\lambda}\right) (F_k +B )f^*\right)\mbox{ mod } m, \\
&L_1 \equiv \frac{1}{w_3}\left( \left( \frac{\lambda^*}{\widetilde{\lambda}}\right)(G_k+B) f^*- \left(\frac{\lambda^*}{\lambda}\right) (F_k +B )f^*\right)\mbox{ mod } m.
\end{align*}
We recall (\ref{equ:relationshipFG}) to observe  that
 \begin{align*}
 L_1 \equiv \left\{\left(\frac{f^*}{w_3} \right)\left(\frac{\lambda^*}{\lambda }\right)\left(\frac{w_2}{ \widetilde{\lambda}}\right)\right\}\left\{\dfrac{1}{w_2}\left(\left(\frac{z_1}{w_1}\right)k + (\lambda-\widetilde{\lambda})B\right)\right\} \mbox{ mod } m.
 \end{align*}
 The residue class in the first pair of $\{\}$-brackets is primitive. We now choose an integer $A_1$  for which $A_1 \equiv (\frac{w_3}{f^*})(\frac{\widetilde{\lambda}}{w_2})(\frac{1}{w_3})^2(\frac{\lambda^*}{\widetilde{\lambda}})A$ mod $m$. Here $(A_1,m)=1$. The enumerator of the fraction inside $e_m()$ in $\Upsilon_4$ now reads $w_2 A_1 L_1$ and we have $ ( (\frac{z_1}{w_1})k+(\lambda-\widetilde{\lambda})B,m ) = (w_2 A_1 L_1,m)$.
 
 \vspace{3mm} Furthermore, we relabel our  coefficient sequences, writing $\psi_{\Delta_2/z_1}(d_3)=\psi_{\frac{\Delta^*}{z_1 f^*}}(d_3)$ and $\psi_{N_2}(n_3)=\psi_{(\frac{\lambda^*}{\lambda w_3})N}( n_3)$. Here $\Delta_2 \leq \Delta^*$ and $N_2 \leq N$.  We recall that $g_1$ is divisible by $q_0$ and  write $g_1 = q_0 q_3$.  Then  
  \begin{align*}
&\Upsilon_4(\dots,d_4,n_4) =\Bigg| \sum_{\substack{d_3 \\  d_3 \equiv d_4  (q_0q_3) }}\sum_{\substack{n_3 \\   n_3 \equiv n_4   (q_0q_3)  }}  \psi_{\Delta_2/z_1}(d_3)\psi_{N_2}( n_3)  e_{m}\left(\dfrac{w_2A_1L_1}{( n_3+B_1 d_3)(  n_3+(B_1+L_1)d_3)}\right)\Bigg|.
\end{align*}
We now see that $\Upsilon_4(f_1,f_2,f_3,\psi_{(\frac{\lambda^*}{\lambda w_3})N},\psi_{\frac{\Delta^*}{z_1 f^*}},e_0,e_1,e_2,e_3,e_4,d_4,n_4)=\Sigma_8(z,w_8,k)$ for some $w_8 \in W_8(z)$. 
Putting (\ref{chain1}), (\ref{chain2}), (\ref{chain3}) and (\ref{chain4}) together, we obtain the proposed bound
$$
 \Sigma_7(z,w_6,w_7, k) \ll_\varepsilon \sup_{w_8 \in W_8(z)} \left( x^{2\varepsilon} q_3 \right) \Sigma_8(z,w_8,k) + O_\varepsilon(x^{-89}). \eqno \qedhere
$$
\end{proof}

\subsection{Exponential sums}

We have finally arrived at an expression to which the exponential sum bounds of Section~8 of  Polymath~\cite{Polymath:2014:EDZ} can be applied. In particular, we use bound (8.21) of Proposition~8.4 of Polymath~\cite{Polymath:2014:EDZ}, a version of which is stated below: 

\begin{lemma}  \label{lem:polymathprop8.4}
Let $\delta>0$ and  $\varepsilon\in(0,\frac{\delta}{10^{100}})$. Let $\psi_{\Delta}$ and $\psi_N$ be shifted smooth coefficient sequences at scales $\Delta$ and $N$ with $0<\Delta(x), N(x) \ll x^{O(1)}$. For all $x>1$, $m \in \mathbb{N}$ with $m \mid P(x^\delta)$  and $m \ll x^{O(1)}$, and residue classes $A$, $B$, $L$, $\gamma_1$ and $\gamma_2$  mod $m$, we have
\begin{align*}
&\sum_d \sum_n \psi_{\Delta}(d;x) \psi_N(n;x) e_m\left(\dfrac{AL}{(n+Bd+\gamma_1)(n+(B+L)d +\gamma_2)} \right) \\ &\ll_\varepsilon  x^\varepsilon (AL,m) \left(\dfrac{N(x)}{m^{1/2}} + m^{1/2} \right)\left(\dfrac{\Delta(x)}{m^{1/2}}  +  m^{1/2}\right).
\end{align*}
\end{lemma}
Polymath~\cite{Polymath:2014:EDZ} deduced this exponential sum estimate from deep results of 
Deligne on  the Riemann hypothesis over finite fields.

\vspace{3mm}
This lemma is not quite applicable to $\Sigma_8(z,w_8,k)$ yet, since $\Sigma_8(z,w_8,k)$ also involves the congruence conditions $d \equiv d^\star$ mod $(q_0q_3)$ and $n \equiv n^\star$ mod $(q_0q_3)$. Proposition~8.4 of~\cite{Polymath:2014:EDZ} actually allows for the presence of such congruence conditions, but its stated bounds are not sufficiently strong for our purposes, and it potentially also contains some small mistakes: some factors $q_0$ appear to be in the wrong place. Hence we give a different version of this proposition below. 

\begin{lemma} \label{lem:alternativepolymathprop8.4}
Let $\delta>0$ and  $\varepsilon \in (0,\frac{\delta}{10^{100}})$. Let $\psi_{\Delta}$ and $\psi_N$ be shifted smooth coefficient sequences at scales $\Delta$ and $N$ with $0< \Delta(x), N(x) \ll x^{O(1)}$. For all $x>1$, $m, q \in \mathbb{N}$ with $q \mid m$, $m \mid P(x^\delta)$ and $m \ll x^{O(1)}$, residue classes  $A, B$ and $L$ mod $m$ and residue classes $d^\star$ and $n^\star$ mod $q$, the following is true:
\begin{align*}
&\sum_{\substack{d \\ d \equiv d^\star (q)}} \sum_{\substack{n \\ n \equiv n^\star (q)}} \!\!\!\! \psi_{\Delta}(d;x) \psi_N(n;x) e_m\!\left(\dfrac{AL}{(n+Bd)(n+(B+L)d)} \right) \! \\ &\ll_\varepsilon \!  \frac{x^\varepsilon (AL,m) }{q}\!\left(\dfrac{N(x)}{m^{1/2}} + m^{1/2}\! \right)\!\left(\dfrac{\Delta(x)}{m^{1/2}}  + m^{1/2}\!\right)\!.
\end{align*}
\end{lemma}

\begin{proof}
To remove the congruence conditions, we write $d = d_1 q + d^\star$ and $n = n_1 q+ n^\star$. Then 
\begin{align*}
\Upsilon&=\sum_{\substack{d \\ d \equiv d^\star (q)}} \sum_{\substack{n \\ n \equiv n^\star (q)}} \psi_{\Delta}(d) \psi_N(n) e_m\left(\dfrac{AL}{(n+Bd)(n+(B+L)d)} \right)\\ &= 
\sum_{\substack{d_1 }} \sum_{\substack{n_1 }} \psi_{\Delta}(q d_1 + d^\star) \psi_N(q n_1 + n^\star) e_m\left(\dfrac{AL}{(q n_1 + n^\star+B(q d_1 + d^\star))(q n_1 + n^\star+(B+L)(q d_1 + d^\star))} \right).
\end{align*}
Write $\psi_{\Delta/q}(d_1)=\psi_{\Delta}(q d_1 + d^\star)$ and $\psi_{N/q}(n_1) =\psi_N(q n_1 + n^\star)$. Observe that $\psi_{\Delta/q}(d_1)$ is a shifted smooth coefficient sequence at scale $\Delta/q$ and $\psi_{N/q}(n_1)$ is a shifted smooth coefficient sequence at scale $N/q$. Set $\gamma_1 = n^\star+Bd^\star$ and $\gamma_2 = n^\star + (B+L)d^\star$. Since $m$ is squarefree, we may write $m = qm_1$, where $(q,m_1)=1$. In particular, $q$ has an inverse mod $m_1$. Then note that
\begin{align*}
&e_m\left(\dfrac{AL}{(q n_1 + n^\star+B(q d_1 + d^\star))(q n_1 + n^\star+(B+L)(q d_1 + d^\star))} \right) \\
&= e_{qm_1}\left(\dfrac{AL}{(q n_1 +Bq d_1 + \gamma_1)(q n_1 + (B+L)q d_1 + \gamma_2)} \right) \\ &= e_q\!\left(\dfrac{AL}{m_1(q n_1 \!+\!Bq d_1 \!+\! \gamma_1)(q n_1 \!+\! (B+L)q d_1 \!+\! \gamma_2)} \right) 
e_{m_1}\!\left(\dfrac{AL}{q(q n_1 \!+\!Bq d_1 \!+\! \gamma_1)(q n_1 \!+\! (B+L)q d_1 \!+\! \gamma_2)} \right) \\ &= e_q\!\left(\dfrac{AL}{m_1 \gamma_1\gamma_2} \right) 
e_{m/q}\left(\dfrac{(AL/q^3)}{(n_1 +B d_1 + \gamma_1/q)( n_1 + (B+L) d_1 + \gamma_2/q)} \right).
\end{align*}
Here $e_q(AL/(m_1\gamma_1\gamma_2))$ is a constant which is either equal to zero or has modulus $1$. So 
\begin{align*}
|\Upsilon| \leq \left|
\sum_{\substack{d_1 }} \sum_{\substack{n_1 }} \psi_{\Delta/q}(d_1)\psi_{N/q}(n_1) e_{m/q}\left(\dfrac{(AL/q^3)}{(n_1 +B d_1 + \gamma_1/q)( n_1 + (B+L) d_1 + \gamma_2/q)} \right) \right|.
\end{align*}
Now Lemma~\ref{lem:polymathprop8.4} is applicable to the expression above. We apply the lemma with $\Delta/q$ in the place of $\Delta$, $N/q$ in the place of $N$, $m/q$ in the place of $m$ and $A_1 \equiv AL/q^3$ mod $m/q$ in the place of $AL$ mod $m$. We also note that $(A_1,m/q)=(AL,m/q) \leq (AL,m)$ since $(q,m/q)=1$.  The result follows.
\end{proof}

\subsubsection{Upper bounds on $\Sigma_6$} We now combine Lemma~\ref{lem:adjustments2a}, Lemma~\ref{lem:adjustments2} and  Lemma~\ref{lem:alternativepolymathprop8.4}. 

\begin{lemma}[Bounding $\Sigma_6$]\label{lem:concretebounds} Let $Z_5$, $W_6$ and $\Sigma_6$ be as given in Lemma~{\rm\ref{lem:summary}} and  Lemma~{\rm\ref{lem:adjustments1}}.   Then for every $x>1$, $z \in Z_5(x)$, $w_6 \in W_6(z;x)$  and  $k \in \mathbb{Z}$ with   $|k| \ll  \frac{w_1\Lambda N}{x^{5\varepsilon}\Delta_1}$,   $((\frac{z_1}{w_1})k+(\lambda-\widetilde{\lambda})B,m) \leq T(z;x)$  and
 $w_2 \mid k$, we have
\begin{align*}
 \Sigma_6(z,w_6,k;x) \ll_\varepsilon ((\tfrac{z_1}{w_1})k + (\lambda - \widetilde{\lambda})B,m)\left(\frac{x^{3\varepsilon}}{q_0} \right)\left(\frac{\Delta_1}{\Delta^*}\right)\left(\dfrac{N}{m^{1/2}} + m^{1/2} \right)\!\left(\dfrac{\Delta^*}{ m^{1/2} }  + m^{1/2}\right). 
\end{align*}
\end{lemma}

\begin{proof} 
Let $z \in Z_5$, $w_6 \in W_6(z)$ and $k \in \mathbb{Z}$ with   $|k| \ll  \frac{w_1\Lambda N}{x^{5\varepsilon}\Delta_1}$,   
 $w_2 \mid k$  and $((\frac{z_1}{w_1})k+(\lambda-\widetilde{\lambda})B,m) \leq T(z)$. We first recall from Lemma~\ref{lem:adjustments2a} that 
\begin{align*}
 \Sigma_6(z,w_6,k) \ll_\varepsilon \sup_{w_7 \in W_7(z,k)} \left( \frac{\Delta_1}{\Delta^*}\right) (\Sigma_7(z,w_6,w_7,k) + O_\varepsilon(x^{-90})).
\end{align*} 
Lemma~\ref{lem:adjustments2} then gives 
 \begin{align*}
 \Sigma_7(z,w_6,w_7, k) \ll_\varepsilon \sup_{w_8 \in W_8(z)} \left( x^{2\varepsilon} q_3 \right) \Sigma_8(z,w_8,k) + O_\varepsilon(x^{-89}).
 \end{align*}
 In combination, we have  the bound
 \begin{align} \label{combination1}
 \Sigma_6(z,w_6,k) \ll_\varepsilon \sup_{w_8 \in W_8(z)} \left( \frac{\Delta_1}{\Delta^*}\right) \left( x^{2\varepsilon} q_3 \right) \Sigma_8(z,w_8,k) + O_\varepsilon(x^{-50}).
\end{align} 
Now we recall that for a given $w_8 \in W_8(z)$, 
 \begin{align*}
&\Sigma_8(z,w_8,k) =  \Bigg| \sum_{\substack{d \\  d \equiv d^\star  (q_0q_3) }} \sum_{\substack{n \\   n \equiv n^\star  (q_0q_3)    }}  \psi_{\Delta_2/z_1}(d)\psi_{N_2}(n)  e_{m}\left(\dfrac{w_2A_1L_1}{( n+B_1 d)(  n+(B_1+L_1)d)}\right)\Bigg|,
\end{align*}  
where $
L_1 \equiv  A_2 \left\{\frac{1}{w_2}\left(\left(\frac{z_1}{w_1}\right)k + (\lambda-\widetilde{\lambda})B\right)\right\} \mbox{\rm mod } m$. By Lemma~\ref{lem:alternativepolymathprop8.4}, then 
\begin{align*}
&\Sigma_8(z,w_8,k) \ll_\varepsilon \frac{x^\varepsilon (w_2A_1 L_1,m) }{q_0q_3}\left(\dfrac{N_2}{m^{1/2}} + m^{1/2} \right)\!\left(\dfrac{\Delta_2}{z_1 m^{1/2} }  + m^{1/2}\right)\!.
\end{align*}
Recall that for every $w_8 \in W_8(z)$, we have $(A_1A_2,m)=1$,  $\Delta_2 \leq \Delta^*$ and $N_2 \leq N$. Thus 
\begin{align} \label{combination2}
\Sigma_8(z,w_8,k) \ll_\varepsilon \frac{x^{\varepsilon} ((\frac{z_1}{w_1})k + (\lambda - \widetilde{\lambda})B,m) }{q_0q_3}\left(\dfrac{N}{m^{1/2}} + m^{1/2} \right)\!\left(\dfrac{\Delta^*}{ m^{1/2} }  + m^{1/2}\right). 
\end{align}
Substituting (\ref{combination2}) into (\ref{combination1}), we then get 
\begin{align*}
\Sigma_6(z,w_6,k) \ll_\varepsilon \frac{x^{3\varepsilon} ((\frac{z_1}{w_1})k + (\lambda - \widetilde{\lambda})B,m) }{q_0}\left(\frac{\Delta_1}{\Delta^*}\right)\left(\dfrac{N}{m^{1/2}} + m^{1/2} \right)\!\left(\dfrac{\Delta^*}{ m^{1/2} }  + m^{1/2}\right). 
\end{align*}
 This is what we wanted to show. 
\end{proof}

\subsubsection{Sums of greatest common divisors} 

Next we will bound $\Sigma_5$. Our bounds on $\Sigma_6$ involve gcds and so we need to bound sums of common divisors. We use the lemma below.

\begin{lemma}[gcd sums]\label{lem:gcdsums}
Let $\varepsilon>0$,  $m \in \mathbb{N}$ and $ K, T>0$ with $T \geq K$. Let $A, B \in \mathbb{Z}$ with $A\neq 0$. Then 
\begin{align*}
\sum_{\substack{  |k| \leq K \\ (Ak+B, m ) \leq T }} (Ak+B,m) \ll_\varepsilon m^\varepsilon (A,m) T.
\end{align*}  
\end{lemma}

\begin{proof}
We have 
\begin{align*}
\sum_{\substack{ |k| \leq K \\ (Ak+B, m ) \leq T }} (Ak+B,m) &\leq \sum_{d \mid m}\sum_{\substack{ |k| \leq K \\ (Ak+B, m ) \leq T \\ d \mid Ak+B}}  d   \leq 
\sum_{\substack{d \mid m \\ d \leq T}}\sum_{\substack{ |k| \leq K  \\ d \mid Ak+B}} d.
\end{align*} 
Now write $w(d) = (A,d)$. We have $Ak +B \equiv 0$ mod $d$ if and only if $(w(d) \mid B$ and $\frac{A}{w(d)}k+\frac{B}{w(d)} = \ell \frac{d}{w(d)}$ for some $\ell \in \mathbb{Z})$ if and only if $(w(d) \mid B$ and $k \equiv - (\frac{A}{w(d)})^{-1} (\frac{B}{w(d)})$ mod $(\frac{d}{w(d)}))$. For $d$ with $ d \mid m$ and $d \leq T$, this implies that
\begin{align*}
\sum_{\substack{ |k| \leq K  \\ d \mid Ak+B}} d \ll d\left(\dfrac{K}{(d/w(d))} + 1 \right) \ll d  \left(\dfrac{(A,m)K}{d} + \dfrac{T}{d}\right) \ll (A,m)T
\end{align*} 
and hence we get 
$$
\sum_{\substack{ |k| \leq K \\ (Ak+B, m ) \leq T }} (Ak+B,m) \ll \sum_{\substack{d \mid m \\ d \leq T}} (A,m)T \ll_\varepsilon m^\varepsilon (A,m)T. \eqno \qedhere
$$
\end{proof}

\begin{lemma}[Finding a sequence  $(\xi_k) \in \Xi$]\label{lem:gcdboundssigma6}
 Let $Z_5$, $\Sigma_6$ and $\Xi$ be as given in Lemma~{\rm\ref{lem:summary}} and  Lemma~{\rm\ref{lem:adjustments1}}.  

\vspace{3mm} Let $x>1$ and $z \in Z_5(x)$. Then for all $w_6 \in W_6(z;x)$ and all $k \in \mathbb{Z}$ with   $|k| \ll  \frac{w_1\Lambda N}{x^{5\varepsilon}\Delta_1}$,   
 $w_2 \mid k$  and $((\frac{z_1}{w_1})k+(\lambda-\widetilde{\lambda})B,m) \leq T(z;x)$, we have
\begin{align*}
 &\Sigma_6(z,w_6,k;x) \ll_\varepsilon \xi_k \left(\frac{x^{4\varepsilon} (w_2,m) T(z;x)}{q_0} \right)\left(\frac{\Delta_1}{\Delta^*}\right)\left(\dfrac{N}{m^{1/2}} + m^{1/2} \right)\!\left(\dfrac{\Delta^*}{ m^{1/2} }  + m^{1/2}\right), \\ 
 &\mbox{where } \xi_k = ((\tfrac{z_1}{w_1})k + (\lambda - \widetilde{\lambda})B,m) (x^\varepsilon (w_2,m) T(z;x))^{-1}.
\end{align*}
Furthermore, $(\xi_k) \in \Xi(z;x)$. 
\end{lemma}

\begin{proof}
Fix $z \in Z_5$. We recall from Lemma~\ref{lem:concretebounds} that for every $w_6 \in W_6(z)$ and  every $k \in \mathbb{Z}$ with   $|k| \ll  \frac{w_1\Lambda N}{x^{5\varepsilon}\Delta_1}$,   
 $w_2 \mid k$  and $((\frac{z_1}{w_1})k+(\lambda-\widetilde{\lambda})B,m) \leq T(z) = \max\{\tfrac{1}{(w_2,m)}, \tfrac{1}{H}\} \tfrac{x^{\delta+100\varepsilon}H^2 N}{(v_1,v_2)\Delta_1}$, 
\begin{align*}
 \Sigma_6(z,w_6,k) \ll_\varepsilon ((\tfrac{z_1}{w_1})k + (\lambda - \widetilde{\lambda})B,m)\left(\frac{x^{3\varepsilon}}{q_0} \right)\left(\frac{\Delta_1}{\Delta^*}\right)\left(\dfrac{N}{m^{1/2}} + m^{1/2} \right)\!\left(\dfrac{\Delta^*}{ m^{1/2} }  + m^{1/2}\right). 
\end{align*} 
Recall $\Lambda \ll \frac{x^{\delta+5\varepsilon}H^2}{w_1(v_1,v_2)}$. Hence observe that $T(z) = \max\{\tfrac{1}{(w_2,m)}, \tfrac{1}{H}\} \tfrac{x^{\delta+100\varepsilon}H^2 N}{(v_1,v_2)\Delta_1}  \geq 
 \tfrac{x^{\delta+100\varepsilon}H^2 N}{w_2(v_1,v_2)\Delta_1} 
\geq \frac{w_1\Lambda N}{ x^{5\varepsilon}w_2\Delta_1}$ when $x$ is large. 
We also recall that $(z_1,m)=1$ and so $((\frac{z_1}{w_1})w_2,m)=(w_2,m)$. Hence by Lemma~\ref{lem:gcdsums}, 
\begin{align*}
\sum_{\substack{|k| \ll  \frac{w_1\Lambda N}{x^{5\varepsilon}\Delta_1} \\ w_2 \mid k \\ ((\frac{z_1}{w_1})k+(\lambda-\widetilde{\lambda})B,m) \leq T(z)}} \!\!\!\!\!\!\!\!\!\!\!\!\!\!\! ((\tfrac{z_1}{w_1})k + (\lambda - \widetilde{\lambda})B,m) = \!\!\!\!\!\!\!\!\!\!\!\!\!\!\!\!\!\!\!\!\! \sum_{\substack{|k_1| \ll  \frac{w_1\Lambda N}{ x^{5\varepsilon}w_2\Delta_1} \\  \\ ((\frac{z_1}{w_1})w_2k_1+(\lambda-\widetilde{\lambda})B,m) \leq T(z)}} \!\!\!\!\!\!\!\!\!\!\!\!\!\! ((\tfrac{z_1}{w_1})w_2k_1 + (\lambda - \widetilde{\lambda})B,m) \ll_\varepsilon x^\varepsilon (w_2,m) T(z). 
\end{align*}
Setting $\xi_k = ((\tfrac{z_1}{w_1})k + (\lambda - \widetilde{\lambda})B,m) (x^\varepsilon (w_2,m) T(z))^{-1}$, we thus get $\sum\limits^\ast \xi_k \ll_\varepsilon 1$ and $(\xi_k ) \in \Xi(z)$. This is what we wanted to show. 
\end{proof}

\subsubsection{Proof of Proposition $1$}
 
To prove Proposition~\ref{prop:typeI/IIestimate}, it remains to combine Lemma~\ref{lem:summary}, Lemma~\ref{lem:adjustments1} and Lemma~\ref{lem:gcdboundssigma6}.  

\begin{lemma}[Small $\gamma$] \label{lem:proofofprop1}
Let $\omega, \delta >0$. Let $\varepsilon \in (0, 10^{-100}\delta)$. Let $\alpha$ and $\beta$ be coefficient sequences at scales $M$ and $N$ with $x \ll M(x)N(x) \ll x$.   Assume $\beta$ has the Siegel-Walfisz property. 

\vspace{3mm}
Write $N(x)=x^{\gamma(x)}$ and suppose that for all $x>1$,
\begin{align*}
\max\left\{\frac{1}{4}+12 \omega+4\delta+100\varepsilon, \, 32\omega + 10\delta+400\varepsilon\right\}\leq \gamma(x)  \leq  \frac{1}{2}-4\omega-2\delta-50\varepsilon.  
\end{align*}
Then for $x>1$ and for every $C >0$ and every $a \in \mathbb{Z}$,
\begin{align} \label{equ:equiestimate}
\sum_{\substack{d \leq x^{1/2+2\omega} \\ q \mid P(x^\delta) \\ (q,a)=1}}\Bigg| \sum_{n \equiv a(q)} (\alpha \star \beta)(n;x) - \dfrac{1}{\phi(q)} \sum_{(n,q)=1} (\alpha \star \beta)(n;x)\Bigg| \ll_{C,\varepsilon}  \dfrac{x}{\log(x)^C}.
\end{align}
\end{lemma}

\begin{proof} Lemma~\ref{lem:gcdboundssigma6} tells us that 
for every $z \in Z_5$, there exists $(\xi_k) \in \Xi(z)$ so that for all $w_6 \in W_6(z)$ and  $k \in \mathbb{Z}$ with   $|k| \ll  \frac{w_1\Lambda N}{x^{5\varepsilon}\Delta_1}$,   
 $w_2 \mid k$  and $((\frac{z_1}{w_1})k+(\lambda-\widetilde{\lambda})B,m) \leq T(z)$, 
\begin{align*}
 &\Sigma_6(z,w_6, k) \ll_\varepsilon \xi_k \left(\frac{x^{4\varepsilon} (w_2,m) T(z)}{q_0} \right)\left(\frac{\Delta_1}{\Delta^*}\right)\left(\dfrac{N}{m^{1/2}} + m^{1/2} \right)\!\left(\dfrac{\Delta^*}{ m^{1/2} }  + m^{1/2}\right).
\end{align*}
Suppose the following bound is satisfied for every $z \in Z_5$: 
\begin{align} \label{wts}
 \left(\frac{x^{4\varepsilon}  T(z)}{q_0} \right)\left(\frac{\Delta_1}{\Delta^*}\right)\left(\dfrac{N}{m^{1/2}} + m^{1/2} \right)\!\left(\dfrac{\Delta^*}{m^{1/2} }  + m^{1/2}\right) \ll  \min\left\{\frac{1}{x^{\delta+10\varepsilon}H^3},\frac{1}{H^4}\right\}\frac{q_0(q_0,\ell)  N^2}{ x^{27\varepsilon}} .
\end{align}
Then by Lemma~\ref{lem:adjustments1}, $\sup_{z \in Z_5}\Sigma_5(z) \ll_\varepsilon  \min\{\tfrac{(w_2,m)}{x^{\delta+10\varepsilon}H^3}, \tfrac{ (w_2,m)}{H^4}\}\tfrac{q_0^2(q_0,\ell)^2 N^2 }{x^{27\varepsilon}} $. However, by Lemma~\ref{lem:summary}, this upper bound  implies~(\ref{equ:equiestimate}). 

\vspace{3mm}
Recall that $T(z) = \max\{\tfrac{1}{(w_2,m)}, \tfrac{1}{H}\} \tfrac{x^{\delta+100\varepsilon}H^2 N}{(v_1,v_2)\Delta_1} \leq \tfrac{x^{\delta+100\varepsilon}H^2 N}{(v_1,v_2)\Delta_1} $ and $\Delta^* \leq N$.   To prove (\ref{wts}) and hence (\ref{equ:equiestimate}), it thus remains to verify the following three bounds:
\begin{align} \label{wts1}
&\dfrac{N\Delta^*}{ m}  \ll  \left(\frac{\Delta^*}{\Delta_1}\right)\min\left\{\frac{1}{x^{\delta+10\varepsilon}H^3},\frac{1}{H^4}\right\}\frac{(v_1,v_2)q_0^2(q_0,\ell)  N \Delta_1}{ x^{\delta+131\varepsilon}H^2}, \\ \label{wts2}
 &N \ll  \left(\frac{\Delta^*}{\Delta_1}\right)\min\left\{\frac{1}{x^{\delta+10\varepsilon}H^3},\frac{1}{H^4}\right\}\frac{(v_1,v_2)q_0^2(q_0,\ell)  N \Delta_1}{ x^{\delta+131\varepsilon}H^2 }, \\ \label{wts3}
&m \ll  \left(\frac{\Delta^*}{\Delta_1}\right)\min\left\{\frac{1}{x^{\delta+10\varepsilon}H^3},\frac{1}{H^4}\right\}\frac{(v_1,v_2)q_0^2(q_0,\ell)  N \Delta_1}{ x^{\delta+131\varepsilon}H^2}.
\end{align}
Recall  that $\tfrac{N}{x^{\delta+55\varepsilon}H^2}\ll \Delta_1 \ll \tfrac{N}{x^{55\varepsilon}H^2}$, $\tfrac{RQ^2H}{q_0 (v_1,v_2)\Delta_1} \ll m \ll \tfrac{x^{\delta}RQ^2H}{q_0 (v_1,v_2)\Delta_1}$,  $\Delta^* = \min\{\frac{N}{x^{5\varepsilon}\Lambda }, \Delta_1\}$, $\Lambda \ll \frac{1}{w_1(v_1,v_2)}x^{\delta+5\varepsilon}H^2 $,  $ RQ^2=x^{-\varepsilon} q_0 M H  $, $H = x^\varepsilon RQ^2 (q_0 M)^{-1}$ and $H \ll \tfrac{x^{4\omega+\delta+7\varepsilon}}{q_0}$.

\vspace{3mm}
We begin with (\ref{wts1}). This condition rearranges to
\begin{align} \label{wts11}
\left( \frac{ x^{\delta+131\varepsilon}}{(v_1,v_2)q_0^2(q_0,\ell)  }\right) \dfrac{\max\left\{x^{\delta+10\varepsilon}H^5,H^6\right\}}{ m}  \ll  1. 
\end{align}
Note that $m \gg \tfrac{RQ^2H}{q_0 (v_1,v_2)\Delta_1} \gg \tfrac{x^{55\varepsilon}RQ^2H^3}{q_0(v_1,v_2)N} \gg \tfrac{x^{54\varepsilon}  M H^4}{(v_1,v_2)N}$. 
Substituting lower bounds on $m$ and upper bounds on  $H$, the LHS of (\ref{wts11}) is bounded  by 
\begin{align*} 
\left( \frac{ x^{\delta+131\varepsilon}}{(v_1,v_2)q_0^2(q_0,\ell)  }\right) \dfrac{\max\left\{x^{\delta+10\varepsilon}H^5,H^6\right\}}{ m}  \ll   x^{-1+\delta+77\varepsilon}N^2 \max\left\{x^{\delta+10\varepsilon}H,H^2\right\} \ll  x^{-1+8\omega+3\delta+100\varepsilon}  N^2. 
\end{align*}
But we assumed $\gamma \leq \frac{1}{2}-4\omega-2\delta-50\varepsilon$. Hence (\ref{wts11}) and (\ref{wts1}) hold. 

\vspace{3mm}  
Next we consider (\ref{wts2}). Rearranging, we find that this condition holds if
\begin{align*}
\frac{x^{\delta+131\varepsilon}}{\Delta^*}\max\left\{x^{\delta+10\varepsilon}H^5,H^6\right\} \ll  1.
\end{align*}

 We have $\Delta^*=\min\{\frac{N}{x^{5\varepsilon}\Lambda }, \Delta_1\} \geq \frac{N}{x^{\delta+55\varepsilon}H^2}$ and $H \ll \tfrac{x^{4\omega+\delta+7\varepsilon}}{q_0}$. Hence (\ref{wts2}) holds if the following is true: 
\begin{align*}
x^{32\omega+10\delta+400\varepsilon}  \ll N.
\end{align*}
It thus suffices to assume that $\gamma \geq 32\omega + 10 \delta +400\varepsilon$. 

\vspace{3mm}
Now we consider the final condition, (\ref{wts3}). Recall that  $m \ll \frac{x^\delta RQ^2 H}{(v_1,v_2)\Delta_1}\ll \frac{q_0 x^\delta M H^2}{(v_1,v_2)\Delta_1}$. So (\ref{wts3}) is true if 
$$\frac{ x^{1+2\delta+131\varepsilon } }{N^2\Delta^*\Delta_1} \max\left\{x^{\delta+10\varepsilon}H^7,H^8\right\}\ll  1.$$
 Using $
 H \ll \tfrac{x^{4\omega+\delta+7\varepsilon}}{q_0}$, we see that (\ref{wts3}) holds if the following two inequalities are satisfied: 
 \begin{align} \label{twofinalinequalities}
 \frac{x^{1+32\omega+  10\delta+200\varepsilon}   }{N^2} \left(\frac{ x^{5\varepsilon}\Lambda}{\Delta_1N}\right)  \ll 1 \quad \mbox{ and } \quad  \frac{x^{1+32\omega+  10\delta+200\varepsilon}   }{N^2} \left(\frac{1}{\Delta_1^2}\right) \ll  1.
\end{align}
Using  $\Lambda \ll x^{\delta+5\varepsilon} H^2$ and  $\Delta_1 \gg \frac{N}{x^{\delta+55\varepsilon}H^2}$, we observe
\begin{align*}
 \frac{x^{1+32\omega+  10\delta+200\varepsilon}   }{N^2} \max\left\{\left(\frac{ x^{5\varepsilon}\Lambda}{\Delta_1N}\right),  \left(\frac{1}{\Delta_1^2}\right)\right\}  \ll  \frac{x^{1+48\omega+  16\delta+400\varepsilon}   }{N^4}.
\end{align*}
Hence (\ref{twofinalinequalities}) and (\ref{wts3}) hold when $\gamma \geq 1/4+12\omega+4\delta+100\varepsilon$. 
  Therefore (\ref{equ:equiestimate}) is true. 
\end{proof}

\textit{Proof of Proposition~$1$ $($Final Steps$):$} We consider a coefficient sequence $\alpha$ at scale $N$ and a coefficient sequence $\beta$ at scale $M$, where $N(x)M(x) \asymp x$ and $N(x)= x^{\gamma(x)}$ with $\gamma(x) \in [1/2-\sigma,1/2]$. We also assume that $\beta$ has the Siegel-Walfisz property and that the following three inequalities are satisfied:
\begin{align} \label{assumptions}
\begin{cases}
\,\,\, 72\omega + 24 \delta < 1, \\
\,\,\, 48 \omega + 16\delta + 4\sigma <1, \\ 
\,\,\, 64 \omega +20\delta+2\sigma<1.
\end{cases}
\end{align}
The second and third inequality rearrange to $1/4+12 \omega + 4\delta <1/2-\sigma $ and $32 \omega +10\delta<1/2-\sigma$. But $\gamma(x) \geq 1/2-\sigma$. Choosing $\varepsilon$ sufficiently small compared to $\delta$, we then find that 
$$\max\{1/4+12 \omega + 4\delta+100\varepsilon, 32 \omega +10\delta+400\varepsilon\}\leq \gamma(x).$$  
By Lemma~\ref{lem:proofofprop1} we thus have the desired equidistribution estimate (\ref{equ:equiestimate}) if also $\gamma(x) \leq 1/2-4\omega-2\delta-50\varepsilon$. 

\vspace{3mm}
 On the other hand, part (ii) and (iv) of Theorem~2.8 of Polymath~\cite{Polymath:2014:EDZ} tell us that (\ref{equ:equiestimate}) holds if $68\omega + 14 \delta<1$ and $\gamma(x)>1/4+14\omega + 4 \delta    $. Observe that $1/4+14\omega + 4 \delta < 1/2-4\omega-2\delta-50\varepsilon$ if $72\omega + 24\delta +200\varepsilon < 1$.  Hence under the assumptions (\ref{assumptions}), the equidistribution estimate (\ref{equ:equiestimate}) holds for every $\gamma(x) \in [1/2-\sigma,1/2]$.  This  concludes the proof of Proposition~\ref{prop:typeI/IIestimate}. \qed

\vspace{3mm}
(Note: In (\ref{equ:equiestimate}), the implied constant depends on $\varepsilon$, but in the proof of Proposition~\ref{prop:typeI/IIestimate}, the choice of $\varepsilon$ depends on $\omega$ and $\delta$, so we have  replaced $\ll_{C,\varepsilon}$ by $\ll_{C, \omega,\delta}$.)

\section{Harman's sieve} \label{sec:harman}

In this section we apply Harman's sieve, constructing a function $\rho(n;x) $ which is a minorant of the prime indicator function $1_{\mathbb{P}}(n)$ for each $x$,  and which has exponent of distribution $\theta = 0.5253 =1/2 + 2 \cdot 0.01265$. 

\subsection{Motivation} 

We first explain our choice of $\theta$. In the arguments of Maynard~\cite{Maynard:2015:SGP} and Polymath~\cite{Polymath:2014:VSS}, bounds on $H_m$ were derived using the exponent of distribution of $\Lambda(n)$. This exponent of distribution was determined by decomposing $\Lambda(n)$ via the Heath-Brown identity and treating various convolutions of coefficient sequences separately. Certain bad coefficient sequences, corresponding to convolutions of $5$ sequences of length about $x^{0.2}$,  formed a barrier to improving the exponent of distribution of $\Lambda(n)$ further. The idea behind a use of Harman's sieve is to discard these bad convolutions: 
Baker and Irving~\cite{Baker:2017:BIP} noticed that if we replace $\Lambda(n)$ by $\log(x) \rho_\omega(n;x)$, where $\rho_\omega(n;x)$ is a minorant for $1_{\mathbb{P}}(n)$ which satisfies 
\begin{align} 
\sum_{\substack{ q \leq x^{1/2 +2\omega } \\ q \mid P(x^\delta) \\ (q,a)=1}} \Bigg| \sum_{\substack{ n \in [x,2x] \\ n \equiv a (q) }} \rho_\omega(n;x) - \dfrac{1}{\phi(q)}\sum_{\substack{ n \in [x,2x] \\ (n,q)=1 }} \rho_\omega(n;x)   \Bigg| \ll_{A,  \delta} \dfrac{x}{\log(x)^A}
\end{align} 
and  $\sum_{x \leq n \leq 2x} \rho_w(n;x) = (1- c_1(\omega) +o(1))\frac{x}{\log(x)}$, then we still have a bound of the form
$$H_m \ll_\epsilon  \exp\left(\frac{(1+\epsilon)m}{(1/4+\omega)(1-c_1(\omega))} \right).$$
 The quantity $c_1(\omega)$ describes the loss incurred by discarding bad convolutions and increases as $\omega$ increases, since fewer and fewer equidistribution estimates are available. (For $\omega = 1/80$ thus $c_1(\omega)=0$.)
Harman's sieve then gives an improvement over a direct treatment of $\Lambda(n)$ whenever the loss $c_1(\omega)$, incurred by discarding bad cases, does not outweigh the gain of a larger $\omega$, so that $(1/4+\omega)(1-c_1(\omega))>1/4+1/80$. 

\subsubsection{Limitations}
 Somewhat unusually, this particular application of Harman's sieve has a discontinuity, present in the same form both in Baker and Irving's and in our work: The minorant $\rho_\omega(n;x)$ is constructed using the Buchstab identity, and in the process one must consider whether or not $\alpha_1 \star \alpha_2 \star \alpha_3$ with $\alpha_i(n)=1_{\mathbb{P}}(n)$ on $[N_i,2N_i]$  and $N_1N_2N_3 \asymp x$ has exponent of distribution $1/2+2\omega$. To do so, each $\alpha_i$ is first decomposed further via the Heath-Brown identity and we consider $\alpha_i = \psi_i \star \beta_i$, where  $\psi_i$ and $\beta_i$ have support $[A_i,2A_i]$ and $[B_i,2B_i]$, $A_iB_i \approx N_i$, $B_i \leq x^{0.1}$ and $\psi_i(n)=1$ on $[A_i, 2A_i]$. Good equidistribution estimates are available if $A_1A_2A_3$ is large or if the product of some subset of $\{A_1,A_2,A_3,B_1,B_2,B_3\}$ is close to $x^{0.5}$. In particular, if $\omega =  1/79 $ then  we either need that $A_1A_2A_3 \geq x^{0.8925 }$ or that a product of $A_1,A_2,A_3,B_1,B_2,B_3$ is contained in $[x^{0.405},x^{0.595}]$, and one of these two cases is always satisfied by $\psi_1 \star \psi_2 \star \psi_3 \star \beta_1 \star \beta_2 \star \beta_3$.

\vspace{3mm} However, as soon as $\omega$ is increased slightly beyond $ 1/79 $, our equidistribution estimates do not cover the case where $A_1A_2A_3$ is very close to $x^{0.8925}$ -- the three smooth factors are too short for Polymaths's Type III estimate, but too long for our Type I/II estimate. When constructing $\rho_\omega(n;x)$, this forces us to discard many products of three primes in  Buchstab's identity, leading to a sudden jump where $c_1(\omega)$ becomes much larger as $\omega$ increases beyond $  1/79 $.

\vspace{3mm} This discontinuity  determines the choice $\omega = 0.01265 < 1/79$. Here we still have $c_1(\omega) < 2 \cdot 10^{-5}$ and thus obtain the proposed bound $H_m = O(\exp(3.8075m))$.

\subsection{Type I, II and III information}

We now record the available Type I, II and III information for the case $\omega =0.01265 \approx \frac{1}{79}$. (Here we are using Type I/II in the sense of Harman rather than Polymath.) 
 
 \begin{lemma} \label{lem:typei/ii/ii}
Let $\delta = 10^{-10}$. Suppose $f: \mathbb{N} \times (1,\infty) \rightarrow \mathbb{C}$  satisfies  one of the following three conditions:  
\begin{enumerate} 
\item[{\rm (I)}]  $f(n;x) = (\alpha \star \beta)(n;x)$, where $\alpha$ is  coefficient sequence at scale $M$ and $\beta$ is a \textbf{smooth} coefficient sequence at scale $N$ with $M(x) N(x) \asymp x$. Additionally, $N(x) = x^{\gamma(x)}$ with $
\gamma(x) \geq  0.33856-\delta.
$
\item[{\rm (II)}]  $f(n;x) = (\alpha \star \beta)(n;x)$, where $\alpha$ is  coefficient sequence at scale $M$, $\beta$ is a  coefficient sequence at scale $N$ with $M(x) N(x) \asymp x$, and $\alpha$ and $\beta$ both have the Siegel--Walfisz property.   Additionally, $N = x^{\gamma(x)}$ with 
$0.40481 -\delta \leq \gamma(x) \leq  0.59519+\delta$. 
\item[{\rm (III)}] $f(n;x) = (\alpha \star \psi_1 \star \psi_2 \star \psi_3)(n;x)$, where $\alpha$ is a coefficient sequence at scale $M$ and $\psi_1$, $\psi_2$ and $\psi_3$ are  \textbf{smooth} coefficient sequences at scales $N_1$, $N_2$ and $N_3$ with $M(x) N_1(x) N_2(x) N_3(x) \asymp x$.   Additionally, $x^{0.19-\delta} \leq N_1(x), N_2(x), N_3(x) \leq x^{0.405+\delta}$ and $N_i(x)N_j(x) \geq x^{0.595-\delta}$ for $i \neq j$.
\end{enumerate} 
Then for any  $x>1$, $a\in \mathbb{Z}$ and   $A>0$, we have
\begin{align*}
\sum_{\substack{q \leq x^{0.5253 } \\ q \mid P(x^\delta) \\ (q,a)=1}} \Bigg| \sum_{\substack{n \in [x,2x] \\ n \equiv a (q)}} f(n;x) - \dfrac{1}{\phi(q)} \sum_{\substack{ n \in [x,2x] \\ (n,q)=1}}f(n;x) \Bigg| \ll_A \dfrac{x}{\log(x)^A}.
\end{align*}
In particular, $f$ has exponent of distribution $0.5253$ to smooth moduli.
\end{lemma}

\begin{proof}
For option (I), we use Lemma~5 of Baker and Irving~\cite{Baker:2017:BIP}. We take $\eta = 0.0138825$ to get $\theta = 0.5253$ and observe that $0.33856 > 199/600+119 \eta/240$. For option (II), we use Proposition~\ref{prop:typeI/IIestimate}. We take $\omega =0.01265$ and observe that $0.40481  > \max\left\{\frac{1}{4}+12\omega, \, 32\omega \right\}$.  Finally, for option (III) we use part (v) of Theorem~2.8 of Polymath~\cite{Polymath:2014:EDZ}. Again we take $\omega = 0.01265$ and observe that $0.595 > 1/2+1/18+28\omega/9$, while $0.405< 1/2-1/18-28\omega/9$ and $0.19>2(1/18+28\omega/9)$. 
\end{proof}

\subsection{Good sifted sets}
From here on out, we  follow the arguments of Baker and Irving~\cite{Baker:2017:BIP} closely, effectively only updating their proofs with better Type~II information. 

\vspace{3mm}
 As a first step, we record a version of Baker and Irving's Lemma~7. Recall that \begin{align*}
\psi(n,y) = \begin{cases} 
\,\, 1  \quad \mbox{ if } p \mid n \rightarrow p \geq y, \\ \,\, 0 \quad \mbox{ otherwise. }
\end{cases}
\end{align*}

\begin{lemma}\label{lem:goodsiftedsets}  Write $p_i =x^{\alpha_i}$. Let $\zeta = 1-0.33856-0.40481 = 0.25663$ and $\lambda = 0.59519-0.40481=0.19038 $. 
Each of the functions given below has exponent of distribution $0.5253$ to smooth moduli. 
\vspace{3mm}

\begin{enumerate}
\item[\rm{(1)}] \hspace{45mm} $\begin{aligned}\theta_1(n;x)=\psi(n, x^\lambda).\end{aligned}$ \vspace{2mm}
\item[\rm{(2)}] \hspace{45mm} $\begin{aligned}\theta_2(n;x)=\sum_{\substack{ n=p_1n_2 \\ 0.19038 \leq \alpha_1 < 0.40481}} \psi(n_2, x^\lambda).\end{aligned}$ \vspace{2mm}
\item[\rm{(3)}] \hspace{45mm} $\begin{aligned}\theta_3(n;x)=\sum_{\substack{ n=p_1p_2n_3 \\ 0.19038 \leq \alpha_2 <  \alpha_1 < 0.40481 \\ \alpha_1 + \alpha_2 < 0.40481}} \psi(n_3, x^\lambda).\end{aligned}$ \vspace{2mm}
\item[\rm{(4)}] \hspace{45mm} $\begin{aligned}\theta_4(n;x)=\sum_{\substack{ n=p_1p_2n_3 \\ 0.19038 \leq \alpha_2 <  \alpha_1 < 0.40481 \\ \alpha_1 + \alpha_2 > 0.59519,\, \alpha_2 <\zeta}} \psi(n_3, x^\lambda).\end{aligned}$ \vspace{2mm}
\item[\rm{(5)}] \hspace{45mm} $\begin{aligned}\theta_5(n;x)=\sum_{\substack{ n=p_1p_2p_3n_4 \\ 0.19038 \leq \alpha_3 < \alpha_2 <  \alpha_1 < 0.40481 \\ \alpha_1 + \alpha_2 < 0.40481,\, \alpha_3 <\zeta}} \psi(n_4, x^\lambda).\end{aligned}$ \vspace{2mm}
\item[\rm{(6)}] \hspace{45mm} $\begin{aligned}\theta_6(n;x)=\sum_{\substack{ n=n_1p_2p_3p_4 \\ 0.19038 \leq \alpha_4 < \alpha_3 <  \alpha_2 < 0.40481 \\ \alpha_3 + \alpha_4 < 0.40481,\, \alpha_2 <\zeta, \, \alpha_4 \geq \alpha_3}} \psi(n_1, x^\lambda).\end{aligned}$ \vspace{2mm}
\end{enumerate}
\end{lemma}
\begin{proof}
This follows from the exact same arguments as Lemma~7 of~\cite{Baker:2017:BIP}. We have simply updated the available Type II information. 

\vspace{3mm} More precisely, we use Lemma~14 of Baker and Weingartner~\cite{Baker:2014:TDI} with $\alpha = 0.40481$, $\beta = 0.59519-0.40481=0.19038 = \lambda$, so that $\alpha + \beta
 = 0.59519$, and $M = x^{1-0.33856} = x^{0.66144 }$. In the notation of that lemma we must then take $R < x^{0.40481}$ and $S< Mx^{-\alpha} = x^{0.25663} = x^\zeta$.  After verifying that the sequences which appear in the proof of Lemma~14 of~\cite{Baker:2014:TDI} have the Siegel-Walfisz property, we  obtain that $\theta_k(n)$ has exponent of distribution $0.5253$ to smooth moduli if the following is true: We can partition $\{1, \dots, j\}$ into $I$ and $J$ so that $\sum_{i \in I} \alpha_i <0.40481$ and $\sum_{i \in J} \alpha_i < \zeta =0.25663$. (Here $j$ is the number of $p_i$ which appear in $\theta_k(n)$. In (6) we partition $\{2,3,4\}$ rather than  $\{1, \dots, j\}$.)
 
 \vspace{3mm} For (1), we  simply take $I=J=\emptyset$. In (2),  we choose $I=\{1\}$ and $J=\emptyset$. In (3) we take $I = \{1\}$ and $J = \{2\}$. We use that $\alpha_2 \leq (\alpha_1+\alpha_2)/2 < 0.40481/2 <\zeta$. In (4) we also take $I = \{1\}$ and $J = \{2\}$. For (5), we take $I = \{1,2\}$ and $J =\{3\} $. Finally, in (6) we use $I=\{3,4\}$ and $J=\{2\}$. 
\end{proof}

\subsection{Buchstab decompositions} \label{ssec:buchstab}

We now decompose the prime indicator function using Buchstab iterations. Write $\lambda = 0.19038$ and $\zeta = 0.25663$ and set $p_i = x^{\alpha_i}$. For $n \in [x,2x] \cap \mathbb{N}$, we have
\begingroup
\allowdisplaybreaks
\begin{align}
1_{\mathbb{P}}(n) &= \psi(n,(3x)^{1/2}) \nonumber \\ &= \psi(n,x^\lambda) - \sum_{\substack{ n = p_1n_2 \\ \lambda \leq \alpha_1 < 0.40481}} \psi(n_2, x^\lambda)- \sum_{\substack{ n = p_1n_2 \\ 0.40481 \leq \alpha_1 < 0.5 }} \psi(n_2, p_1) \label{equ:buchstab1}\\
&+ \sum_{\substack{ n = p_1p_2n_3 \\ \lambda \leq \alpha_2< \alpha_1 < 0.40481 \\ \alpha_1+\alpha_2 < 0.40481}} \!\!\!\! \psi(n_3, x^\lambda) + \!\!\!\!\!\!\!\!\!\!\! \!\!\!\sum_{\substack{ n = p_1p_2n_3 \\ \lambda \leq \alpha_2 < \alpha_1 < 0.40481 \\ 0.40481 \leq \alpha_1+\alpha_2 \leq  0.59519}} \!\!\!\!\!\!\!\!\!\!\!\!\psi(n_3, p_2)+ \!\!\!\!\!\!\! \sum_{\substack{ n = p_1p_2n_3 \\ \lambda \leq \alpha_2<\alpha_1 < 0.40481 \\ \alpha_1+\alpha_2 > 0.59519 \\ \alpha_2 < \zeta}}\!\!\!\!\!\! \psi(n_3, x^\lambda) + \!\!\!\! \sum_{\substack{ n = p_1p_2n_3 \\ \lambda \leq \alpha_2 < \alpha_1 < 0.40481 \\ \alpha_1+\alpha_2 > 0.59519 \\ \alpha_2 \geq \zeta}}\!\!\!\!\!\! \psi(n_3, p_2) \label{equ:buchstab2}\\
&- \sum_{\substack{ n = p_1p_2p_3n_4 \\ \lambda \leq \alpha_3 < \alpha_2 < \alpha_1 < 0.40481 \\ \alpha_1+\alpha_2 < 0.40481 \\ \alpha_3 < \zeta}} \!\!\!\! \psi(n_4, x^\lambda) - \sum_{\substack{ n = p_1p_2p_3n_4 \\ \lambda \leq \alpha_3 < \alpha_2 < \alpha_1 < 0.40481 \\ \alpha_1+\alpha_2 < 0.40481 \\ \alpha_3 \geq \zeta}} \!\!\!\! \psi(n_4, p_3) + \sum_{\substack{ n = p_1p_2p_3p_4n_5 \\ \lambda \leq \alpha_4 < \alpha_3 < \alpha_2 < \alpha_1 < 0.40481 \\ \alpha_1+\alpha_2 < 0.40481 \\ \alpha_3 < \zeta}} \!\!\!\! \psi(n_5, p_4) \label{equ:buchstab3}\\ 
&- \sum_{\substack{ n = p_1p_2p_3n_4 \\ \lambda \leq \alpha_3 <\alpha_2<\alpha_1 < 0.40481 \\ \alpha_1+\alpha_2 > 0.59519 \\ \alpha_2 < \zeta}}\!\!\!\!\!\! \psi(n_4, p_3). \label{equ:buchstab4}
\end{align}
\endgroup
We may apply part (1) of Lemma~\ref{lem:goodsiftedsets} to the first function in (\ref{equ:buchstab1}) and part (2) to the second function in (\ref{equ:buchstab1}). Further, the third function can be treated directly, using our Type II information, recorded in Lemma~\ref{lem:typei/ii/ii}. The first and third function in (\ref{equ:buchstab2}) can be treated via part (3) and part (4) of Lemma~\ref{lem:goodsiftedsets}, respectively. On the other hand, the second function is treated directly, using Type II information from Lemma~\ref{lem:typei/ii/ii}.  This leaves only the fourth expression. For (\ref{equ:buchstab3}), part (5) of Lemma~\ref{lem:goodsiftedsets} can be applied to the first expression. If $\alpha_3 \geq \zeta$, then $\alpha_1+\alpha_2 \geq 2 \cdot 0.25663 > 0.40481$. Hence the second sum is zero. This leaves only the third function in (\ref{equ:buchstab3}). Finally, for the moment, we leave the function in (\ref{equ:buchstab4}) unchanged.

\vspace{3mm} Overall we have found that there exists a function $\theta_0(n;x)$ which has exponent of distribution $0.5253$ to smooth moduli and which satisfies that for all $n \in [x,2x] \cap \mathbb{N}$,
\begin{align*}
1_{\mathbb{P}}(n) = \theta_0(n;x) + \!\!\!\!\!\!\!\!\!\!\sum_{\substack{ n = p_1p_2n_3 \\ 0.19038 \leq \alpha_2 < \alpha_1 < 0.40481 \\ \alpha_1+\alpha_2 > 0.59519 \\ \alpha_2 \geq 0.25663 }} \!\!\!\!\!\!\!\!\!\!\psi(n_3, p_2) + \!\!\!\!\!\!\!\!\!\!\sum_{\substack{ n = p_1p_2p_3p_4n_5 \\ 0.19038\leq \alpha_4 < \alpha_3 < \alpha_2 < \alpha_1 < 0.40481 \\ \alpha_1+\alpha_2 < 0.40481 \\ \alpha_3 < 0.25663}} \!\!\!\!\!\!\!\!\!\!\!\!\!\! \psi(n_5, p_4)- \!\!\!\!\!\!\!\!\!\!\sum_{\substack{ n = p_1p_2p_3n_4 \\ 0.19038 \leq \alpha_3 <\alpha_2<\alpha_1 < 0.40481 \\ \alpha_1+\alpha_2 > 0.59519 \\ \alpha_2 < 0.25663 }}\!\!\!\!\!\!\!\!\!\!\!\!\!\!\!\! \psi(n_4, p_3).
\end{align*}

\subsection{3 prime factors}
Next we inspect the first sum in the above expression for $1_{\mathbb{P}}(n)$ more closely. This corresponds to Lemma~10 of~\cite{Baker:2017:BIP}.

\begin{lemma} \label{lem:3primes}
For $x>1$ and $n \in [x,2x] \cap \mathbb{N}$, define the function
\begin{align*}
\Gamma(n;x)=\sum_{\substack{ n = p_1p_2n_3 \\ 0.19038 \leq \alpha_2 < \alpha_1 < 0.40481 \\ \alpha_1+\alpha_2 > 0.59519 \\ \alpha_2 \geq 0.25663}} \psi(n_3, p_2).
\end{align*}
Then $\Gamma(n;x)$ has exponent of distribution  $0.5253$ to smooth moduli. 
\end{lemma}

\begin{proof}
Consider $n=p_1p_2n_3$ with $n \in [x,2x]$ and  $0.19038 \leq \alpha_2 < \alpha_1 < 0.40481$ and  $\alpha_2 \geq 0.25663$. Suppose that $n_3$ has no prime factor less than $p_2$.  Assume for a contradiction that $n_3$ is not prime. Then $n$ has at least 4 prime factors, all of which are of size at least $x^{0.25663}$. But then $n$ is certainly larger than $2x$, assuming that $x$ is large. We work with $n\leq 2x$, so this gives a contradiction. Hence it suffices to consider $n=p_1p_2p_3$.  We now look at the sum 
\begin{align*}
\Gamma_0(n;x)=\sum_{\substack{ n = p_1p_2p_3 \\ 0.19038 \leq \alpha_2 < \alpha_1 < 0.40481 \\ \alpha_1+\alpha_2 > 0.59519 \\ \alpha_3 \geq \alpha_2 \geq 0.25663 }} 1.
\end{align*}
Applying the Heath-Brown identity like in Section~3 of Polymath~\cite{Polymath:2014:EDZ}, we find that $\Gamma_0(n;x)$ has exponent of distribution $0.5253$ to smooth moduli provided that each coefficient sequence $\beta \in \mathcal{B}$, with $\mathcal{B}$ as described below, has exponent of distribution $0.5253$ to smooth moduli.

\vspace{3mm} The set $\mathcal{B}$ denotes the set of coefficient sequences 
\begin{align*}
\beta=\beta_1 \star \dots \star \beta_{k}
\end{align*}
with $k \leq 40$ which satisfy the following properties:
\begin{enumerate}
\item Each $\beta_i$ is a coefficient sequence, located at some scale $N_i$, with $\prod_{ i =1}^k N_i(x) \asymp x$ and $\log_x(N_i(x))$ constant.   Furthermore, if $S \subseteq \{1, \dots, k\}$, then $\bigstar_{i \in S} \beta_i$ is a coefficient sequence at scale $\prod_{i \in S} N_i$.
\item If $N_i(x) \geq x^{0.1}$, then $\beta_i$ is smooth at scale $N_i$. 
\item If $S \subseteq \{1, \dots, k\}$ with $\prod_{i \in S} N_i(x) \geq x^{\delta} = x^{10^{-10}}$, then $\bigstar_{i \in S} \beta_i$ has the Siegel-Walfisz property. 
\item There exits a partition $I_1, I_2, I_3$ of $\{1, \dots, k\}$ such that $x^{\gamma_j}=\prod_{i \in I_j} N_i(x)$ satisfy  $0.19038-\frac{\delta}{10} \leq \gamma_1, \gamma_2 \leq 0.40481+\frac{\delta}{10}$, $\gamma_1 + \gamma_2 \geq 0.59519 - \frac{\delta}{10}$, $\gamma_2 \geq 0.25663-\frac{\delta}{10}$ and $\gamma_1, \gamma_3 \geq \gamma_2 - \frac{\delta}{10}$. 
\end{enumerate}

\vspace{3mm} It remains to show that the functions in $\mathcal{B}$ have the desired exponent of distribution.   As a first step, partition $\{1, \dots, k\}$ into sets $J_1, \dots, J_j$ such that $x^{\mu_s}=\prod_{i \in J_s} N_i(x)$ satisfies $\mu_s \geq 0.05$ for all but at most one $s \in \{1,\dots j\}$ and such that $J_s$ is a singleton whenever $\mu_s \geq 0.1$. In particular, $\bigstar_{i \in J_s} \beta_i$ is smooth when $\mu_s \geq 0.1$.  We also assume $\mu_1 \geq \dots \geq \mu_j$. 

\vspace{3mm} If there exits $s$ with $\mu_s \geq 0.33856-\delta$, then we may use the  Type I information available in  Lemma~\ref{lem:typei/ii/ii} to deduce that $\beta$ has exponent of distribution $0.5253$ to smooth moduli. If there exists $K \subseteq \{1, \dots, j\}$ with $\sum_{s \in K} \mu_s \in [0.40481-\delta,0.59519+\delta]$, we use Type II information to deduce  the same. Hence we may assume that neither of these two statements holds. 

\vspace{3mm} Now consider $\mu_2+\mu_3$. Since $\sum_{s \in K} \mu_s \not\in [0.40481-\delta,0.59519+\delta]$, we either have  $\mu_2+\mu_3 \geq 0.59519+ \delta$ or $\mu_2+\mu_3 \leq 0.40481-\delta$. Assume first that  $\mu_2+\mu_3 \geq 0.59519+ \delta$. Then also  $\mu_1+\mu_2, \mu_1+\mu_3 \geq 0.59519+ \delta$. We certainly must have $\mu_1 \leq 0.405+\delta$, since otherwise $\mu_1 \geq 0.59519+\delta$ and $\mu_1+\mu_2+\mu_3>1$.  Hence if also $\mu_3 \geq 0.19$, then $\bigstar_{i \in J_s} \beta_i$ is smooth for $s \leq 3$, and by part (III) of Lemma~\ref{lem:typei/ii/ii}, the sequence
\begin{align*}
(\bigstar_{i \in J_1} \beta_i) \star (\bigstar_{i \in J_2} \beta_i) \star (\bigstar_{i \in J_3} \beta_i) \star (\bigstar_{i \not\in J_1\cup J_2 \cup J_3} \beta_i).
\end{align*} 
has exponent of distribution $0.5253$ to smooth moduli. (The set $\{1,\dots,j\}\setminus (J_1 \cup J_2 \cup J_3)$ may be empty. In that case we replace the final (empty) convolution by the delta function.)  
 Thus we may assume that $\mu_3 <0.19$. But this implies $\mu_2 \geq 0.40519 +\delta$, contradicting our earlier assumptions. This concludes the treatment of the case $\mu_2+\mu_3 \geq 0.59519+ \delta$. 
 
 \vspace{3mm} Finally, we consider the case $\mu_2+\mu_3 \leq 0.40481-\delta$. We have $\mu_2 + \dots + \mu_j \geq 1-0.33856=0.66144$, so if $\mu_4 \leq 0.59519-0.40481=0.19038$, then there exists $K$ with $\sum_{s \in K} \mu_s \in [0.40481-\delta,0.59519+\delta]$, contradicting our earlier assumption. Similarly, since $\mu_4 \leq \mu_3 \leq 0.21$, we have $\mu_2 + \mu_3 +\mu_5 + \dots +\mu_j \geq 1-0.33856-0.21 =0.45144 > 0.40481$, and thus may also assume that $\mu_5 \geq 0.19038$.  Notice that $\mu_1, \dots, \mu_5 > 0.1$ and thus $J_1, \dots, J_5$ are singletons. In particular, for each $i \leq 5$, there exists $l \in \{1,\dots, k\}$ with  $N_l(x) = x^{\mu_i} \geq x^{0.19038}$. We also have $\mu_i + (\mu_6 + \dots + \mu_j)\leq \mu_1 + (\mu_6 + \dots + \mu_j) \leq 1-4\cdot 0.19038 =0.23848 $ for $1 \leq i \leq 5$.

\vspace{3mm} 
  Looking at condition (4) of $\mathcal{B}$, we have $\gamma_1, \gamma_2, \gamma_3 \geq 0.25663-\delta$. Each of $\{N_i: i \in I_1\}$ and $\{N_i: i \in I_2\}$ and $\{N_i: i \in I_3\}$ must then contain at least two choices of $x^{\mu_i}$ with $1 \leq i \leq 5$. This is because $\mu_i + ( \mu_6+ \dots + \mu_j) < 0.25663-\delta$ for $1 \leq i \leq 5$. But $2 \cdot 3 = 6 >5$, so there are not enough $x^{\mu_i}$ available to satisfy this requirement. Hence condition (4) tells us that $\mu_2 + \mu_3 \leq 0.40481-\delta$ is impossible unless there exists $K$ with $\sum_{s \in K} \mu_s \in [0.40481-\delta,0.59519+\delta]$. This concludes the proof. 
\end{proof}

\subsection{5 prime factors}

Next we consider  products of $5$ primes. This corresponds to Lemma~11 of~\cite{Baker:2017:BIP}. 

\begin{lemma} \label{lem:5primes}
For $x>1$ and  $n \in [x,2x]\cap \mathbb{N}$, define the function 
\begin{align*}
\Gamma(n;x)=\sum_{\substack{ n = p_1p_2p_3p_4n_5 \\ 0.19038 \leq \alpha_4 < \alpha_3 < \alpha_2 < \alpha_1 < 0.40481 \\ \alpha_1+\alpha_2 < 0.40481 \\ \alpha_3 < 0.25663}}  \psi(n_5, p_4). 
\end{align*}
If $x$ is sufficiently large, we can write $\Gamma(n;x) = \Gamma_1(n;x) + \Gamma_2(n;x)$, where $\Gamma_1(n;x)$ has exponent of distribution $0.5253$ to smooth moduli and $\Gamma_2(n;x)$ is non-negative with 
\begin{align*}
\sum_{n \in [x,2x]} \Gamma_2(n;x) = \sum_{n \in [x,2x]} \sum_{\substack{ n = p_1p_2p_3p_4n_5 \\ 0.19038\leq \alpha_4 < \alpha_3 < \alpha_2 < \alpha_1 < 0.40481 \\ \alpha_1+\alpha_2 < 0.40481 \\ \alpha_2+\alpha_3 +\alpha_4 > 0.59519}}  \psi(n_5, p_4)\leq \frac{3x}{10^7\log(x)}.
\end{align*}
\end{lemma}

\begin{proof} The inequality $\alpha_3 < 0.25663$ automatically follows from $\alpha_1+\alpha_2 < 0.40481$ and $\alpha_3 < \alpha_2$, so we may remove it. We now split $\Gamma(n;x)$ up into 
 \begin{align*}
 \Gamma(n;x) = \Gamma_1(n;x) + \Gamma_2(n;x)=   \sum_{\substack{ n = p_1p_2p_3p_4n_5 \\ 0.19038\leq \alpha_4 < \alpha_3 < \alpha_2 < \alpha_1 < 0.40481 \\ \alpha_1+\alpha_2 < 0.40481 \\ \alpha_2+\alpha_3 +\alpha_4 \leq 0.59519}} \!\!\!\!\!\! \psi(n_5, p_4) + \sum_{\substack{ n = p_1p_2p_3p_4n_5 \\ 0.19038\leq \alpha_4 < \alpha_3 < \alpha_2 < \alpha_1 < 0.40481 \\ \alpha_1+\alpha_2 < 0.40481 \\ \alpha_2+\alpha_3 +\alpha_4 > 0.59519}} \!\!\!\!\!\! \psi(n_5, p_4).
 \end{align*}
 Notice here that $\alpha_2 + \alpha_3 + \alpha_4 \geq 3 \cdot 0.19038   > 0.40481$. Hence our Type II  information tells us that $\Gamma_1(n;x)$ has exponent of distribution $0.5253$ to smooth moduli.  On the other hand, if $\alpha_2+\alpha_3 +\alpha_4 > 0.59519$, then $\alpha_1 + \alpha_2+\alpha_3 +3\alpha_4 > (4/3) \cdot 0.59519 + 2 \cdot 0.19038 >1$. Hence if $n_5$ contributes positively to $\Gamma_2(n;x)$, then $n_5$ is prime. Using standard techniques, and in particular the prime number theorem for short intervals, we find that $\sum_{n \in [x,2x]} \Gamma_2(n;x)$ equals 
\begin{align*}
\dfrac{x(1+o(1))}{\log(x)}\int^{0.40481}_{0.19038}\int^{\min\{\alpha_1, 0.40481-\alpha_1\}}_{0.19038}\int^{\alpha_2}_{0.19038}\int^{\min\{\alpha_3,(1-\alpha_1-\alpha_2-\alpha_3)/2\}}_{\max\{0.19038, 0.59519-\alpha_2-\alpha_3\}} \frac{1}{\alpha_1\alpha_2\alpha_3\alpha_4(1-\sum_{i=1}^4\alpha_i)} d\underline{\alpha}.
\end{align*} 
In particular, we have $\sum_{n \in [x,2x]} \Gamma_2(n;x) \leq \frac{3x}{10^7\log(x)}$ when $x$ is large. 
\end{proof}

\subsection{4 prime factors} It remains to look at products of $4$ primes. This corresponds to Lemma~13 of~\cite{Baker:2017:BIP}, and we closely follow the steps of that proof.

\begin{lemma}\label{lem:4primes}
For $x>1$ and $n \in [x,2x]\cap \mathbb{N}$, define the function 
\begin{align*}
\Gamma(n;x)=\sum_{\substack{ n = p_1p_2p_3n_4 \\ 0.19038 \leq \alpha_3 <\alpha_2<\alpha_1 < 0.40481 \\ \alpha_1+\alpha_2 > 0.59519 \\ \alpha_2 < 0.25663}} \psi(n_4, p_3). 
\end{align*}
If $x$ is sufficiently large, we can write $\Gamma(n;x) = \Gamma_1(n;x) - \Gamma_2(n;x)$, where $\Gamma_1(n;x)$ has exponent of distribution $0.5253$ to smooth moduli and $\Gamma_2(n;x)$ is non-negative with 
\begin{align*}
\sum_{n \in [x,2x]} \Gamma_2(n;x)  = \sum_{n \in [x,2x]} \sum_{\substack{ n = p_2p_3p_4p_5p_6 \\ 0.19038 \leq \alpha_2, \alpha_3, \alpha_4, \alpha_5, \alpha_6 \leq 0.23848  \\ \alpha_2, \alpha_4 > \alpha_3  \\ \alpha_2+\alpha_4 < 0.40481\\ \alpha_3 + \alpha_4 + \alpha_5 >0.59519 \\   \alpha_6 \geq \alpha_5 }} 1 \leq \frac{0.00001x}{\log(x)}.
\end{align*}
\end{lemma}

\begin{proof}
Notice that $\alpha_1+\alpha_2 +3\alpha_3 > 0.59519 + 3 \cdot 0.19038 >1$. Hence if $n_4$ contributes to $\Gamma(n;x)$, then $n_4$ is prime. Notice also that $\alpha_1 + \alpha_2 > 0.59519$ and $\alpha_2 < 0.25663$ imply $\alpha_1 > \alpha_2$. Further, $\alpha_i \geq 0.19038$ implies $\alpha_1 < 1- 3 \cdot 0.19038 < 0.59519$. Hence we may write  
\begin{align*}
\Gamma(n;x)= \Upsilon_1(n;x)+\Upsilon_2(n;x)= \sum_{\substack{ n = p_1p_2p_3p_4 \\ 0.19038 \leq \alpha_3 <\alpha_2 < 0.40481 \\ \alpha_1+\alpha_2 > 0.59519 \\ \alpha_2 < 0.25663 \\ \alpha_4 \geq \alpha_3}} 1 - \sum_{\substack{ n = p_1p_2p_3p_4 \\ 0.19038 \leq \alpha_3 <\alpha_2 < 0.40481 \\ \alpha_1+\alpha_2 > 0.59519 \\ \alpha_2 < 0.25663 \\ \alpha_4 \geq \alpha_3 \\ \alpha_1 \in [0.40481,0.59519]}} 1.
\end{align*}
Using the available Type II information, we immediately see that the second sum has exponent of distribution $0.5253$ to smooth moduli.  On the other hand, $\Upsilon_1(n;x)$ may be rewritten as 
\begin{align*}
\Upsilon_1(n;x)= \sum_{\substack{ n = p_1p_2p_3p_4 \\ 0.19038 \leq \alpha_3 <\alpha_2 < 0.40481 \\ \alpha_3+\alpha_4 < 0.40481 \\ \alpha_2 < 0.25663 \\ \alpha_4 \geq \alpha_3}} 1.
\end{align*}
First applying reversal of roles to $p_1$, and then using the Buchstab identity, we get 
\begin{align*}
\Upsilon_1(n;x)&= \sum_{\substack{ n = n_1p_2p_3p_4 \\ 0.19038 \leq \alpha_3 <\alpha_2 < 0.40481 \\ \alpha_3+\alpha_4 < 0.40481 \\ \alpha_2 < 0.25663 \\ \alpha_4 \geq \alpha_3}} \psi(n_1,(\tfrac{2x}{p_2p_3p_4})^{1/2}) \\
&= \sum_{\substack{ n = n_1p_2p_3p_4 \\ 0.19038 \leq \alpha_3 <\alpha_2 < 0.40481 \\ \alpha_3+\alpha_4 < 0.40481 \\ \alpha_2 < 0.25663 \\ \alpha_4 \geq \alpha_3}} \psi(n_1,x^\lambda) -
 \sum_{\substack{ n = p_2p_3p_4p_5n_6 \\ 0.19038 \leq \alpha_3 <\alpha_2 < 0.40481 \\ \alpha_3+\alpha_4 < 0.40481 \\ \alpha_2 < 0.25663 \\ \alpha_4 \geq \alpha_3 \\ 0.19038 \leq \alpha_5 < (1-\alpha_2-\alpha_3-\alpha_4 )/2}} \psi(n_6,p_5). 
\end{align*}
Notice that part (6) of Lemma~\ref{lem:goodsiftedsets} applies to the first sum in the second line, which hence has exponent of distribution $0.5253$ to smooth moduli. Hence we may focus on the second sum. Observe that $6 \cdot 0.19038>1$. Hence if $n_6$ contributes a non-zero amount to $\Upsilon_1(n;x)$, then we must have $n_6$ prime. Let 
\begin{align*}
\Upsilon_3(n;x) = \sum_{\substack{ n = p_2p_3p_4p_5p_6 \\ 0.19038 \leq \alpha_3 <\alpha_2 < 0.40481 \\ \alpha_3+\alpha_4 < 0.40481 \\ \alpha_2 < 0.25663 \\ \alpha_4 \geq \alpha_3 \\ 0.19038 \leq \alpha_5 \leq \alpha_6}} 1. 
\end{align*} 
Our proof is finished if we can show that $\Upsilon_3(n;x)$ is the sum of a function $\Upsilon_4(n;x)$ which has exponent of distribution $0.5253$  to smooth moduli and a non-negative function $\Gamma_2(n;x)$ which satisfies $\sum_{n \in [x,2x]} \Gamma_2(n;x) \leq \frac{0.00001 x}{\log(x)}$. 

\vspace{3mm} The contribution of $(\alpha_2,\dots, \alpha_6)$ with $\sum_{i \in S} \alpha_i \in [0.40481,0.59519]$ to $\Upsilon_3(n;x)$ can be treated via the available Type II information. Notice also that ($\alpha_i \geq 0.19038$ for all $i$)  implies $\alpha_i \leq 1- 4 \cdot 0.19038= 0.23848  < 0.25663$ for $i \in \{2,3,4,5\}$. So we consider 
\begin{align*}
\Upsilon_3^*(n;x) = \sum^*_{\substack{ n = p_2p_3p_4p_5p_6 \\ 0.19038 \leq \alpha_2, \alpha_3, \alpha_4, \alpha_5, \alpha_6 \leq 0.23848  \\ \alpha_2, \alpha_4 > \alpha_3  \\ \alpha_3+\alpha_4 < 0.40481 \\   \alpha_6 \geq \alpha_5 }} 1,
\end{align*}
where $\sum^*$ denotes the sum over $(\alpha_2,\dots, \alpha_6)$ with $\sum_{i \in S} \alpha_i \not\in [0.40481,0.59519]$ for all $S \subseteq \{2, \dots, 6\}$. 
Observe that $\alpha_2 + \alpha_4 < 2 \cdot 0.24 < 0.59519$. Hence we may assume $\alpha_2+\alpha_4 < 0.40481$. This immediately also implies $\alpha_3+\alpha_4 < 0.40481$. Further, $\alpha_3+\alpha_4+\alpha_5 \geq 3 \cdot 0.19038 >0.40481$ and so $\alpha_3+\alpha_4+\alpha_5 \geq 0.59519$. 
\begin{align*}
\Upsilon_3^*(n;x) \leq \Gamma_2(n;x) =  \sum_{\substack{ n = p_2p_3p_4p_5p_6 \\ 0.19038 \leq \alpha_2, \alpha_3, \alpha_4, \alpha_5, \alpha_6 \leq 0.23848  \\ \alpha_2, \alpha_4 > \alpha_3  \\ \alpha_2+\alpha_4 < 0.40481\\ \alpha_3 + \alpha_4 + \alpha_5 >0.59519 \\   \alpha_6 \geq \alpha_5 }} 1.
\end{align*}
Using standard techniques, we find that $\sum_{ n \in[x,2x]}\Gamma_2(n;x)$ is asymptotically equal to 
\begin{align*}
\dfrac{x}{\log(x)}\!\!\int^{0.23848 }_{0.19038}\int^{\min\{0.23848 , \alpha_2\}}_{0.19038}\!\!\!\!\int_{\max\{\alpha_3, 0.19038\}}^{\min\{0.23848 , 0.40481-\alpha_2\}}\!\!\!\!\int^{\min\{0.23848 ,(1-\sum_{i=2}^4\alpha_i)/2\}}_{\max\{0.19038, 0.59519-\alpha_3-\alpha_4\}} \frac{1}{\prod_{i=2}^5\alpha_i(1-\sum_{i=2}^5\alpha_i)} d\underline{\alpha}.
\end{align*} 
In particular, we have $\sum_{n \in [x,2x]} \Gamma_2(n;x) \leq \frac{0.00001x}{\log(x)}$ when $x$ is large. 
\end{proof}

\subsection{Construction of a minorant}

Combining the Buchstab decomposition given in Section~\ref{ssec:buchstab} with Lemma~\ref{lem:3primes}, Lemma~\ref{lem:5primes} and Lemma~\ref{lem:4primes}, we have that for $n \in [x,2x]\cap \mathbb{N}$, 
\begin{align*}
1_{\mathbb{P}}(n) = \theta^*(n;x) + \sum_{\substack{ n = p_1p_2p_3p_4n_5 \\ 0.19038\leq \alpha_4 < \alpha_3 < \alpha_2 < \alpha_1 < 0.40481 \\ \alpha_1+\alpha_2 < 0.40481 \\ \alpha_2+\alpha_3+\alpha_4 >0.59519}} \psi(n_5, p_4)+ \sum_{\substack{ n = p_2p_3p_4p_5p_6 \\ 0.19038 \leq \alpha_2, \alpha_3, \alpha_4, \alpha_5, \alpha_6 \leq 0.23848  \\ \alpha_2, \alpha_4 > \alpha_3  \\ \alpha_2+\alpha_4 < 0.40481\\ \alpha_2 + \alpha_3 + \alpha_5 >0.59519 \\   \alpha_6 \geq \alpha_5 }} 1,
\end{align*}
where $\theta^*(n;x)$ has exponent of distribution $0.5253$ to smooth moduli and where 
\begin{align*}
\sum_{n \in [x,2x]} (1_{\mathbb{P}}(n) - \theta^*(n;x)) \leq  \dfrac{0.0000103x}{\log(x)}
\end{align*}
when $x$ is sufficiently large. Setting $\rho(n;x) =  \theta^*(n;x)$, this concludes the proof of Proposition~\ref{prop:construction}. \qed

\appendix

\section{Proof of Lemma~\ref{lem:summaryofpolymath}} 

We now give a proof of Lemma~\ref{lem:summaryofpolymath}. It follows the arguments of Section~5  and Section~8A  of Polymath~\cite{Polymath:2014:EDZ} and only small changes need to be made. As a first step, we summarize  the first pages of Section~5 of~\cite{Polymath:2014:EDZ}, which can be used largely unchanged.

\begin{lemma}\label{lem:PolymathSection5}
Let $\omega, \delta >0$ and $\varepsilon  \in (0, 10^{-100} \delta)$. Let $\alpha$ and $\beta$ be coefficient sequences at scales $M$ and $N$ with $x \ll M(x)N(x) \ll x$ and $N(x)= x^{\gamma(x)}$, where $\gamma(x) \in (4\omega+2\delta,1/2]$.    Assume $\beta$ has the Siegel-Walfisz property. 

\vspace{3mm} We denote by $Z_0(x)$ the set of tuples  $(Q,R,a,\psi_M,(c_{q,r}) )$ which have the following properties:
\begin{enumerate}[{\rm (i)}]
\item $a$ is an arbitrary integer.
\item  $\psi_M(m)$ is a real, non-negative coefficient sequence  which is smooth at scale $M$.
\item For all $q,r \in \mathbb{N}$, $c_{q,r}$ is a complex number with $|c_{q,r}|=1$.
\item  $Q$ and $R$ are positive real numbers with
\begin{align} 
&x^{-4\varepsilon-\delta}N \ll R \ll x^{-2\varepsilon}N, \label{appendix:R}\\
&x^{1/2-\varepsilon} \ll QR \ll x^{1/2+2\omega+\varepsilon}.  \label{appendix:Q}
\end{align}
\end{enumerate}
We also write $C^*(n) = 1_{\substack{\frac{b_1}{n}  \equiv \frac{b_2}{n+\ell r} ((q_1,q_2))}}
$ for $n$ with  $(n,q_1 r) = (n+\ell r, q_2)=1$  
 and  set $D_0 = \exp(\log(x)^{1/3})$.

\vspace{3mm} 
For $x>1$, $z \in Z_0(x)$ and 
 $b_1, b_2 \in \mathbb{Z}$ with $(b_1b_2,P(x^\delta))=1$, we then define
\begin{align*}
&\Sigma(b_1,b_2;z,x)= \!\!\! \sum_{\substack{r \asymp R \\ (r,a)=1}} \sum_{\substack{ |\ell| \ll \frac{N}{R} \\ \ell \neq 0}} \!\!\! \sum_{\substack{q_1 \asymp Q \\  q_1r \mid P(x^\delta) \\ (q_1,a )=1\\ (q_1,P(D_0))=1}}\sum_{\substack{q_2\asymp Q \\  q_2r \mid P(x^\delta) \\ (q_2,a )=1 \\ (q_2,P(D_0))=1}} \!\!\!\!\!\!\!\!\!\! c_{q_1,r} \overline{c_{q_2,r}} \!\!\!\!\!\!\sum_{\substack{n \\ (n,q_1r)=1 \\ (n+\ell r,q_2)=1}} \!\!\!\!\!\!\!\!C^*(n) \beta(n) \overline{\beta(n+\ell r)} \sum_m \psi_M(m)1_{m \equiv y \, ([q_1,q_2]r)}\\
&\mbox{where } y \mbox{ is a function of }  n, \ell, a, b_1, b_2, r, q_1, q_2 \mbox{ with } y \equiv \tfrac{b_1}{n}(q_1), y \equiv \tfrac{b_2}{n+\ell r}(q_2) \mbox{ and } y \equiv \tfrac{a}{n}(r).
\end{align*}
Suppose that for all $x>1$, $z \in Z_0(x)$ and 
 $b_1, b_2 \in \mathbb{Z}$ with $(b_1b_2,P(x^\delta))=1$,
 $$\Sigma(b_1,b_2;z,x) = X(z,x)+O_{A, \varepsilon}(N^2MR^{-1}\log(x)^{-A}),$$ where $X(z,x)$ is some function that does not depend on $b_1$ and $b_2$, and where the implied constant does not depend on $x$, $z$, $b_1$ and $b_2$. Then for all $x>1$, $a_0 \in \mathbb{Z}$  and $A>0$,
\begin{align} \label{finalresultappendix}
\sum_{\substack{q \leq x^{1/2+2\omega} \\ q \mid P(x^\delta) \\ (q,a_0)=1}}\Bigg| \sum_{n \equiv a_0(q)} (\alpha \star \beta)(n;x) - \dfrac{1}{\phi(q)} \sum_{(n,q)=1} (\alpha \star \beta)(n;x)\Bigg| \ll_{A, \varepsilon}  \dfrac{x}{\log(x)^A}.
\end{align}
\end{lemma}

\begin{proof}
This is essentially a summary of pages 2118 to 2125 of Section~5C of~\cite{Polymath:2014:EDZ}, which concluded with an application of the dispersion method of Linnik. 
\end{proof}

The next step is to complete sums in   $\Sigma(b_1,b_2)$, with the goal of reducing the proof of (\ref{finalresultappendix}) to bounding certain exponential sums. For this purpose we record the following minor amendment of Lemma~4.9 of~\cite{Polymath:2014:EDZ}, which (unlike Lemma~4.9) preserves the smooth coefficients of $h$: 

\begin{lemma}[Completion of sums]\label{lem:amendment4.9}
Let  $\psi_M$ be a shifted smooth coefficient sequence at scale $M$.

Let $x>1$,  
let $I$ be a finite set, let $c_i \in \mathbb{C}$, let $q \in \mathbb{N}$ and let $a_i$ be a residue class mod $q$.  Set 
\begin{align} \label{def:psitophi}
\varphi_H(h;x)= \frac{1}{M(x)}
\sum_m \psi_M(m;x)   e\left(\dfrac{-mh}{q} \right). 
\end{align}
Then for $A > 0$,  $\varepsilon>0$ and $B\geq qM(x)^{-1+\varepsilon}$,
\begin{align*}
\sum_{i \in I} c_i \sum_m \psi_M(m;x) 1_{m =a_i (q)} &=\dfrac{1}{q} \sum_{m \geq 1} \psi_M(m;x) \sum_{i \in I}c_i + \dfrac{M(x)}{q}\sum_{0 < |h| \leq B} \varphi_H(h;x) \sum_{i \in I}c_i  e\left(\dfrac{a_i h}{q} \right) \nonumber \\ &+ O_{A,\varepsilon}\Big(M(x)^{-A} \sum_{i \in I} |c_i|\Big)
. 
\end{align*}
\end{lemma}
\begin{proof}
This is a direct consequence of the proof of Lemma~4.9 of~\cite{Polymath:2014:EDZ}.
\end{proof}

(Note: The implied constant in Lemma~\ref{lem:amendment4.9} does not depend on $x$, $I$, $c_i$, $q$, $a_i$ or $B$.)

\vspace{3mm}

Before we apply Lemma~\ref{lem:amendment4.9} to $\Sigma(b_1,b_2)$,  it is convenient to record one more result of Polymath~\cite{Polymath:2014:EDZ}, which will be used to bound some of the error terms which soon arise: 

\begin{lemma} \label{lem:polymathconvenienterrors}
Let $\omega, \delta >0$. Let $M$, $N$, $Q$ and $R$ be positive real numbers  such that $MN \asymp x$ and such that  conditions {\rm(\ref{appendix:R})} and {\rm(\ref{appendix:Q})} of Lemma~{\rm\ref{lem:PolymathSection5}} hold. Let $a, b_1, b_2 \in \mathbb{Z}$ and let $C^*(n)$ be as defined in Lemma~{\rm\ref{lem:PolymathSection5}}. Write $D_0 = \exp(\log(x)^{1/3})$. Then 
\begin{align*}
\dfrac{M \log(x)^{O(1)}}{R}\sum_{r \asymp R} \sum_{\substack{|\ell| \ll \frac{N}{R} \\ \ell \neq 0}} \sum_{\substack{1 \neq q_0 \ll Q \\ (q_0, P(D_0))=1 }} \mathop{\sum\sum}\limits_{\substack{q_1,q_2 \asymp Q \\ q_1r, q_2r \mid P(x^\delta) \\ (q_1,q_2)=q_0}} \dfrac{1}{[q_1,q_2]} \!\!\!\!\! \sum_{\substack{n \\ (n,q_1r)=1 \\ (n+\ell r,q_2)=1}} \!\!\!\!\!(\tau(n)\tau(n+\ell r))^{O(1)}C^*(n) \ll_A \dfrac{MN^2}{R\log(x)^{A}}.
\end{align*}
\end{lemma}

\begin{proof}
This is shown on page~2126 of~\cite{Polymath:2014:EDZ}.
\end{proof}

Applying Lemma~\ref{lem:amendment4.9} and Lemma~\ref{lem:polymathconvenienterrors} to Lemma~\ref{lem:PolymathSection5}, we obtain Lemma~\ref{lem:amendmenttheorem5.8}, stated below. This is an amended version of Theorem~5.8 of~\cite{Polymath:2014:EDZ}, in which the smooth coefficients of variable $h$ have been preserved. 
\begin{lemma}[Exponential sum estimates]\label{lem:amendmenttheorem5.8} 
Let $\omega, \delta >0$ and $\varepsilon  \in (0, 10^{-100} \delta)$. Let $\alpha$ and $\beta$ be coefficient sequences at scales $M$ and $N$ with $x \ll M(x)N(x) \ll x$ and $N(x)= x^{\gamma(x)}$, where $\gamma(x) \in (4\omega+2\delta,1/2]$.    Assume $\beta$ has the Siegel-Walfisz property. 

\vspace{3mm} We denote by $Z_0^*(x)$ the set of tuples  $(Q,R,q_0,a,b_1,b_2,\ell,\psi_M,)$ which have the following properties:
\begin{enumerate}[{\rm (i)}]
\item $Q, R \in (0,\infty)$ satisfy conditions {\rm(\ref{appendix:R})} and {\rm(\ref{appendix:Q})} of Lemma~{\rm\ref{lem:PolymathSection5}}.
\item  $q_0 \in \mathbb{N}$  with $q_0 \ll Q$, $q_0 \mid P(x^\delta)$ and $(q_0, P(D_0))=1$, where $D_0 = \exp(\log(x)^{1/3})$.
\item $a,  b_1, b_2 \in \mathbb{Z}$ with $(q_0,ab_1b_2)=1$. 
\item  $\ell \in \mathbb{Z}$ with $0 \neq |\ell| \ll \frac{N}{R}$.
\item  $\psi_M(m)$ is a real, non-negative coefficient sequence  which is smooth at scale $M$.
\end{enumerate}
We also write $C_0(n)  = 1_{\substack{\frac{b_1}{n}  \equiv \frac{b_2}{n+\ell r} (q_0)}}
$ and  $H(q_0) = \frac{x^\varepsilon RQ^2}{q_0 M}$. 
As  in Lemma~{\rm\ref{lem:amendment4.9}} $($taking  $q = rq_0q_1q_2)$, we set
$$\varphi_H(h)= \frac{1}{M}
\sum_m \psi_M(m)   e\left(\dfrac{-mh}{rq_0q_1q_2} \right). $$
For $x>1$ and  $z \in Z_0^*(x)$, we then define
\begin{align*}
&\Sigma^*(z;x)= \sum_{\substack{ r \asymp R \\ (r,ab_1b_2)=1}}  \mathop{\sum\sum}\limits_{\substack{q_1,q_2 \asymp Q/q_0\\ (q_1,q_2)=1 \\ (q_1q_2,ab_1b_2)=1\\ q_0q_1r, q_0q_2r \mid P(x^\delta)
}} \Bigg| \sum_{0<|h| \leq H(q_0)} \!\!\! \varphi_H(h)\!\!\sum_{\substack{n \\ (n,rq_0q_1 )=1 \\ (n+\ell r, q_0q_2)=1 }} \!\!\!\!C_0(n)\beta(n)\overline{\beta(n+\ell r)}  \Phi_\ell(h,n,r,q_0,q_1,q_2) \Bigg|,\\
&\mbox{where } \Phi_\ell(h,n,r,q_0,q_1,q_2) =  e_r\left(\dfrac{ah}{nq_0q_1q_2}\right)e_{q_0q_1}\left(\dfrac{b_1h}{nrq_2}\right)e_{q_2}\left(\dfrac{b_2h}{(n+\ell r)rq_0q_1}\right).
\end{align*}

\vspace{3mm}
Suppose that $\Sigma^*(z;x) \ll_{\varepsilon} (q_0,\ell) (x^\varepsilon q_0^2)^{-1}RQ^2N $ for $x>1$ and $z \in Z^*_0(x)$. Then {\rm(\ref{finalresultappendix})} is satisfied.
\end{lemma} 

\begin{proof} We fix some $x>1$ and $z \in Z_0(x)$ and recall that for $b_1,b_2 \in \mathbb{Z}$ with $(b_1b_2,P(x^\delta))=1$, 
\begin{align*}
&\Sigma(b_1,b_2;z,x)= \!\!\! \sum_{\substack{r \asymp R \\ (r,a)=1}} \sum_{\substack{ |\ell| \ll \frac{N}{R} \\ \ell \neq 0}} \!\!\! \sum_{\substack{q_1 \asymp Q \\  q_1r \mid P(x^\delta) \\ (q_1,a )=1\\ (q_1,P(D_0))=1}}\sum_{\substack{q_2\asymp Q \\  q_2r \mid P(x^\delta) \\ (q_2,a )=1 \\ (q_2,P(D_0))=1}} \!\!\!\!\!\!\!\!\!\! c_{q_1,r} \overline{c_{q_2,r}} \!\!\!\!\!\!\sum_{\substack{n \\ (n,q_1r)=1 \\ (n+\ell r,q_2)=1}} \!\!\!\!\!\!\!\!C^*(n) \beta(n) \overline{\beta(n+\ell r)} \sum_m \psi_M(m)1_{m \equiv y \, ([q_1,q_2]r)}\\
&\mbox{where } y \mbox{ is a function of }  n, \ell, a, b_1, b_2, r, q_1, q_2 \mbox{ with } y \equiv \tfrac{b_1}{n}(q_1), y \equiv \tfrac{b_2}{n+\ell r}(q_2) \mbox{ and } y \equiv \tfrac{a}{n}(r).
\end{align*}
 For fixed $r$, $\ell$, $q_1$ and $q_2$, we now apply Lemma~\ref{lem:amendment4.9} to the remaining variables $m$ and $n$ and find that  $\Sigma(b_1,b_2) = \Sigma_0^*(b_1,b_2) + \Sigma_1^*(b_1,b_2) + O_A(N^2MR^{-1}\log(x)^{-A})$, where  $\Sigma_0^*(b_1,b_2)$ and $\Sigma_1^*(b_1,b_2)$ are given by
\begin{align*}
&\Sigma_0^*(b_1,b_2)= \left(\sum_m \psi_M(m) \right)\sum_{\substack{r \asymp R \\ (r,a)=1}} \sum_{\substack{ |\ell| \ll \frac{N}{R} \\ \ell \neq 0}}  \mathop{\sum\sum}\limits _{\substack{q_1, q_2\asymp Q \\  q_1r, q_2r \mid P(x^\delta)\\ (q_1q_2,a)=1 \\ (q_1q_2,P(D_0))=1}}\!\!\!\!\!\! \dfrac{c_{q_1,r} \overline{c_{q_2,r}}}{r[q_1,q_2]} \!\!\!\sum_{\substack{n \\ (n,q_1r)=1 \\ (n+\ell r,q_2)=1}} \!\!\!C^*(n) \beta(n) \overline{\beta(n+\ell r)}, \\
&\Sigma_1^*(b_1,b_2)= \!\!\!\!\! \sum_{\substack{r \asymp R \\ (r,a)=1}} \sum_{\substack{ |\ell| \ll \frac{N}{R} \\ \ell \neq 0}}  \mathop{\sum\sum}\limits _{\substack{q_1, q_2 \asymp Q \\  q_1r, q_2r \mid P(x^\delta) \\ (q_1q_2,a)=1 \\ (q_1q_2,P(D_0))=1}} \!\!\!\!\!\! \dfrac{c_{q_1,r} \overline{c_{q_2,r}}}{r[q_1,q_2]} \!\sum_{0< |h| \leq \frac{x^\varepsilon RQ^2}{(q_1,q_2)M}}\!\!\!\!\!\!\!\!\!\!\!\! M \varphi_H(h) \!\!\!\!\!\!\sum_{\substack{n \\ (n,q_1r)=1 \\ (n+\ell r,q_2)=1}} \!\!\!C^*(n) \beta(n) \overline{\beta(n+\ell r)} e_{[q_1,q_2]r}(y h).
\end{align*} (Here the modulus of $\varphi_H$ is $[q_1,q_2]r$.) 
The component of $\Sigma_0^*(b_1,b_2)$ which has $(q_1,q_2)=1$ is given by 
\begin{align*}
X=\left(\sum_m \psi_M(m) \right)\!\!\sum_{\substack{r \asymp R \\ (r,a)=1}} \sum_{\substack{ |\ell| \ll \frac{N}{R} \\ \ell \neq 0}}  \mathop{\sum\sum}\limits _{\substack{q_1, q_2 \\  q_1r, q_2r \mid P(x^\delta) \\ (q_1q_2,a)=1\\ (q_1q_2,P(D_0))=1\\ (q_1,q_2)=1}}\!\!\!\!\!\! \dfrac{c_{q_1,r} \overline{c_{q_2,r}}}{r[q_1,q_2]} \!\!\!\sum_{\substack{n \\ (n,q_1r)=1 \\ (n+\ell r,q_2)=1}}  \beta(n) \overline{\beta(n+\ell r)}.
\end{align*}
(Here we used that $C^*(n)=1$ on this component.) Observe that $X$ does not depend on $b_1$ and $b_2$.

\vspace{3mm} 
Writing $(q_1,q_2)=q_0$ and using Lemma~\ref{lem:polymathconvenienterrors}, the remaining terms of $\Sigma_0^*(b_1,b_2)$, for which $q_0\neq 1$, are bounded in absolute value by 
\begin{align*}
\dfrac{M \log(x)^{O(1)}}{R}\sum_{r \asymp R} \sum_{\substack{|\ell| \ll \frac{N}{R} \\ \ell \neq 0}} \sum_{\substack{1 \neq q_0 \ll Q \\ (q_0, P(D_0))=1 }} \mathop{\sum\sum}\limits_{\substack{q_1,q_2 \asymp Q \\ q_1r, q_2r \mid P(x^\delta) \\ (q_1,q_2)=q_0}} \dfrac{1}{[q_1,q_2]} \!\!\!\sum_{\substack{n \\ (n,q_1r)=1 \\ (n+\ell r,q_2)=1}} \!\!\!\!\!C^*(n)(\tau(n)\tau(n+\ell r))^{O(1)} \ll_A \dfrac{N^2M}{R\log(x)^{A}}.
\end{align*}
Now we look at $\Sigma_1^*(b_1,b_2)$. We split the sum up according to the value of $(q_1,q_2)$ and deduce that
\begin{align}
|\Sigma_1^*(b_1,b_2)| \ll \!\!\! \sum_{\substack{ |\ell| \ll \frac{N}{R} \\ \ell \neq 0}} \sum_{q_0 \ll Q}\dfrac{q_0 M}{RQ^2}  \sum_{r \asymp R}\!\!\!\!\!\!\mathop{\sum\sum}\limits _{\substack{q_1, q_2 \asymp Q/q_0 \\ (q_1,q_2)=1 \\  (rq_0q_1q_2,ab_1b_2)=1 \\  q_0q_1r, q_0q_2r \mid P(x^\delta) }} \!\!\!\!\!\!\Bigg|\sum_{\substack{ |h| \leq H(q_0) \\ h \neq 0}}\!\!\!\!\! \varphi_H(h)\!\!\!\! \!\!\!\!\sum_{\substack{n \\ (n,rq_0q_1)=1 \\ (n+\ell r,q_2)=1}}\!\!\!\!\! \!\!\!C_0(n) \beta(n) \overline{\beta(n+\ell r)} \Phi_\ell(h,n,r,q_0,q_1,q_2) \Bigg|. \label{oneofthelastequations}
\end{align}
(Notice that we introduced the condition $(rq_0q_1q_2,b_1b_2)=1$. Since $(b_1b_2,P(x^\delta))=1$, while $rq_0[q_1,q_2]$ divides $P(x^\delta)$,  this condition is  trivially satisfied and does not change the summation range. This is merely a
 technical restriction, which now allows us to drop the requirement  $(b_1b_2,P(x^\delta))=1$.)

\vspace{3mm}
Note that the inner sums of (\ref{oneofthelastequations}) are of the form $\Sigma^*$ described in the statement of  Lemma~\ref{lem:amendmenttheorem5.8}, and recall that we assumed that  $\Sigma^*(z_*) \ll_{ \varepsilon} (q_0,\ell) RQ^2N(x^\varepsilon q_0^2)^{-1} $ for any $z_* \in Z_0^*$.  Hence,
\begin{align*}
|\Sigma_1^*(b_1,b_2)| \ll_{ \varepsilon}  \sum_{q_0 \ll Q}\dfrac{q_0 M}{RQ^2}\sum_{\substack{ |\ell| \ll \frac{N}{R} \\ \ell \neq 0}}\dfrac{(q_0,\ell)RQ^2N}{x^{\varepsilon}q_0^{2}}\ll_{A, \varepsilon} \dfrac{N^2M}{R\log(x)^A}.
\end{align*}
Overall we have shown that $\Sigma(b_1,b_2) = X +O_{A, \varepsilon}(N^2MR^{-1}\log(x)^{-A})$, with $X$ not dependent on $b_1$ and $b_2$. By Lemma~\ref{lem:PolymathSection5} this implies that (\ref{finalresultappendix}) is satisfied.
\end{proof}
 
For the final step of our proof of Lemma~\ref{lem:summaryofpolymath}, we use the Cauchy-Schwarz inequality to remove $\beta$. We  proceed like in Section~8A of~\cite{Polymath:2014:EDZ}.  

\begin{lemma}  \label{lem:proofoflemma3-finalstep}
Let $\omega, \delta >0$ and $\varepsilon  \in (0, 10^{-100} \delta)$. Let $\alpha$ and $\beta$ be coefficient sequences at scales $M$ and $N$ with $x \ll M(x)N(x) \ll x$ and $N(x)= x^{\gamma(x)}$, where $\gamma(x) \in (4\omega+2\delta,\frac{1}{2}-2\omega-8\varepsilon)$.    Assume $\beta$ has the Siegel-Walfisz property. 

\vspace{3mm}
Let $Z_1$ and $\Sigma_1$ be as described in Lemma~{\rm\ref{lem:summaryofpolymath}} and let $\Sigma^*$ be as described in Lemma~{\rm\ref{lem:amendmenttheorem5.8}}. Suppose that $$\Sigma_1(z;x) = O_{\varepsilon}((q_0,\ell)RQN UV^2x^{-4\varepsilon}) $$ for  $x>1$ and $z \in Z_1(x)$.  Then also $\Sigma^*(z;x) \ll_{ \varepsilon} (q_0,\ell)RQ^2 N(x^\varepsilon q_0^2)^{-1}$ for   $x>1$ and $z \in Z_0^*(x)$.  
\end{lemma}

\begin{proof}
We consider $\Sigma^*$, described in Lemma~\ref{lem:amendmenttheorem5.8}. We wish to show that $\Sigma^* \ll_\varepsilon \frac{(q_0,\ell)RQ^2 N}{x^\varepsilon q_0^2}$. Recall that
\begin{align*}
&\Sigma^*= \sum_{\substack{r \asymp R \\ (r,ab_1b_2)=1}} \mathop{\sum\sum}\limits_{\substack{q_1,q_2 \asymp Q/q_0\\ (q_1,q_2)=1 \\ (q_1q_2,ab_1b_2)=1 \\  q_0q_1r, q_0q_2r \mid P(x^\delta)
}} \Bigg| \sum_{0<|h| \leq H} \!\!\! \varphi_H(h)\sum_{\substack{n \\ (n,rq_0q_1 )=1 \\ (n+\ell r, q_0q_2)=1 }} C_0(n)\beta(n)\overline{\beta(n+\ell r)}  \Phi_\ell(h,n,r,q_0,q_1,q_2) \Bigg|,
\end{align*}
where $H = \frac{x^\varepsilon RQ^2}{q_0 M}$. (Since the value of $q_0$ has been fixed in $\Sigma^*$, it is no longer necessary to indicate the dependence of $H$ on $q_0$.) 
Observe first that if $H < 1$, then the sum over $h$ is empty and $\Sigma^* = 0$. Hence we may assume that $H \geq  1$.  But then $q_0^{-1}x^{-5\varepsilon}Q/H \leq x^{-5\varepsilon}(Q/q_0)$. Further, since $RQ \ll x^{1/2+2\omega+\varepsilon}$, $MN \asymp x$ and $N=x^\gamma$ with $\gamma \leq 1/2-2\omega-8\varepsilon$, we also have $q_0^{-1}x^{-\delta-5\varepsilon}Q/H \gg x^{-\delta+\varepsilon} $. Hence for any $q_1 \asymp Q/q_0$ with $q_1 \mid P(x^\delta)$ we can find $u_1$ and $v_1$ with $q_1=u_1v_1$ and $q_0^{-1}x^{-\delta-5\varepsilon}Q/H \leq u_1 \leq q_0^{-1}x^{-5\varepsilon}Q/H$. 

\vspace{3mm}
 So, using dyadic decomposition, we deduce that the bound $\Sigma^* \ll_\varepsilon \frac{(q_0,\ell)RQ^2 N}{x^\varepsilon q_0^2}$ holds if the following is true for every $H^* \in \mathbb{R} \setminus \{0\}$ with $1 \ll |H^*| \ll H$ and all  $U$ and $V$ which satisfy  (\ref{def:UV1}), (\ref{def:UV2}) and (\ref{def:UV3}):
\begin{align*} 
\sum_{\substack{r \asymp R }} \!\!\sum_{\substack{u_1, v_1 \\ u_1 \asymp U, v_1 \asymp V \\ q_0u_1v_1r  \mid P(x^\delta) }} \!\!\!\!\!\sum_{\substack{q_2 \asymp Q/q_0\\ (u_1v_1,q_2)=1 \\ q_0q_2r  \mid P(x^\delta)\\ (rq_0u_1v_1q_2, ab_1b_2)=1}} \!\!\!\!\!\!\!\!\!\!\!\!\ \Bigg|\sum_{h \sim H^*}\!\! \varphi^{\star}_H(h)\!\!\!\!\!\!\!\!\!\!\!\!\sum_{\substack{n \\ (n,rq_0u_1v_1 )=1 \\ (n+\ell r, q_0q_2)=1 }}\!\!\!\!\!\!\!\!\!\! C_0(n)\beta(n)\overline{\beta(n+\ell r)}  \Phi_\ell(h,n,r,q_0,u_1v_1,q_2) \Bigg|\ll_\varepsilon \dfrac{(q_0,\ell)RQ^2N}{x^{2\varepsilon}q_0^{2}}. 
 \end{align*}
(Recall here that $\varphi_H(h)$ has modulus $rq_0q_1q_2$, $\varphi_H^{\star}(h)$ has modulus $rq_0u_1v_1q_2$ and $\varphi_H^{\star\star}(h)$ has modulus  $rq_0u_1v_2q_2$.)  
 Denote the LHS of the above expression by $\Sigma_1^*$. We may write 
 \begin{align*}
\Sigma_1^*&=  \sum_{\substack{r \asymp R  }}\sum_{\substack{u_1, q_2 \\ u_1 \asymp U  \\ q_2 \asymp Q/q_0 }}  \sum_{\substack{n \\ (n,rq_0u_1)=1 \\ (n+\ell r, q_0q_2)=1 }} C_0(n)\beta(n)\overline{\beta(n+\ell r)}   \sum_{h \sim H^*} \sum_{ \substack{v_1 \asymp V \\ (v_1,n)=1} }   c_{r,u_1,v_1,q_2} \varphi^{\star}_H(h)  \Phi_\ell(h,n,r,q_0,u_1v_1,q_2), 
\end{align*}
where $c_{r,u_1,v_1,q_2}$ are $1$-bounded constants, supported on $r$, $u_1$, $v_1$ and $q_2$ with  $q_0u_1v_1 q_2r \mid P(x^\delta)$ and $(rq_0u_1v_1q_2, ab_1b_2)=1$. 
Next we apply the Cauchy-Schwarz inequality. We find  $ (\Sigma_1^*)^2 \leq \Upsilon_1 \Upsilon_2$, where 
\begin{align*}
\Upsilon_1 &= \sum_{r \asymp R}\sum_{\substack{ u_1 \asymp U  }}\sum_{\substack{  q_2 \asymp Q/q_0 }}  \sum_{\substack{n \\ (n,rq_0u_1 )=1 \\ (n+\ell r, q_0q_2)=1 }} C_0(n)|\beta(n)\overline{\beta(n+\ell r)}|^2  \ll \dfrac{(q_0,\ell)RQUN}{q_0^2}, \\
\Upsilon_2 &= \sum_{r \asymp R}\sum_{\substack{ u_1 \asymp U  }}\sum_{\substack{  q_2 \asymp Q/q_0 }}  \sum_{\substack{n \\ (n,rq_0u_1 )=1 \\ (n+\ell r, q_0q_2)=1 }}  C_0(n)\psi_N(n)\Bigg|\sum_{h \sim H^*} \sum_{\substack{v_1 \asymp V \\ (v_1,n)=1}  }   c_{r,u_1,v_1,q_2} \varphi^{\star}_H(h)  \Phi_\ell(h,n,r,q_0,u_1v_1,q_2) \Bigg|^2\\
&\leq \! \!\!\!\!\!\!\!\!\!\!\!\!\!\!\!\! \sum_{\substack{r,u_1,v_1,v_2,q_2 \\ r \asymp R \\ u_1 \asymp U\\ v_1, v_2 \asymp V  \\ q_2 \asymp Q/q_0 \\ (rq_0u_1v_1v_2q_2, ab_1b_2)=1 \\
q_0u_1v_1q_2r  \mid P(x^\delta)\\
q_0u_1v_2q_2r  \mid P(x^\delta)}}\!\!\!\! \!\! \!\!\!\!\Big| \sum_{\substack{h_1, h_2 \\ h_1 \sim H^* \\ h_2 \sim H^* }} \sum_{\substack{n \\ (n,rq_0 )=1 \\ (n,u_1v_1v_2 )=1 \\ (n+\ell r, q_0q_2)=1 }}  \!\!\!\!\!\!\!\! \! C_0(n) \varphi^{\star}_H(h_1) \overline{\varphi^{\star\star}_H(h_2)}\psi_N(n) \Phi_\ell(h_1,n,r,q_0,u_1v_1,q_2) \overline{\Phi_\ell(h_2,n,r,q_0,u_1v_2,q_2)} \Big|
\end{align*}
for some smooth  coefficient sequence $\psi_N$, located at scale $N$, with $\psi_N(n) \geq 0$ for all $n$ and $\psi_N(n) \geq 1$ on the support of $\beta$. Then $\Sigma_1^* \ll_{ \varepsilon} (q_0,\ell)RQ^2N(x^{2\varepsilon}q_0^2)^{-1}$ provided that $\Upsilon_2 \ll_{ \varepsilon} (q_0,\ell)RQN UV^2x^{-4\varepsilon}$. However, looking at the notation of Lemma~\ref{lem:summaryofpolymath}, we see that $\Upsilon_2 \ll  \Sigma_1(z)$ for a suitable choice of $z \in Z_1$. Recall the assumption  $\Sigma_1(z) = O_{\varepsilon}((q_0,\ell)RQN UV^2x^{-4\varepsilon}) $ to conclude the proof.
\end{proof} 

Lemma~\ref{lem:summaryofpolymath} is simply a summary of Lemma~\ref{lem:PolymathSection5}, Lemma~\ref{lem:amendmenttheorem5.8} and Lemma~\ref{lem:proofoflemma3-finalstep}.

\section*{Acknowledgements}

I would like to thank my supervisor, James Maynard, for  many insightful conversations and very helpful comments and suggestions. Further, I would like to thank Andrew Sutherland for computing short admissible $k$-tuples.  My work was supported by an EPSRC studentship. 

\bibliography{bibstad}
\bibliographystyle{acm} 

\end{document}